\documentclass[a4paper,11pt]{article}

\usepackage[utf8]{inputenc} 

\usepackage{amssymb} 
\usepackage{amsmath} 
\usepackage{amsthm} 
\usepackage{graphicx}
\usepackage{epsfig,color}
\usepackage{enumerate,subfigure}
\usepackage{rotating}
\usepackage{multirow}
\usepackage[table]{xcolor}
\usepackage{boldline,natbib}
 \usepackage{arydshln}

\definecolor{k4}{rgb}{0.8,0.8,0.8}
\definecolor{k3}{rgb}{0.6,0.6,0.6}
\definecolor{k2}{rgb}{0.4,0.4,0.4}
\definecolor{k1}{rgb}{0.2,0.2,0.2}

\newcommand{\fb}{\mathbf{f}}
\newcommand{\Fb}{\mathbf{F}}
\newcommand{\mub}{\boldsymbol{\mu}}
\newcommand{\betab}{\boldsymbol{\beta}}
\newcommand{\alphab}{\boldsymbol{\alpha}}
\newcommand{\sigmab}{\boldsymbol{\sigma}}
\newcommand{\zb}{\mathbf{z}}
\newcommand{\xb}{\mathbf{x}}
\newcommand{\tb}{\mathbf{t}}
\newcommand{\cb}{\mathbf{c}}
\newcommand{\mb}{\mathbf{m}}
\newcommand{\pb}{\mathbf{p}}

\newcommand{\lambdab}{\boldsymbol{\lambda}}
\newcommand{\rhob}{\boldsymbol{\rho}}

\newcommand{\R}{\mathbb{R}}

\newcommand{\eps}{\varepsilon}
\newcommand{\f}{\dfrac}

\newcommand{\G}{\mathbf{G}}

\newcommand{\Zc}{Z_\mathrm{sel}}
\newcommand{\Zv}{Z_\mathrm{div}}
\renewcommand{\epsilon}{\varepsilon}

\newcommand{\columnname}[1]
{\makebox[\tempwidth][c]{\textbf{#1}}}

\renewcommand{\leq}{\leqslant}
\renewcommand{\geq}{\geqslant}

\renewcommand{\eps}{\varepsilon}

\renewcommand{\hat}{\widehat}

\newcommand{\m}{m}
\newcommand{\Var}{\mathrm{Var}}
\newcommand{\bz}{\mathbf{z}}

\newtheorem{thm}{Theorem} 
\newtheorem{remark}[thm]{Remark} 

\numberwithin{equation}{section}

\usepackage{hyperref}
  \definecolor{gr}{RGB}   {18,   154,   47 } 
  \definecolor{orange}{RGB} { 255, 127, 0 }
  \definecolor{bl}{rgb}   {0.,   0.,   1. } 
  \definecolor{mg}{rgb}   {0.5,  0.,    0.7} 

 \hypersetup{
     colorlinks=true,       
     linkcolor=red,          
     citecolor=orange,        
     filecolor=mg,      
     urlcolor=cyan           
 }

\def\ds{\displaystyle}
\usepackage{epstopdf}
\usepackage{geometry}
\usepackage[left]{lineno}
\usepackage{setspace}
\title{
Adaptation to a changing environment: what me Normal? 
}
\author{J. Garnier$^{1,*}$, O. Cotto$^2$, T. Bourgeron$^3$, E. Bouin$^4$,  T. Lepoutre$^{5,6}$, O. Ronce$^{7,8}$ and V. Calvez$^{5,6}$\\
\small{$^1$LAMA, UMR 5127, CNRS, Univ. Grenoble Alpes, Univ. Savoie Mont Blanc, Chambéry, France}\\
\small{$^2$ PHIM Plant Health Institute, INRAE, Univ Montpellier, CIRAD, Institut Agro, IRD, Montpellier, France} \\
\small{$^3$ ADIA, Abu Dhabi, United Arab Emirates}
\\
\small{$^4$ CEREMADE, UMR 7534, CNRS, Univ. Paris Dauphine, Paris, France}\\
\small{$^5$ ICJ, UMR 5208, CNRS, Univ. Claude Bernard Lyon 1, Lyon, France}\\
\small{$^6$Equipe-projet Inria Dracula, Lyon, France}\\
\small{$^7$ISEM, Univ Montpellier, CNRS, IRD, Montpellier, France}\\
\small{$^8$CNRS, Biodiversity Research Center, Univ. British Columbia, Vancouver, BC, Canada}
}
\date{}

\begin{document}

\maketitle



 

\noindent{\it $^*$ Corresponding author: {\tt jimmy.garnier@univ-smb.fr}, LAMA, UMR 5127, Univ. Savoie Mont Blanc,  Bâtiment Le Chablais, Campus Scientifique, 73376 Le Bourget du Lac, France}

\

\noindent\emph{Declaration of interest: none}

\

\noindent\emph{Credit author statement: 
VC, OC and OR originally formulated the project; TB, VC, and JG mathematically analysed the model and performed the numerical simulations, with specific contributions from EB and TL; VC, JG, OC and OR wrote the manuscript.
} 

\

\noindent\emph{Fundings: The author(s) acknowledge support of the Institut Henri Poincaré (UAR 839 CNRS-Sorbonne Université), and LabEx CARMIN (ANR-10-LABX-59-01). JG acknowledges GLOBNETS project (ANR-16-CE02-0009) and ModEcoEvo project funded by the Univ. Savoie Mont-Blanc. This project has received funding from the European Research Council 
(ERC) under the European Union’s Horizon 2020 research and innovation 
programme (grant agreement No 865711). OR acknowledges support from the Peter Wall Institute of Advanced Studies and from the France Canada Research Funds.}

\


\abstract{ 
{
Predicting the adaptation of populations to a changing environment is crucial to assess the impact of human activities on biodiversity. Many theoretical studies have tackled this issue by modeling the evolution of quantitative traits subject to stabilizing selection around an optimum phenotype, whose value is shifted continuously through time. In this context, the population fate results from the equilibrium distribution of the trait, relative to the moving optimum. Such a distribution may vary with the shape of selection, the system of reproduction,  the number of loci, the mutation kernel or their interactions. Here, we develop a methodology that provides quantitative measures of population maladaptation and potential of survival directly from the entire profile of the phenotypic distribution, without any a priori on its shape.
We investigate two different models of reproduction (asexual and infinitesimal sexual models of inheritance), with general forms of selection.
{ In particular, we recover that fitness functions such that selection weakens away from the optimum lead to evolutionary tipping points, 
with an abrupt collapse of the population when the speed of environmental change is too high. 
Our unified framework furthermore allows highlighting the underlying mechanisms that lead to this phenomenon.} More generally, it allows discussing similarities and discrepancies between the two reproduction models, the latter being ultimately explained by different constraints on the evolution of the phenotypic variance. We demonstrate that the mean fitness in the population crucially depends on the shape of the selection function in the sexual infinitesimal model, in contrast with the asexual model.
In the asexual model, we also investigate the effect of the mutation kernel and we show that kernels with higher kurtosis tend to reduce maladaptation and improve fitness, especially in fast changing environments.
}

}


\section{Introduction}\label{sec:1}
Rapid environmental changes resulting from human activities have motivated the development of a theory to understand and predict the corresponding response of populations. Efforts have specially been focused on identifying conditions that allow populations to adapt and survive in changing environments ~\cite[e.g.][for  pioneering work]{LynGabWoo91,LynLan93,BurLyn95}. To this aim, most theoretical studies have modeled the evolution of polygenic quantitative traits subject to stabilizing selection around some optimal phenotype, whose value is shifted continuously through time ~\citep[see][]{KopMat14,Wal12,Ale14}.
A major prediction of these early models is that when the optimal phenotype changes linearly with time, it will be tracked by the mean phenotype in the population with 
{ a lag that eventually stabilizes over time. This evolutionary lag, which quantifies the maladaptation induced by the environmental change,
is predicted to depend on the rate of the change, on the genetic standing variance for the phenotypic trait and on the strength of stabilizing selection on the trait. The maladaptation of the population
due to the environmental change, 
also decreases the mean fitness of the population, which is commonly defined as the \emph{lag load} or\emph{evolutionary load} ~\citep{LynLan93,LanSha96}.}
Thus, above a critical rate of change of the optimal phenotype with time, the evolutionary lag  is so large that { the lag-load of the population will rise above}
the value that allows its persistence and the population will be doomed to extinction. 

{ These predictions have typically been derived under the assumptions of (i) a particular form of selection, (ii) a constant genetic variance for the evolving trait, (iii) a Gaussian distribution of phenotypes and breeding values in the population.  The selection function, describing how the Malthusian fitness declines away from the optimum, has typically a quadratic shape in many models}~\citep{Bur99,KopMat14}.
{ However, the shape of selection functions is difficult to estimate and some studies suggest that it strongly deviates from a quadratic shape in the case of phenological traits involved in climate adaptation}~\citep{Gau20}.
{ Moreover, the dynamics of adaptation} become more complex when the shape of selection deviates from the quadratic form and few theoretical studies have addressed this issue. Recently, \citep{OsmKla17,Kla20} have shown that ``evolutionary tipping points'' occur when the strength of selection weakens away from the optimum. In this situation, the population abruptly collapses when the speed of environmental change is too large. { In this paper, we  aim to investigate in a general setting the effects of the shape of selection functions on the adaptation of the population under environmental changes.}


The genetic standing variance also plays a key role in the adaptation to changing environments and the determination of the critical rate of change. In many quantitative genetic models, this variance is assumed to be constant. Although it is approximately true at a short time scale, over a longer time scale the variance in the population is also subject to evolutionary change. More generally, obtaining mathematical predictions for the dynamics and the equilibrium value of the variance
remains a notoriously difficult issue for many theoretical population genetics models~\citep{BarTur89,Bur00,BarKei02,JohBar05,Hil10}. How the genetic variance evolves in a changing environment has therefore been explored mostly through simulations~\citep{Jon12, Bur99, WaxPec99}.
In our paper,  we overcome this problem by modeling the evolution of the entire phenotype distribution,
{  which in particular provides insights on the effect of maladaptation, induced by environmental changes, on the evolution of genetic standing variance.

}

{ Many theoretical works assumed that the phenotype distribution is Gaussian}~\citep{LynGabWoo91}.
In the absence of environmental change, there are indeed many circumstances where the phenotypic distribution in the population is well captured by Gaussian distributions in  quantitative genetics models.
For example, in asexual populations, the distribution of a polygenic trait is Gaussian at mutation-selection equilibrium,
providing that mutation effects are weak and selection is { quadratic}~\citep{Kim65,Lan75,Fle79}.
{ In the case of sexual reproduction, similar outcomes are expected with the celebrated Fisher infinitesimal model of inheritance introduced by~\citet{Fis18}. In this model, quantitative traits are under the control of many additive loci and each allele has a relatively small contribution on the character~\citep{Fis18}. Within this framework, offspring are normally distributed  within families around the mean of the two parental trait values, with fixed standing variance~\citep[and references therein]{TurBar94,Tur17,BarEthVeb17}.}
As a result, the phenotype distribution of the full population is a Gaussian under various assumptions on the selection function (see \citealt{TurBar94} under truncation selection, or see \citep{Raoul21} and \citep{CalGarPat19} for a wider class of selection functions).
{In the process of adaptation to environmental change, since the mean phenotype is lagging behind the optimum,  selection however may induce a skew in the distribution~\citep{Jon12}.} The distribution of the mutational effects can also have a strong influence on the distribution as well, in particular when the evolutionary lag is large~\citep{WaxPec99}. 
{ The Gaussian approximation of the phenotypic distribution should therefore naturally be questioned for both model of inheritance (asexual and sexual infinitesimal model).}


The main objective of our work is to derive signatures of maladaptation at equilibrium, {\em e.g.} the evolutionary lag, the mean fitness and the genetic standing variance, depending on some general shape of selection and features of trait inheritance.
{ Those three components are linked by two generic identities describing the demographic equilibrium and  the genetic equilibrium. Would the genetic variance be known, it would be possible to identify both the evolutionary lag and  the mean fitness~\citep{KopMat14}. In the general case, a third relationship is, however, needed. }
To this aim, we shall compute accurate approximations of the phenotypic distribution. 
Several methodological alternatives have been developed to unravel the phenotypic distribution, without any {\em a priori} on its shape. First methods attempted to derive the equations describing the dynamics of the mean, the variance and the higher moments of the distribution~\citep{Lan75,BarTur87,TurBar90,FraSla90}.
Then in his pioneering work,~\citet{Bur91} derived relationships between the cumulants of the distribution, which are functions of the moments. However this system of equations is not closed, as the cumulants influence each other in cascade.
More recently, \citet{MarRoq16}  analyzed a large class of  integro-differential models where the trait coincides with the fitness, through the partial differential equation (PDE) satisfied by the cumulant generating function (CGF). They applied their approach to the adaptation of asexual populations facing environmental change, using the Fisher Geometric Model for selection and specific assumptions on trait inheritance (diffusion approximation for the mutational effects)~\citep{RoqPat20}. However, the extension of their method to different models of selection or trait inheritance (general mutational kernel) seems difficult mainly because it relies on specific algebraic identities to reduce the complexity of the problem.

Here, we use quantitative genetics models based on integro-differential equations to handle various shapes of stabilizing selection, and trait inheritance mechanisms. {While we deal with a large class of mutational kernels (including thin- and fat-tailed kernels) in the asexual model, we consider the Fisher infinitesimal model as a mechanism of trait inheritance in sexually reproducing populations.} We assume that the environment is changing linearly with time, as in the classical studies reviewed in \citep{KopMat14}. 
In order to provide quantitative results, we assume that very little variance in fitness is introduced in the population through either mutation or recombination events during reproduction. It allows some flexibility about the trait inheritance process as well as the shape of the selection function.
Under this assumption of \emph{small variance regime}, a recent mathematical methodology had been developed to derive analytical features in models of quantitative genetics in asexual populations in fixed phenotypic environment~\citep{DieJabMisPer05,PerBar08,LorMirPer11,MirRoq16,Mirrahimi-2017,Calvez-Lam-2020}. This asymptotic method  was first introduced by~\citet{DieJabMisPer05} and~\citet{Per07b} in the context of evolutionary biology as an alternative formulation of adaptive dynamics, when the mutational effects are supposed to be small, but relatively frequent.
{Recently, this methodology has been  also applied to the infinitesimal model for sexual reproduction in a stationary fitness landscape~\citep{CalGarPat19,Pat20}. In the present paper, we apply this methodology to the case of a moving optimum. In this context, the extension to the infinitesimal model for sexual reproduction is new to the best of our knowledge.}

From a mathematical perspective, the assumption of small variance regime is analogous to some asymptotic analysis performed in mathematical physics, such as the approximation of geometric optics for the wave equation at high frequency~\citep{evans10,rauch2012hyperbolic},  semi--classical analysis for the Schr\"odinger equation in quantum mechanics~\citep{DisSjor99b,zworski2012semiclassical}, and also the large deviation principle for stochastic processes \citep{Fle77,EvaIsh85,FreWen12}.

Conversely to previous methods focusing on the moments of the phenotypic distribution, our approach focuses on the entire phenotypic distribution and it provides an accurate approximation of the phenotypic distribution even if it deviates significantly from the Gaussian distribution. As a result, our method allows deriving analytical formulas for biologically relevant quantities, that is for instance the evolutionary lag measuring maladaptation, the genetic standing variance of the population, the { lag-load depressing the population mean fitness and critical rates of environmental changes in a gradually changing environment, without solving the complete profile of the distribution and for a general class of selection functions and reproduction models. We are consequently able to answer the following issues
\begin{itemize}
    \item What is the effect of the shape of selection on the adaptation of a population to continuously changing environment?
    \item How does the distribution of mutational effects affect the adaptation dynamics?
    \item Does the choice of a particular reproduction model influence predictions about the adaptation dynamics of a population?
\end{itemize}
}



\section{Models and methodology}\label{sec:model}
\begin{figure}[h]
\begin{center}
\includegraphics[width = 0.85\linewidth]{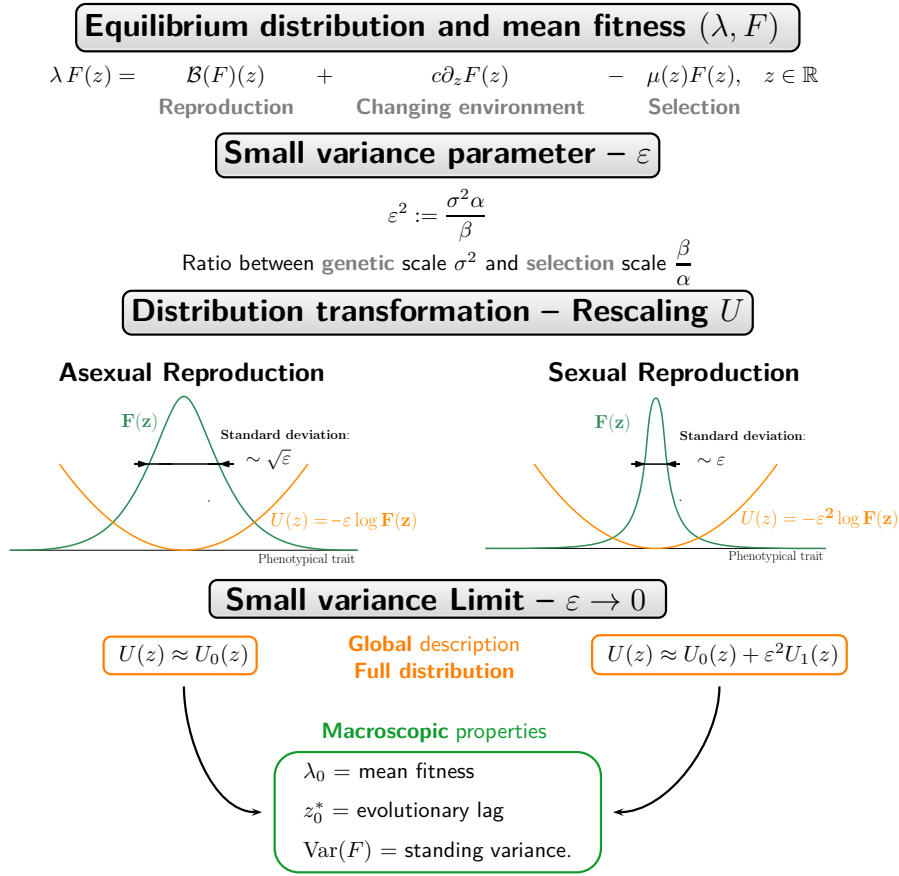} 
\caption{Schematic description of our methodology. To describe the equilibrium $F$ we need the following steps: (1) Identify the scaling parameter $\epsilon$ and rescale the equation satisfied by the distribution $F$; (2) Transform the distribution $F$ into $U.$ The transformed distribution $U$ is the logarithmic of the  density $F$, normalized by the ratio $\epsilon$ in the asexual reproduction case and by $\epsilon^2$ in the infinitesimal sexual reproduction case; (3) Identify the limit equation for $U$ as $\epsilon\to0$ (orange boxes) and deduce macroscopic properties (green box) such as the { mean fitness} $\lambda_0$, the evolutionary lag $z^*_0$ in the population or the phenotypic variance at equilibrium $\Var(F)$.}
\label{fig:Diag_methodo}
\end{center}
\end{figure}
First, we describe in detail our general model of mutation-selection under changing environment with two different reproduction models (asexual and infinitesimal sexual) (Section~\ref{subsec:model}). Then, we introduce the rescaled model including the variance parameter $\epsilon$ (Section~\ref{sec:adim}) and we describe our methodology to investigate the regime of small variance (see Figure~\ref{fig:Diag_methodo} for a sketch of the methodology). It is based on the asymptotic analysis with respect to this small parameter (Section~\ref{sec:small_varinace}). In Section~\ref{sec:results_adim}, we provide, in the regime of small variance, analytical formula for the different characteristic quantities of the evolving population  --- mean fitness, evolutionary lag or standing variance --- for the two different reproduction models: asexual model (Section~\ref{sec:asex}) and infinitesimal sexual model (Section~\ref{sec:sex}).  After scaling back our results in the original units, we can compare the two reproduction regimes and discuss the effect of changing environment on different characteristics of the population: the lag (Section~\ref{sec:lag}), the mean fitness (Section~\ref{sec:fitness}) and the standing variance (Section~\ref{sec:variance}). Furthermore, we discuss the persistence of the population according to the speed of the changing environment (Section~\ref{sec:persitence}) and we describe numerically the whole distribution of the population (Section~\ref{sec:asex_num}).

\subsection{The general model under changing environment}\label{subsec:model}
We consider a continuously growing population with overlapping generations and density dependence. The population is structured by a one--dimensional phenotypic trait, denoted by $\xb\in \R$. The density of individuals with trait $\xb$ is  $\fb(\tb,\xb)$  at time $\tb>0$. For the sake of simplicity, the birth rate $\betab(\xb) = \betab$ is a positive constant. Selection acts through the intrinsic mortality rate $\mub(\xb),$ by means of stabilizing selection around some optimal value. In order to capture the dynamics of the population under a gradual environmental change, we assume that the optimal trait $\xb=0$ is shifted at a constant speed $\cb>0$. We define the phenotypic lag  as the difference between the phenotypic value $\xb$ and the optimal value at time $\tb$: $\zb=\xb-\cb\tb$. It quantifies the maladaptation of an individual of trait $\xb$ in the changing environment. 
The intrinsic mortality rate $\mub$ is decomposed as follows 
\begin{equation}\label{eq:mu}
\mub(\zb) = \mub_0 + \mb(\zb)\,,
\end{equation}
where $\mub_0$ is the basal mortality rate at the optimum at low density. We assume that $\betab>\mub_0$ to ensure that the population at the optimum will not go extinct in the absence of environmental change. The function $\mb(\zb)=\mb(\xb-\cb\tb)$ is the increment of mortality due to maladaptation. 
The function $\mb\geq 0$ attains its unique minimum value at $\zb = 0$ where $\mb(0)=0$, and it is increasing with respect to $|\zb|$: $\mb$ is decreasing on $(-\infty,0)$ and increasing on $(0,\infty)$.
The strength of stabilizing selection is captured by the positive parameter 
\begin{equation}\label{eq:scale_alpha}
\alphab =  \mub''(0) =  \mb''(0)>0\, .
\end{equation}

The dynamics of the density $\fb(\tb,\xb)$ is given by the following equation:
\begin{equation} \label{eq:f z} 
 \partial_{\tb} \fb(\tb,\xb) + \Big(\mub(\xb - \cb\tb) + (\betab-\mub_0)\rhob(\tb) \Big)  \fb(\tb,\xb) = \betab \mathcal{B}( \fb(\tb,\cdot) )(\xb)\, , 
\end{equation}
where the term $\rhob(\tb) = \int_\R \fb(\tb,\xb')d\xb'$ corresponds to the size of the population. 
This nonlinear term introduces density--dependent mortality in the model. 

The operator $\mathcal{B}$ describes how new individuals with phenotype $\xb$ are generated depending on the whole phenotypic density. For simplicity, we assume no environmental effects on the expression of the phenotype and phenotypic values equal to breeding values.
We consider the two following  scenarios for the reproduction operator $\mathcal{B}$:

\paragraph*{Asexual genetic model of reproduction with mutations.}
We first consider the case of asexual reproduction where the phenotype of an offspring $\xb$ is drawn randomly around the phenotype of its single parent $\xb'$. The mutation kernel $K$ describes the distribution of the effects of mutations.   The reproduction operator has then the following expression: 
\begin{equation}\label{eq:asexual B}
\mathcal{B}(\Fb)(\xb) =  \dfrac1{\sigmab} \int_\R  K\left(\dfrac{\xb-\xb'}\sigmab\right)\Fb(\xb')\, d\xb' \, ,
\end{equation}
where $\sigmab^2$ is the mutational variance and $K$ is a symmetric normalized probability density function.
We furthermore assume that $K$ decays faster than some exponential function for large $|y|$. 
This is usually called  a {\em thin--tailed kernel}. 
This  corresponds to the scenario where the mutations with large effect on phenotypic traits are rare. 

The extremal case corresponding to accumulation of infinitesimal effects of mutations is referred to as {\em diffusion approximation}. This translates into the following formula
\begin{equation}\label{eq:asexual B diffusion}
\mathcal{B}(\Fb)(\xb) = \Fb(\xb) +  \dfrac{\sigmab^2}2 \partial_\xb^2\Fb(\xb)\, ,
\end{equation}
In this case, the shape of the mutation kernel does not matter and only the variance remains. 

In the absence of selection, the standing variance of the phenotypic distribution thus increases indefinitely in this model. This asexual model does not impose any strong constraint on the standing variance of the phenotypic distribution of the population, contrary to the next case we consider (see below).

\paragraph*{Infinitesimal model of sexual reproduction.}
Secondly, we consider the case where the phenotype of the offspring $\xb$ is drawn randomly around the mean trait of its parents $(\xb_1,\xb_2)$, following a Gaussian distribution $G_{\sigmab^2}$ with a given variance $\sigmab^2/2$. This is known as the Fisher infinitesimal model 
\citep{Fis18,Bul80,TurBar94,Tuf00,BarEthVeb17}.
The reproduction operator has then the following expression:
\begin{equation}\label{eq:sexual B}
\mathcal{B}(\Fb)(\xb) = \iint_{\mathbb{R}^2} G_{\sigmab^2}\left(\xb - \dfrac{\xb_1 + \xb_2}2 \right)\Fb(\xb_1)\left (\dfrac{\Fb(\xb_2)}{ \int_\R \Fb(\xb_2')\, d\xb_2'}\right )\, d\xb_1 d\xb_2 \, ,
\end{equation}
where $G_{\sigmab^2}$ denotes the centered Gaussian distribution with variance $\sigmab^2/2$. Here the parameter $\sigmab^2$ corresponds to  the genetic variance at linkage equilibrium in the absence of selection~\citep{Bul71,Lan78,Bul74,San98,TurBar94}. We can observe that, conversely to the previous asexual model, the infinitesimal model generates a finite standing variance in the absence of selection. Thus the dynamics of the standing variance are more constrained in the infinitesimal model than in the asexual model.
However, the genetic variance at linkage equilibrium $\sigmab^2/2$ plays an analogous role in our analysis, as the variance of the mutation kernel $\sigmab^2$ from the asexual model (that is why we use the same notation). In particular, they both determine the  phenotypic variance among offspring born to the same parents and therefore scale the input of phenotypic diversity through reproduction in the population.

\paragraph*{Equilibrium in a changing environment}
In this paper, we focus on the asymptotic behavior of the model, studying whether the population will persist or go extinct in the long term. In order to mathematically address  the problem, we seek special solutions  of the form $\fb(\tb,\xb) = \Fb(\xb-\cb\tb)$. These solutions correspond to a situation where the phenotypic distribution $\Fb$ has reached an equilibrium, which is shifted at the same speed $\cb$ as the environmental change. {This distribution of phenotypic lag $\zb:=\xb-\cb\tb$ quantifies \emph{maladaptation}.} One can easily observe from equation~\eqref{eq:f z} that the trivial solution, which corresponds to $\Fb = 0$, always exists. Our aim is first to decipher when non trivial equilibrium $\Fb$ exists. Secondly, we characterize in detail the distribution $\Fb$ when it exists.

A straightforward computation implies that a non trivial equilibrium $\Fb$ is the solution to the following eigenvalue problem,
\begin{equation}\label{eq:equilibrium c}
\lambdab \Fb(\zb) - \cb \partial_{\zb} \Fb(\zb)  + \mub(\zb) \Fb(\zb) = \betab \mathcal{B}( \Fb )(\zb)\, 
\end{equation}
with
\begin{equation}\label{eq:rhob}
   \lambdab=(\betab-\mub_0)\rhob  =(\betab-\mub_0)\int_\R\Fb(\zb')d\zb', 
\end{equation}
such that the eigenvalue $\lambdab>0$ .
The backward transport term $-\cb\partial_\zb \Fb$ corresponds to the effect of the moving optimum on the phenotypic distribution $\Fb$ at equilibrium. 
A formal integration of equation~\eqref{eq:equilibrium c} shows that 
$$\lambdab = \int_\R (\betab - \mub(\zb)) \dfrac{\Fb(\zb)}{\int_\R \Fb(\zb')d\zb'} d\zb. $$
The eigenvalue $\lambdab$ can thus be interpreted as a measure of the mean fitness of the population, or its mean intrinsic rate of increase, where $ \betab - \mub(\zb)$ is the contribution to population growth rate of an individual with phenotypic lag $\zb$ at low density. Thus a precise description of $\lambdab$ will provide a precise analytical formula for the critical speed of environmental change above which extinction is predicted  corresponding to the case where $\lambdab$ is negative.
The value $\lambdab$ also informs us on the size of the population at equilibrium in presence of a changing environment, $\rhob$ (see equation~\eqref{eq:rhob}).

Our aim is to describe accurately the couple solutions $(\lambdab,\Fb)$ in presence of a moving optimum with constant speed $\cb$ in both reproduction scenarios. To do so, we compute formal asymptotics of $(\lambdab, \Fb)$ at a weak selection or slow evolution limit when little variance in fitness is generated by mutation or sexual reproduction per generation.
Note that the shape of $\Fb$ is not prescribed {\em a priori} and our methodology presented here can  handle significantly large deviations from Gaussian distributions.

Noteworthy, the equation~\eqref{eq:equilibrium c} with asexual reproduction operators defined by~\eqref{eq:asexual B} or \eqref{eq:asexual B diffusion} admits solutions under suitable conditions. \cite{CloGab19} proved that solutions exist for any speed $\cb$ if the selection function $\mub$ goes to $\infty$ when $|z|\to\infty$. Furthermore, \cite{CovHam19} proved that solutions also exist for more general selection functions $\mub$ as soon as the speed $\cb$ remains below a critical threshold. For the infinitesimal operator~\eqref{eq:sexual B}, \cite{Pat20} proved the existence of solutions without changing environment and in the special regime of small variance described below. The existence of a pair $(\lambdab, \Fb)$ for positive speed $\cb$ will be the topic of a future mathematical paper.



\subsection{Rescaling the model}\label{sec:adim}
In order to compute asymptotics of the solution of our model, we first need to rescale the model with dimensionless parameters (see Table~\ref{table:scaling-var} for the relationship between original variables and their values after rescaling and Appendix~\ref{app:scaled_pb} for mathematical details).

\begin{table}
\begin{center}
\renewcommand{\arraystretch}{2.5}
\begin{tabular}{l l l}
{\sffamily  Parameters} & \sffamily Description &\sffamily Rescaled parameters  \\\hlineB{3}
\rowcolor{gray!30}
$\zb$        & phenotypic lag  & $z = \zb\sqrt{\alphab/\betab}$ \\
$\Fb \left ( \zb  \right )$        & phenotypic density & $F(z) = \Fb \left ( \zb  \right )$ \\
\rowcolor{gray!30}
$\betab $ & fertility rate &  $  1$\\
$\mb(\zb)$     & increment of mortality rate & $  m(z) =    \mb(\zb)/\betab$   \\ 
\rowcolor{gray!30}
$\alphab$ & strength of stabilizing selection & $ 1$\\
$\lambdab$  & { mean fitness} & $\lambda =(\lambdab + \mub_0)/\betab$ \\
\rowcolor{gray!30}
$\cb$        & speed of environmental change   &  $c = \left\{ \begin{array}{ll}
             \ds \cb/(\sigmab\betab) & \hbox{(asexual)}  \\[1mm]
             \ds \cb/(\sigmab^2\sqrt{\alphab\betab}) & \hbox{(infinitesimal sexual)}
           \end{array}\right. $  \\\hlineB{3}
$\sigmab^2$ & phenotypic variance parameter & $\eps^2 =\sigmab^2\alphab/\betab\ll1$\\\hlineB{3}
\end{tabular}
\end{center}
\caption{Biological parameters and their formula after rescaling for both the asexual and infinitesimal sexual model. 
Our methodology relies on the assumption that the adimentional parameter $\epsilon$ is small, $\epsilon\ll1$, while the other rescaled parameters can be chosen of order $1$.}
\label{table:scaling-var}
\end{table}

\paragraph*{Time scale.} All rates (and in particular fitness) have been rescaled so that time is expressed in the number of generations. We divide the equations by the generation time $1/ \betab$. 

\paragraph*{Phenotypic scale.}
All measures depending on phenotypic units have been rescaled to be dimensionless. We divide the equations by a phenotypic scale inversely related to the strength of stabilizing selection around the optimal phenotype $\alphab$ (see Table~\ref{table:scaling-var}), such that the strength of selection in the rescaled system is equal to unity:
\begin{equation}
m''(0)=\mb''(0)/\alphab = 1\, .
\end{equation} 

\paragraph*{Input of phenotypic variance through reproduction.}
Similarly, in both asexual and infinitesimal sexual models, the parameter describing how much phenotypic variance is introduced in the population by either mutation or recombination during reproduction has been rendered dimensionless. This phenotypic variance parameter $\sigma^2$ is divided by the same (squared) phenotypic scale $\betab/\alphab$, inversely related to the strength of stabilizing selection (see Table~\ref{table:scaling-var}). The ratio $\eps^2=\sigma^2\alphab/\betab$ appears naturally in the expression of the reproduction operator $\mathcal{B}$ in the scaled variables:
\begin{equation}\label{eq:asex sex}
   \mathcal{B}(F)(z) = 
   \left\{ \begin{array}{ll}
             \ds \dfrac{1}{\eps}\int_\R  K\left(\dfrac{z-z'}{\eps}\right)  F(   z')\, d  z' & \hbox{(asexual)}  \\[3mm]
             \ds \dfrac{1}{\eps \sqrt{\pi}}   \iint_{\mathbb{R}^2}  \exp\left ( -  \dfrac{1}{\eps^2} \left( z - \dfrac{ z_1 +   z_2}2 \right)^2 \right )  F( z_1)\dfrac{  F(  z_2)}{ \int_\R   F(  z_2')\, d  z_2'}\, d  z_1 d  z_2  & \hbox{(infinitesimal sexual)}
           \end{array}\right.
\end{equation}
Note that the parameter $\epsilon^2$ has a form similar to the standing load depressing mean fitness, as defined for instance by~\citep{LanSha96}.

\paragraph*{Speed of environmental change.}
The rescaling of the speed of environmental change however differs in the asexual and infinitesimal sexual versions of our model. In both models, the ability to evolve fast enough to track the moving optimum depends critically on the input of phenotypic variation fueling evolutionary change. It therefore makes sense that we may want to scale the speed of environmental change with respect to a measure of phenotypic diversity. 
However, we found that the scaling allowing analytical insights in the limit of small scaled genetic variance (see next section) differs in the two models.  We introduce the speed parameter $c$ so that the speed $\cb$ is scaled to, either $\epsilon$ in the scenario of asexual reproduction, or to $\epsilon^2$ in the  scenario of infinitesimal sexual reproduction (see Table~\ref{table:scaling-var}).

\paragraph*{Rescaled model.}
Using these rescaled variables,  we obtain the following equations: 
\medskip

\noindent{\bf Asexual reproduction}
\begin{equation}\label{eq:asexual no age}
  \lambda F(z) - \eps  c  \partial_z F(z)  + m(z) F(z) =  \dfrac{1}{\eps}\int_\R  K\left(\dfrac{z-z'}{\eps}\right)  F(   z')\, d  z'\, . 
\end{equation}
\noindent{\bf Infinitesimal sexual reproduction}
\begin{multline}
    \label{eq:sexual no age}
 \lambda F(z) - \eps^2  c  \partial_z F(z)  + m(z) F(z) =   \\
 \dfrac{1}{\eps \sqrt{\pi}}   \iint_{\mathbb{R}^2} \hspace{-3mm} \exp\left ( -  \dfrac{1}{\eps^2} \left( z - \dfrac{ z_1 +   z_2}2 \right)^2 \right )  F( z_1)\dfrac{  F(  z_2)}{  \int_\R   F(  z_2')\, d  z_2'}\, d  z_1 d  z_2 \, . 
\end{multline}

\subsection{Small variance asymptotics}\label{sec:small_varinace}
In the following, we further assume that the parameter $\eps$ is small, which means that very little variance in fitness is introduced in the population through either mutation or recombination during reproduction.
This is what we call {\em the small variance regime}. This situation may happen either when the input of phenotypic variation  is small or because stabilizing selection is weak. In other words, we can deal with relatively strong selection in amplitude as soon as the mean effect of mutations is assumed to be relatively small. 
Under such a regime ($\eps\ll 1$), we expect the  equilibrium  $F$ to be concentrated around a mean value for the rescaled phenotypic lag $z^*$, the \emph{evolutionary lag}, meaning that the rescaled standing variance at equilibrium is small.  
The core of our approach consists in the accurate description of the phenotypic distribution  $F$ in the limit of small standing variance, that is $\epsilon\ll 1$. 

The main ingredient is a suitable transformation of the phenotypic distribution $F$. As mentioned above, $F$ is expected to concentrate around a mean value, {the \emph{evolutionary lag} $z^*$} with standing variance depending on the ratio $\epsilon$. Then, it is natural to take the logarithm of the  density $F$, multiplied by a small parameter related to the expected standing variance (this would be straightforward if the distribution would be Gaussian, actually).  Accordingly, the following quantities are introduced, depending on the scenario:  
\begin{equation}\label{eq:log transform}
\begin{cases}
U = -\eps\log F \quad & \text{(asexual)}\\
U = -\eps^2 \log F \quad & \text{(infinitesimal sexual)}
\end{cases}
\end{equation}
Again, the discrepancy between the two scenarios is an outcome. This is the only possible scaling that gives rise to a non trivial limit in the regime $\eps\ll1$.
 
In order to describe $U$, we expand it with respect to $\epsilon$ as follows:
\begin{equation}\label{eq:ansatz} 
\left\{\begin{array}{l} 
U(z) = U_0(z) \textcolor{gray}{+ \eps^\gamma U_1(z) + o(\eps^\gamma)} \\
\lambda  = \lambda_0 \textcolor{gray}{+ \eps^\gamma \lambda_1  + o(\eps^\gamma)}
\end{array}\right. \ \hbox{ where } \ \gamma=\left\{\begin{array}{l}
                                                1 \ \hbox{(asexual)} \\
                                                2 \ \hbox{(infinitesimal sexual)} \\
                                             \end{array}
 \right.
\end{equation}
and $(\lambda_0,U_0)$ is the limit shape as $\epsilon\to0$, and $(\lambda_1,U_1)$ is the correction for small $\epsilon>0$. 
In the next sections~\ref{sec:asex} and~\ref{sec:sex}, we show, by formal arguments, that the function $U$ and the mean fitness $\lambda$ converge towards some non trivial function $U_0$ and some value $\lambda_0$ as $\eps\to 0$. 

The main advantage of our methodology is to bypass the resolution of the limit equation solved by $(\lambda_0,U_0)$, in order to compute directly relevant quantitative features, such as the { mean fitness} $\lambda_0$, the {evolutionary lag}  $z^*_0$, and the standing variance $\Var(F)$, which is related to $U_0$ by the following formula (derived in Appendix~\ref{sec:var_derivation}):
\begin{equation}\label{eq:variance}
\Var(F) = \dfrac{\eps^\gamma}{\partial_z^2 U_0(z_0^*)} + o(\eps^\gamma)\, .
\end{equation}

\section{Results in the regime of small variance}\label{sec:results_adim}

\subsection{The asexual model}\label{sec:asex}
Using the the logarithmic transformation \eqref{eq:log transform} to reformulate our problem~\eqref{eq:asexual no age} and the Taylor expansion of the pair $(\lambda,U)$ with $\gamma=1$, we show that the limit shape $(\lambda_0,U_0)$ satisfies the following problem (see appendix~\ref{app:Mastre_eq_asex}):
\begin{equation}
\lambda_0 + c \partial_z U_0(z) + \m(z) = 1 + H  \left(\partial_z  U_0(z) \right)  \, ,\label{eq:HJ asexual U0}
\end{equation}
where the Hamiltonian function $H$ is the two-sided Laplace transform of the mutation kernel $K$ up to a unit constant:
\begin{equation}\label{eq:H}
H(p) =  \int_\R   K\left(y\right) \exp\left( y p \right)\, dy - 1  \, .
\end{equation}
It is a convex function that satisfies $H(0) = H'(0) = 0$, and $H''(0) = 1$ from hypothesis~\eqref{eq:asexual B} on the mutation kernel $K$.

We can remark that the shape of the equation also holds true for {\em the diffusion approximation model} where the reproduction operator is approximated by a diffusion operator~\eqref{eq:asexual B diffusion}. For the diffusion approximation, we find that the Hamiltonian satisfies $H(p) = p^2/2$ (see appendix~\ref{app:diffusion_approx})

\paragraph*{Computation of the {mean fitness} $\lambda_0$.}
We find that (see Appendix~\ref{app:mean_fitness_asex} for details)
\begin{equation}\label{eq:lambda0}
    \lambda_0   = 1 - L(c)\,,
\end{equation}
where the Lagrangian function $L$ known as the Legendre transform of the Hamiltonian function $H$, is defined as: 
\begin{equation} \label{eq:L}
L(v) = \max_{p\in \R} \left( p v - H(p)\right)\,. 
\end{equation}
It is a convex function satisfying $L(0) = L'(0) = 0$, and $ L''(0) = 1$. Moreover, we always have $L(v)\leq |v|^2/2$ where $L(v) = |v|^2/2$ corresponds to the \emph{diffusion approximation} case. 
{ 

Since the mean fitness $\lambda_0=1$ in absence of environmental change, the quantity $L(c)$ represents the \emph{lag-load} in the rescaled units, which is induced by the moving optimum~\citep{LynLan93,LanSha96}. Moreover, if we push the expansion to the higher order we are able to compute the following mean fitness
\begin{equation}\label{eq:lambda_1}
    \lambda = 1  - L(c) \textcolor{gray}{- \dfrac\eps2   \left (\dfrac{ 1}{L''(c)}\right )^{1/2}   + o(\eps)} 
\end{equation}
The new term of order $\epsilon$ can be seen as the \emph{standing load}, i.e. a reduction in mean fitness due to segregating variance for the trait in the population, which has been introduced in~\citep{LynLan93,BurLyn95,KopMat14}. 
}


\paragraph*{Computation of the evolutionary lag $z^*_0$.} We obtain from the main equation~\eqref{eq:HJ asexual U0}, evaluated at $z = z^*_0,$ that $\lambda_0 + \m(z^*_0) = 1$. Thus, combining with equation~\eqref{eq:lambda0}, we deduce that $z^*_0$ is a root of
\begin{equation} \label{eq: def z0star}
\m(z^*_0) =  L(c)
\end{equation}
with the appropriate  sign, that is $\partial_z \m(z^*_0)$ and $c$ have opposite signs: $z^*_0<0$ if $c>0$ and vice-versa. 


\paragraph*{Computation of the standing variance.}
From equation~\eqref{eq:variance}, we need to compute the second derivative of $U_0$ at the  evolutionary lag $z^*_0.$ We can derive it from the differentiation of  equation~\eqref{eq:HJ asexual U0} evaluated at $z^*_0$ 
(recall that $H'(0) = 0$ by symmetry of the mutation kernel $K$):
\begin{equation}\label{eq:d2U_0} 
\partial^2_z U_0(z^*_0) + \dfrac{\m'(z^*_0)}{c} = 0\, . 
\end{equation}
We deduce the following first order approximation  of the standing  variance:
\begin{equation}
\Var(F) =  - \dfrac{\eps  c}{m'(z_0^*)} + o(\eps)\, .
\end{equation}

\begin{remark}\label{rem:local shape}
The expressions obtained in this section are still valid when $c=0.$ A direct evaluation gives that $\lambda_0=1$ and $z^*_0=0.$ Moreover, we show in Appendix \ref{app:asex_c0} that in the limit $c\to0,$ the previous formula~\eqref{eq:d2U_0} becomes 
\begin{equation}\label{eq:eq:d2U_0 c0} \partial^2_z U_0(0) = 1\, . \end{equation}   
\end{remark}




We will discuss the biological implications of these predictions after expressing them in the original units in the section~\ref{sec:dim}.

\subsection{The infinitesimal model of sexual reproduction in the regime of small variance}\label{sec:sex}

\paragraph*{The limiting problem formulation.}
Remarkably enough, a similar mathematical analysis can be performed when the convolution operator is replaced with the infinitesimal  model for reproduction \eqref{eq:asex sex}. However, the calculations are slightly more involved than the former case, but the final result is somewhat simpler. Here, the suitable logarithmic transformation of the phenotypic distribution $F$ is $U=-\epsilon^2\log(F)$. 
The equation for the new unknown function $U$ is:
\begin{multline}
    \label{eq:U sexual}
\lambda + c\partial_z U(z) + m(z)  =    \\
\f{ \dfrac{1}{
\eps \sqrt{\pi}}\displaystyle \iint_{\R^2} \exp\left(-\dfrac1{\eps^2} \left[ \left(z-\dfrac{z_1+z_2}{2}\right)^2 + U(z_1) + U(z_2)  - U(z) - \min U \right] \right) d z_1 d z_2}{\displaystyle \int_\R \exp\left(-\dfrac{U(z') - \min U}{\eps^2} \right) d z'}  \, , 
\end{multline}
where $\min U$ has been subtracted both in the numerator and the denominator. 
The specific form of the right-hand-side characterizes the shape of $U$. Indeed, the quantity between brackets must remain non negative, unless the integral  takes arbitrarily large values as $\eps \to 0$. Moreover, its minimum value over $(z_1,z_2)\in \R^2$ must be zero, unless the integral vanishes. As a consequence, the function $U$ must be a quadratic function of the form $\frac12(z - z^*_0)^2$ where the evolutionary lag of the distribution $z_0^*$ will be determined aside (see Appendix~\ref{app:sex_U0} for details).
To describe $z^*_0$, we expand the pair $(\lambda,U)$, in a power series  with respect to $\eps^2$:
\begin{equation}\label{eq:ansatz sexual} \begin{cases} 
U(z) = \dfrac12\left (z - z_0^*\right )^2 + \eps^2 U_1(z) \textcolor{gray}{+ \eps^4 U_2(z) + o(\eps^4)} \smallskip\\
\lambda = \lambda_0 + \eps^2 \lambda_1 \textcolor{gray}{+ \eps^4\lambda_ 2 + o(\eps^4)}
\end{cases}
\end{equation}
Plugging this expansion into~\eqref{eq:U sexual}, we obtain the following equation on the corrector $U_1$:
\begin{equation}\label{eq:U1 main}
\lambda_0 + c (z - z_0^*)  + \m(z) = \exp\left( U_1(z_0^*) -2 U_1\left ( \frac{z + z_0^*}{2} \right ) + U_1(z) \right)\, ,
\end{equation}
which contains as a by--product the value of some quantities of interest, such as the mean fitness $\lambda_0$, and the evolutionary lag $z^*_0$. Moreover, we can solve this equation only if $\lambda_0$ and $z_0^*$ takes specific values that we identify below.

\paragraph*{Computation of macroscopic quantities.}
Let us first observe that equation \eqref{eq:U1 main} is equivalent to the following one:  
\begin{equation}\label{eq: U1 log}
\log \left(\lambda_0 + c (z - z_0^*)  + \m(z)\right) =   U_1(z_0^*) -2 U_1\left ( \frac{z + z_0^*}{2} \right ) + U_1(z) \,  .
\end{equation}
The key observation is that the expression on the right hand side  vanishes at $z = z_0^*$, and so does its first derivative with respect to $z$ at $z=z_0^*$. 
This provides two equations for the two unknowns $\lambda_0,z_0^*$, without computing the exact form of $U_1$: 
\begin{equation}\label{eq:sex_lambda_z}
\begin{cases}
\lambda_0 + \m(z_0^*) = 1 \medskip\\
c + \m'(z_0^*) = 0\,.
\end{cases}
\end{equation}
These two relationships are necessary and sufficient conditions, meaning that they guarantee that equation \eqref{eq:U1 main} admits at least one solution $U_1$ (see Appendix~\ref{app:sex_U0} for mathematical details). 
In addition, we can push the expansion further and we can gain access to the higher order of approximation for the quantities of interest (see Appendix~\ref{app:sex_U1}). 

\begin{equation}\label{eq:sexual formula}
\begin{array}{ll}
\text{\bf Evolutionary lag  } & z^* = z^*_0 \textcolor{gray}{- \epsilon^2 \left(\dfrac{\m'''(z_0^*)}{2 \m''(z_0^*) } + 2c \right) + o(\eps^2)}
\, ,\quad \text{such that}\quad  m'(z_0^*) = -c \medskip\\
\text{\bf Mean fitness } & 
\ds \lambda =1 - m(z^*_0) \textcolor{gray}{- \epsilon^2 \left(  2c^2 + c\dfrac{ \m'''(z_0^*)}{ 2 \m''(z_0^*)} + \frac12m''( z_0^*) \right) + o(\eps^2)}
  \medskip\\

\text{\bf Standing variance} & \Var(F) = \dfrac{\epsilon^2}{1 \textcolor{gray}{+ 2\epsilon^2 m''(z^*_0) + o(\eps^2) }}
\end{array}
\end{equation}

\section{Comparison of predictions of the asexual and infinitesimal models}\label{sec:dim}
To discuss our mathematical results from a biological perspective, we need to scale back the results in the original units (see Table~\ref{table:scaling-var} for the link between the scaled parameters and the parameters in the original units). Our general predictions for macroscopic quantities in the original units are shown in Table ~\ref{tab-summary}. For ease of comparison with previous literature, which has generally assumed a quadratic form for the selection function, we present our predictions in Table~\ref{tab-summary-quadra} under this special assumption and with the diffusion approximation. 

\paragraph*{Numerical simulations.}
To illustrate our discussion, we also perform numerical simulations. The simulated stationary distribution is obtained through long time simulations of a suitable numerical scheme for \eqref{eq:f z} (details in Appendix \ref{app:num_sol_F_lambda}). Using this numerical expression, we compute the lag, the mean fitness and the standing variance of the distribution. In the asexual model, the function $U_0$ is obtained  from the direct resolution of the ordinary differential equation~\eqref{eq:HJ asexual U0} using classical integration methods -- see Appendix~\ref{app:num_sol_asex_U0_U1}. In the infinitesimal model, the correction $U_1$ is computed directly from its analytical expression given in Appendix~\ref{app:p*}. The macroscopic quantities in the regime of small variance are directly computed from their analytical expressions given in the Table~\ref{tab-summary} and~\ref{tab-summary-quadra}.



\begin{table}
\begin{center}
\renewcommand{\arraystretch}{2.5}
\begin{tabular}{p{2cm}|l|l}
\parbox{3cm}{\sffamily  Macroscopic \\ quantities} &  \sffamily Asexual model &\sffamily  Infinitesimal sexual model \\\hlineB{3}
\rowcolor{gray!30}
Evolutionary lag       
& \parbox{5.5cm}{\begin{equation*}\label{eq:asex_zstar_original}
     \begin{array}{l}
        \zb^* \approx \zb^*_0 \\
        \hbox{ with } \ \mb(\zb^*_0) = \betab L\left(\frac{\cb}{\sigmab\betab}\right) 
     \end{array}
  \end{equation*} }
& \parbox{7cm}{\begin{equation*}\label{eq:sex_zstar_original}
    \begin{array}{l}
        \zb^* \approx \zb^*_0 \textcolor{gray}{ -\sigmab^2  \dfrac{ \mb'''(  \zb_0^*)}{ 2 \mb''(  \zb_0^*)} - 2\dfrac{\cb}{\betab}} \\
        \hbox{ with } \ \mb'(  \zb^*_0) = - \frac{\cb}{\sigmab^2} 
     \end{array}
  \end{equation*}} \\

Mean fitness        
& \parbox{5.5cm}{\begin{equation*}\label{eq:asex_lambda_original}
     \begin{array}{rl}
       \lambdab \approx &  \betab -\mub_0  -  \betab L\left (\frac{\cb}{\sigmab\betab} \right )\\
                        &\textcolor{gray}{- \dfrac{1}{2}  \left ( \frac{\sigmab^2 \alphab\betab}{L''\left (   \frac{\cb}{\sigmab\betab}\right )}   \right )^{1/2}}
     \end{array}
  \end{equation*} }
& \parbox{7cm}{
\begin{equation*}\label{eq:sex_lambda_original}
     \begin{array}{l}
       \lambdab \approx  \betab -\mub_0  - \mb(\zb^*_0) \\
                         \textcolor{gray}{  -   \left ( \frac{2\cb^2}{\sigmab^2 \betab}+  \cb \dfrac{ \mb'''(  \zb_0^*)}{ 2 \mb''(  \zb_0^*)}   + \frac{\sigmab^2\mb''(  \zb_0^*)}{2}  \right )}
     \end{array}
  \end{equation*}}\\

\rowcolor{gray!30}
\parbox{3cm}{Standing \\ variance} 
& \parbox{5.5cm}{\begin{equation*}\label{eq:asex_Var_original}
     \Var(\Fb) \approx - \dfrac{\cb}{\mb'( \zb_0^*)} 
  \end{equation*} }
& \parbox{7cm}{\begin{equation*}\label{eq:sex_Var_original}
     \Var(\Fb) \approx \dfrac{\sigmab^2 }{1  \textcolor{gray}{+ 2\dfrac{\sigmab^2 }{\betab}\mb''(\zb_0^*) }}
  \end{equation*}}
 \end{tabular}
\end{center}
\caption{Analytical predictions for the evolutionary lag $\zb^*$, the mean fitness $\lambdab$ and the standing  variance $\Var(\Fb)$ for both the asexual and infinitesimal sexual model in the original variable. In the asexual model, $L$ is the Lagrangian defined by~\eqref{eq:L} and it is associated to the mutation kernel $K$.
}
\label{tab-summary}
\end{table}

\begin{table}
\begin{center}
\renewcommand{\arraystretch}{2.5}
\begin{tabular}{p{2cm}|l|l}
\parbox{2cm}{\sffamily  Macroscopic quantities} & \parbox{6cm}{\sffamily Asexual model \\ \textit{(quadratic selection / diffusive approx)}} &\parbox{6cm}{\sffamily Infinitesimal sexual  model\\ \textit{(quadratic selection)}} \\\hlineB{3}
\rowcolor{gray!30}
Evolutionary lag      
& \parbox{6cm}{\begin{equation*} \label{eq:asex_zstar_quad_diffapprox}
\zb^* \approx  - \dfrac{\cb}{\sigmab \left ( \alphab\betab  \right )^{1/2}}                     
                 \end{equation*}  }
& \parbox{7cm}{\begin{equation*} \label{eq:sex_zstar_quad_diffapprox}
\zb^* \approx -\dfrac{\cb}{\sigmab^2 \alphab} \textcolor{gray}{  - 2\dfrac{\cb}{\betab} }                     
                 \end{equation*}}\\
  
Mean fitness        
& \parbox{6cm}{\begin{equation*}\label{eq:asex_lambda_quad_diffapprox}
                    \lambdab \approx \betab -\mub_0 -   \frac{\cb^2}{2\sigmab^2\betab}  \textcolor{gray}{- \dfrac{\sigmab\left ( \alphab\betab  \right )^{1/2}}{2}}
                 \end{equation*}}
& \parbox{7cm}{\begin{equation*} \label{eq:sex_lambda_quad_diffapprox}
\lambdab \approx \betab -\mub_0  - \frac{\cb^2}{2  \sigmab^4\alphab} \textcolor{gray}{  -   \left (  \frac{2\cb^2}{\sigmab^2 \betab}   + \frac{\sigmab^2\alphab}{2}  \right )}                    
                 \end{equation*}} \\
                 
\rowcolor{gray!30}
Standing variance 
& \parbox{6cm}{\begin{equation*} \label{eq:asex_Var_quad_diffapprox}
\Var(\Fb) \approx  \sigmab\left (\dfrac{\betab}{\alphab}\right )^{1/2}                      
                 \end{equation*}}
& \parbox{7cm}{\begin{equation*}\label{eq:sex_Var_quad_diffapprox}
                    \Var(\Fb) \approx \dfrac{\sigmab^2}{1 \textcolor{gray}{  + 2\dfrac{\sigmab^2\alphab}{\betab} }}
                 \end{equation*}}
\end{tabular}
\end{center}
\caption{Analytical predictions for the evolutionary lag $\zb^*$, the mean fitness $\lambdab$ and the standing  variance $\Var(\Fb)$ for both the asexual and infinitesimal sexual model in the original variable when assuming a quadratic form of selection $\mb(\zb) = \alphab|\zb|^2/2$. In the asexual model, we are under the diffusion approximation: $L(v) = |v|^2/2$.}
\label{tab-summary-quadra}
\end{table}

\newpage

\subsection{Evolutionary lag}\label{sec:lag}
The lag is defined as the absolute value of the difference between the optimal trait (here set to $0$) and the evolutionary lag of the population $\zb^*$, thus the lag is equal to $|\zb^*|$. 
 
  \begin{figure}[t!]
      \subfigure[
      ]{\includegraphics[width=0.45\linewidth]{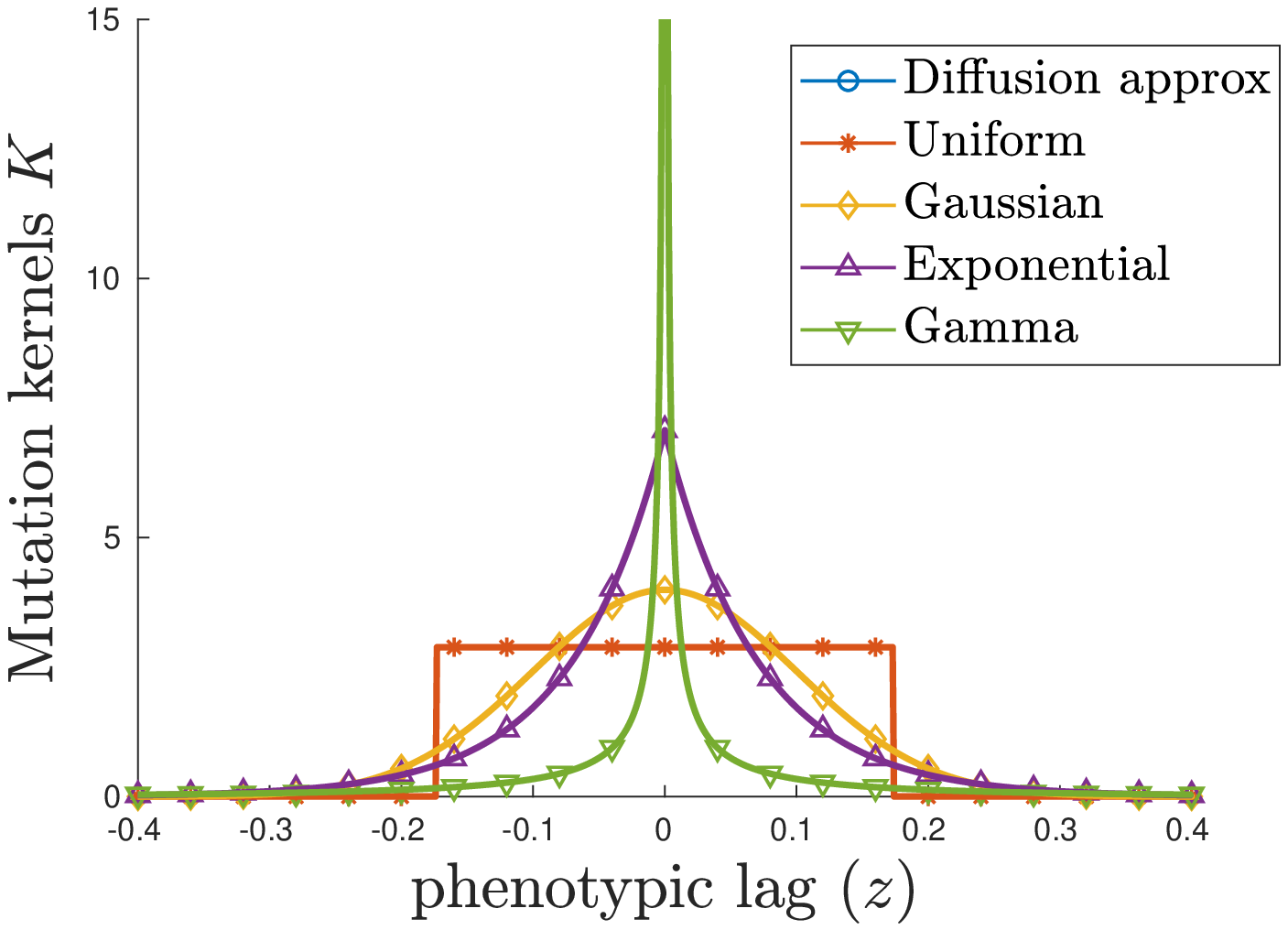}} \hspace{7mm}
    \subfigure[
    ]{\includegraphics[width=0.45\linewidth]{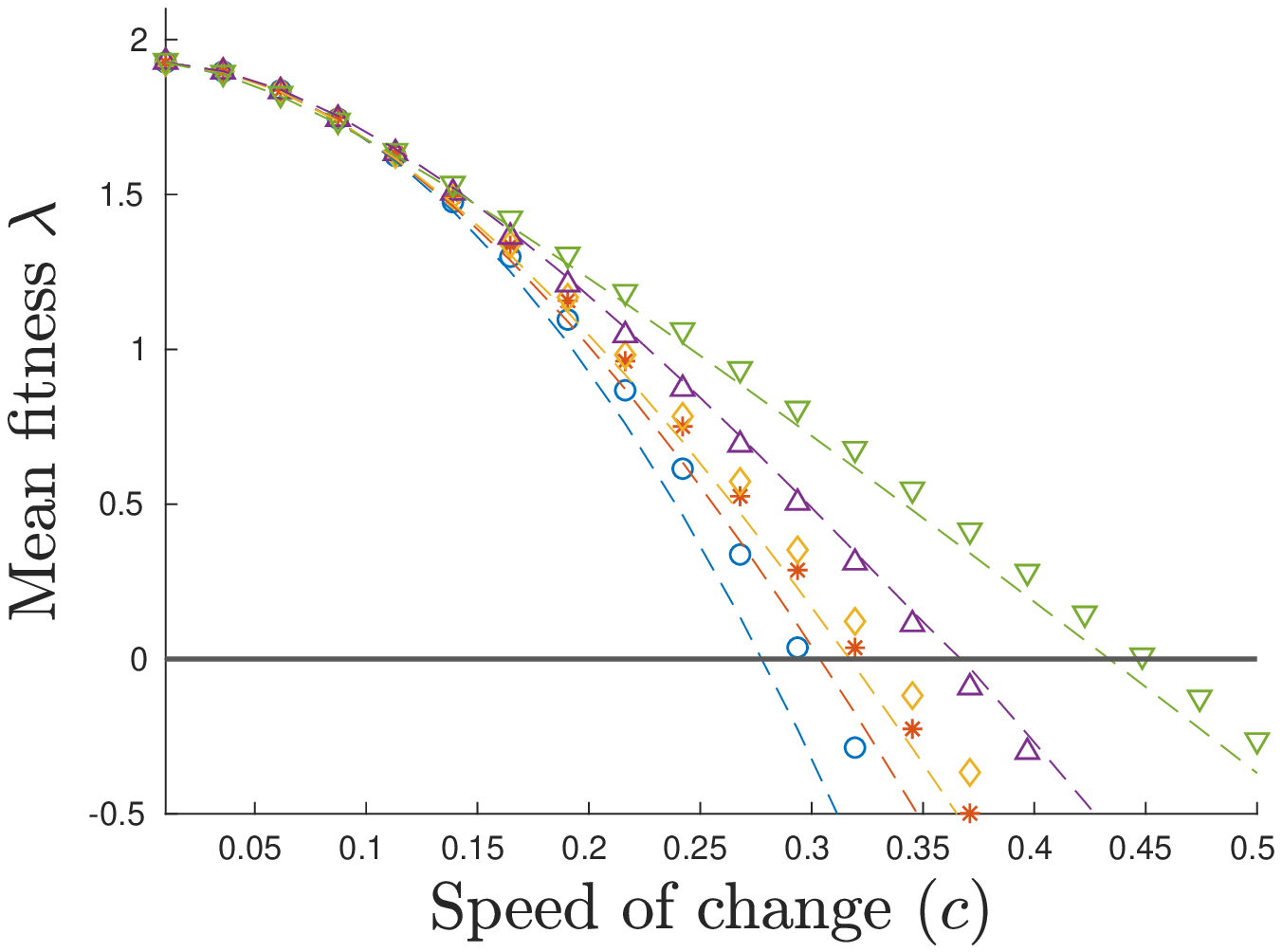}}
    \subfigure[
    ]{\includegraphics[width=0.45\linewidth]{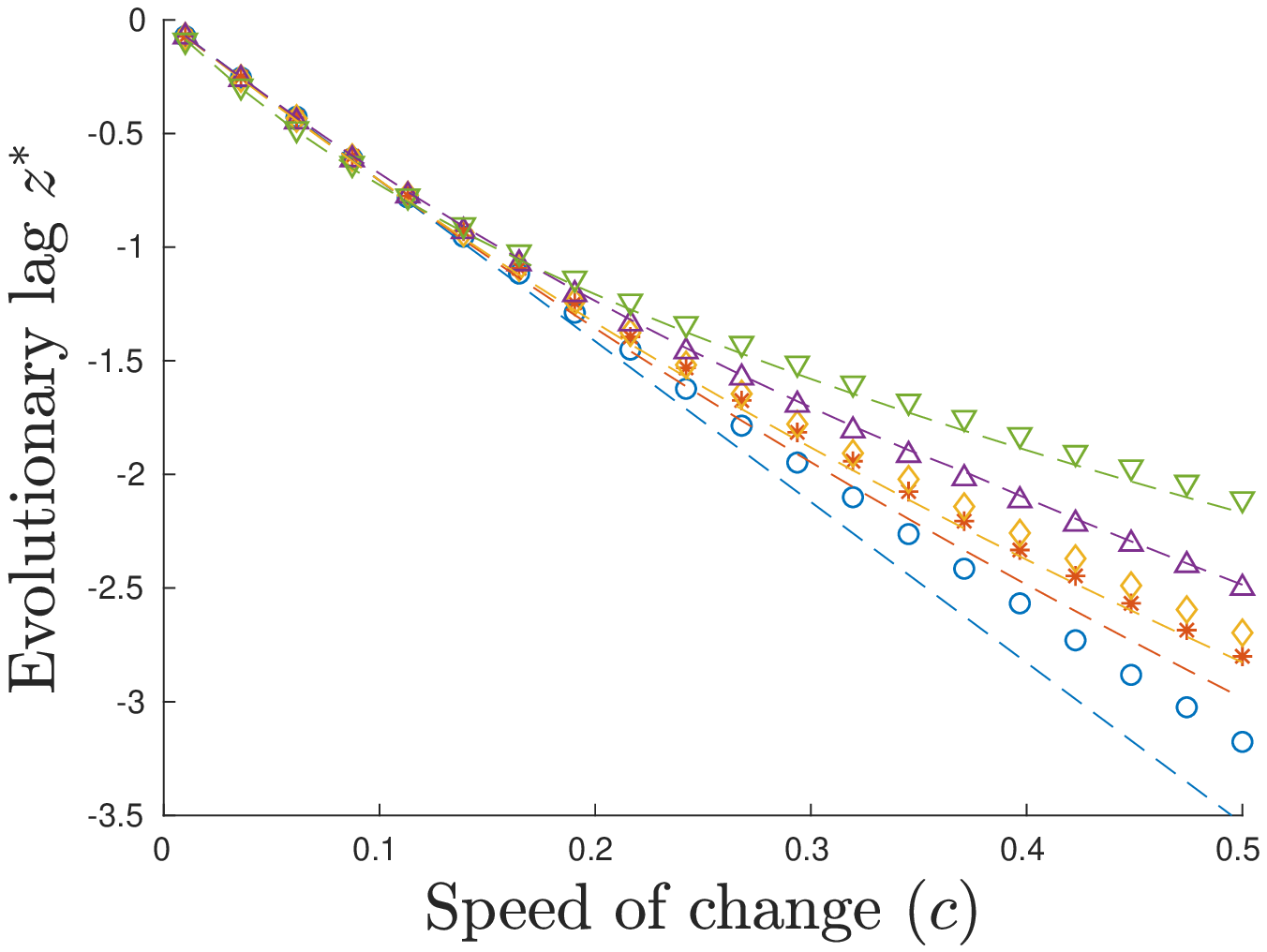}} \hspace{7mm}
    \subfigure[
    ]{\includegraphics[width=0.45\linewidth]{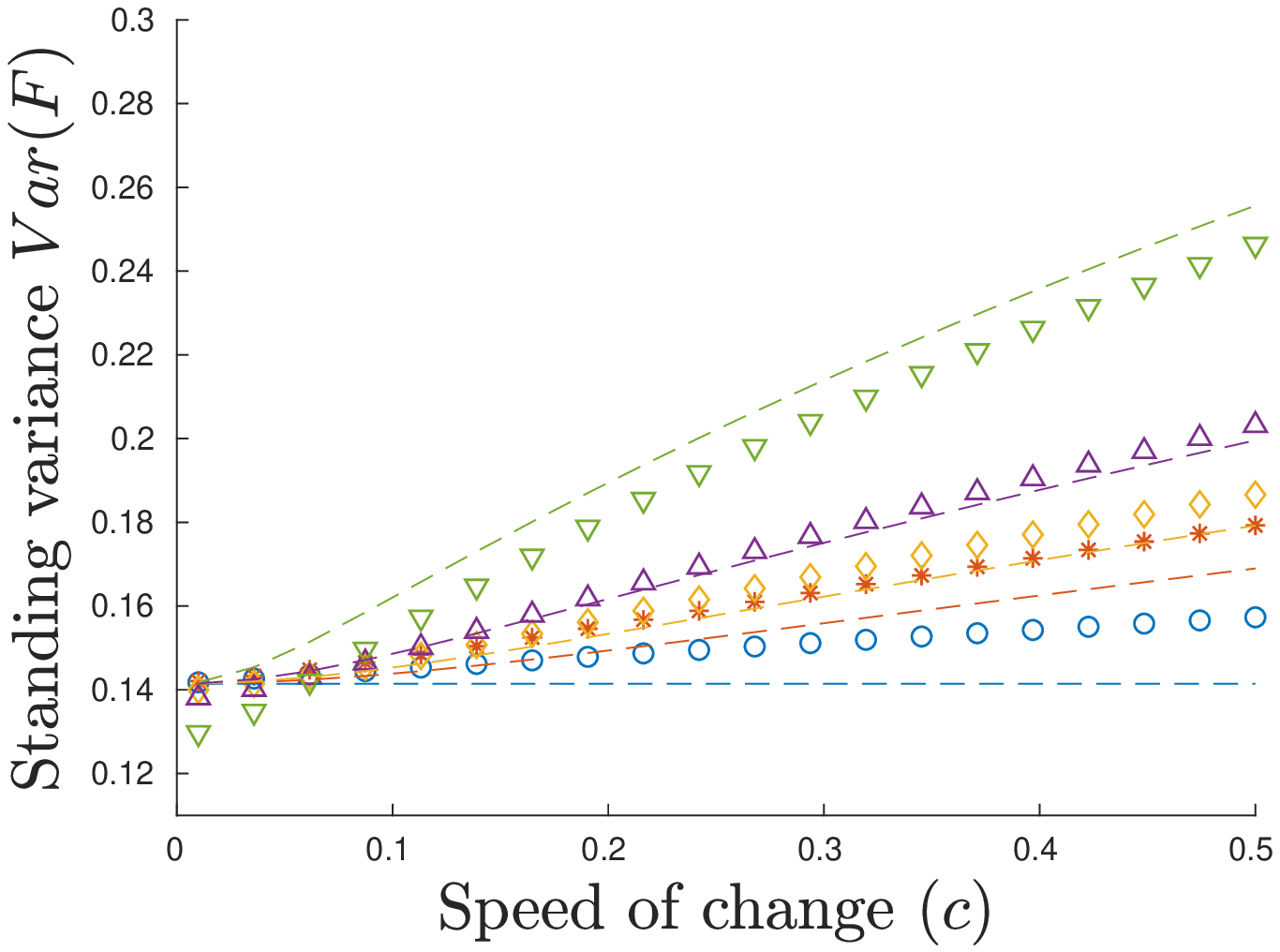}}

    \caption{Influence of the mutational kernel $K$, described in panel (a), on (b)  the mean fitness $\lambda$, (c) the evolutionary lag $z^*$ and (d) the standing phenotypic variance $\Var(F)$ at equilibrium in an environment changing at rate $c$ ranging in $(0,0.5)$. We compare the diffusion approximation (blue curves) with four different mutation kernels with the same variance $\sigmab=0.1$: the Uniform distribution (red curves), the Gaussian distribution (orange curves), Exponential distribution (purple curves) and Gamma distribution (green curves). Other parameters are: $\alphab = 1$; $\betab = 2$. For each case we compare our analytical results (dashed lines) with the simulation results (marked symbol).}\label{fig:asex_kernel_C}
 \end{figure}
 
  \begin{figure}[t!]
  \makebox[0.45\linewidth][c]{Asexual model} \hfil
  \makebox[0.45\linewidth][c]{Infinitesimal sexual model} \\[3mm]
  \subfigure[]{\includegraphics[width = 0.45\linewidth]{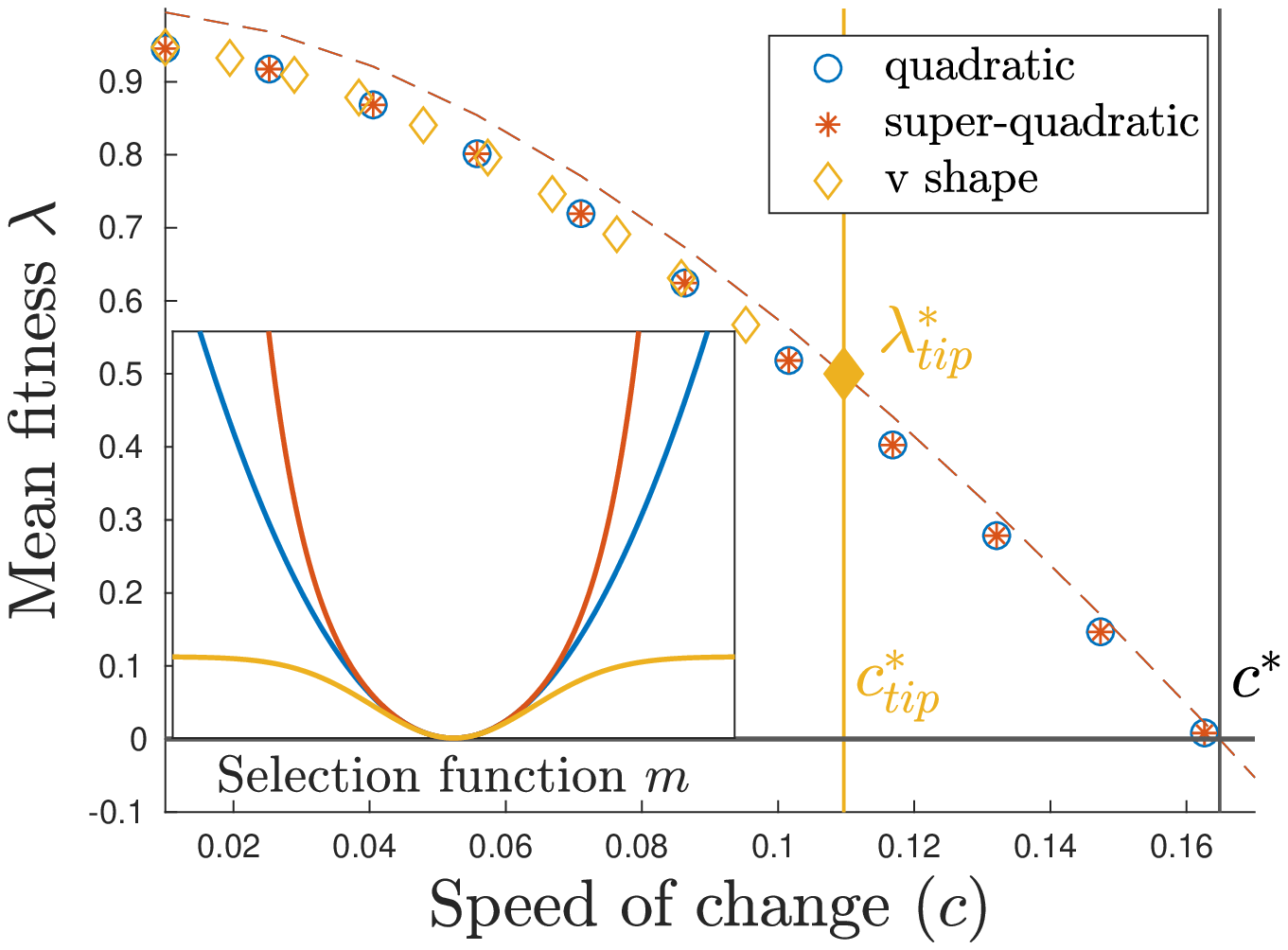}}
       \hfil
   \subfigure[]{\includegraphics[width = 0.45\linewidth]{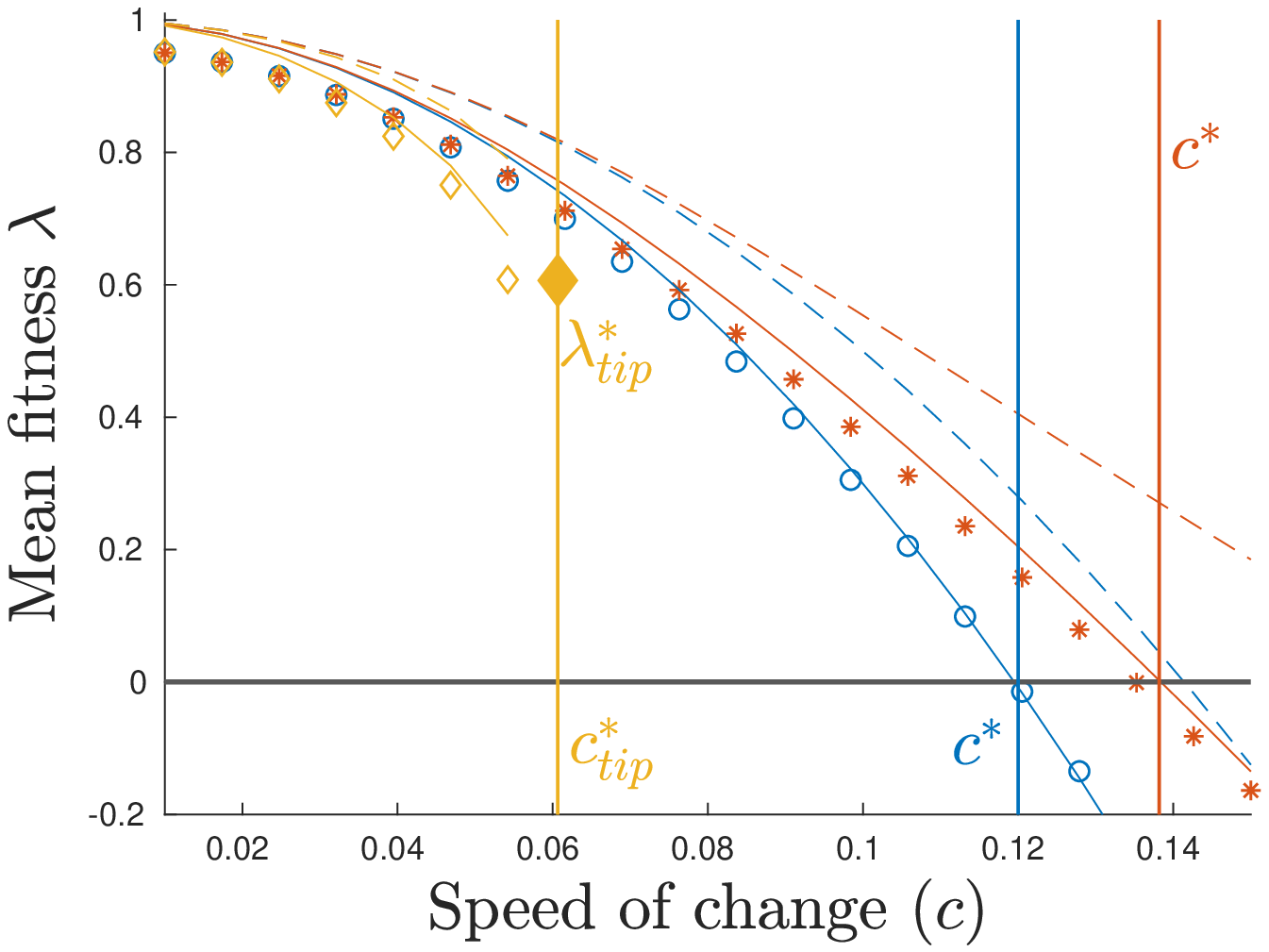}}\\ 
   \subfigure[]{\includegraphics[width = 0.45\linewidth]{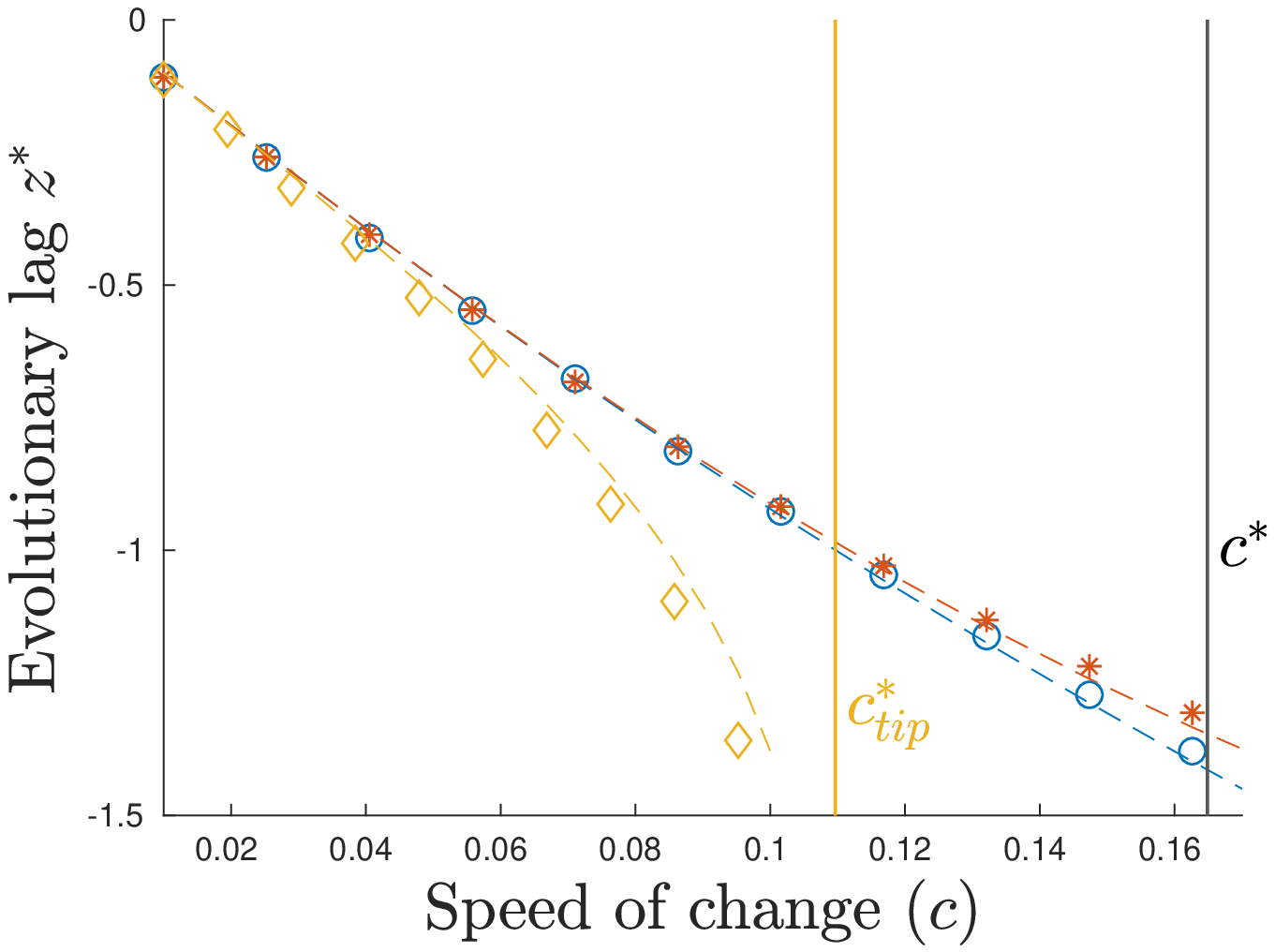}}
   \hfil
   \subfigure[]{\includegraphics[width = 0.45\linewidth]{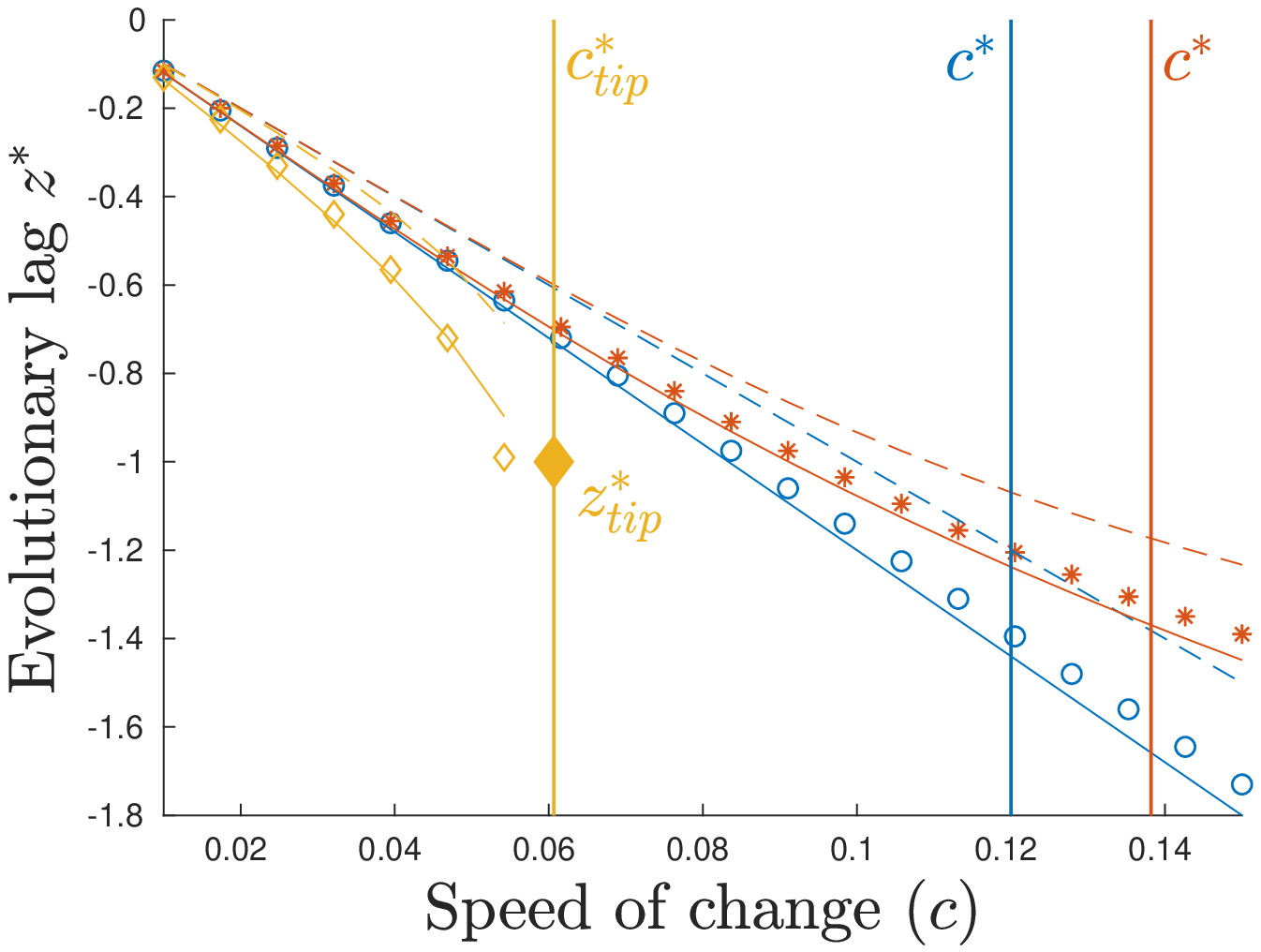}}\\
   \subfigure[]{\includegraphics[width = 0.45\linewidth]{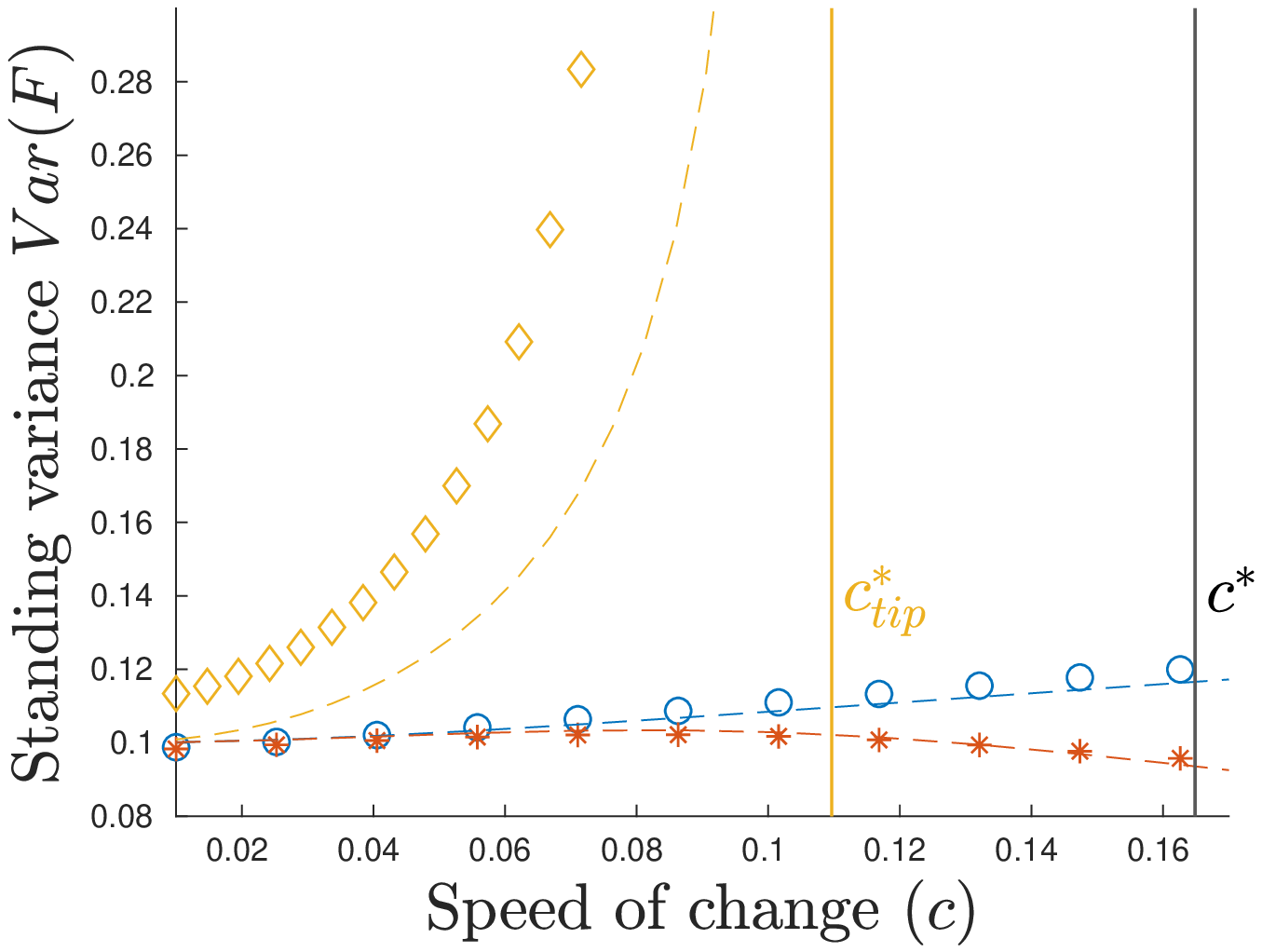}}
   \hfil
   \subfigure[]{\includegraphics[width = 0.45\linewidth]{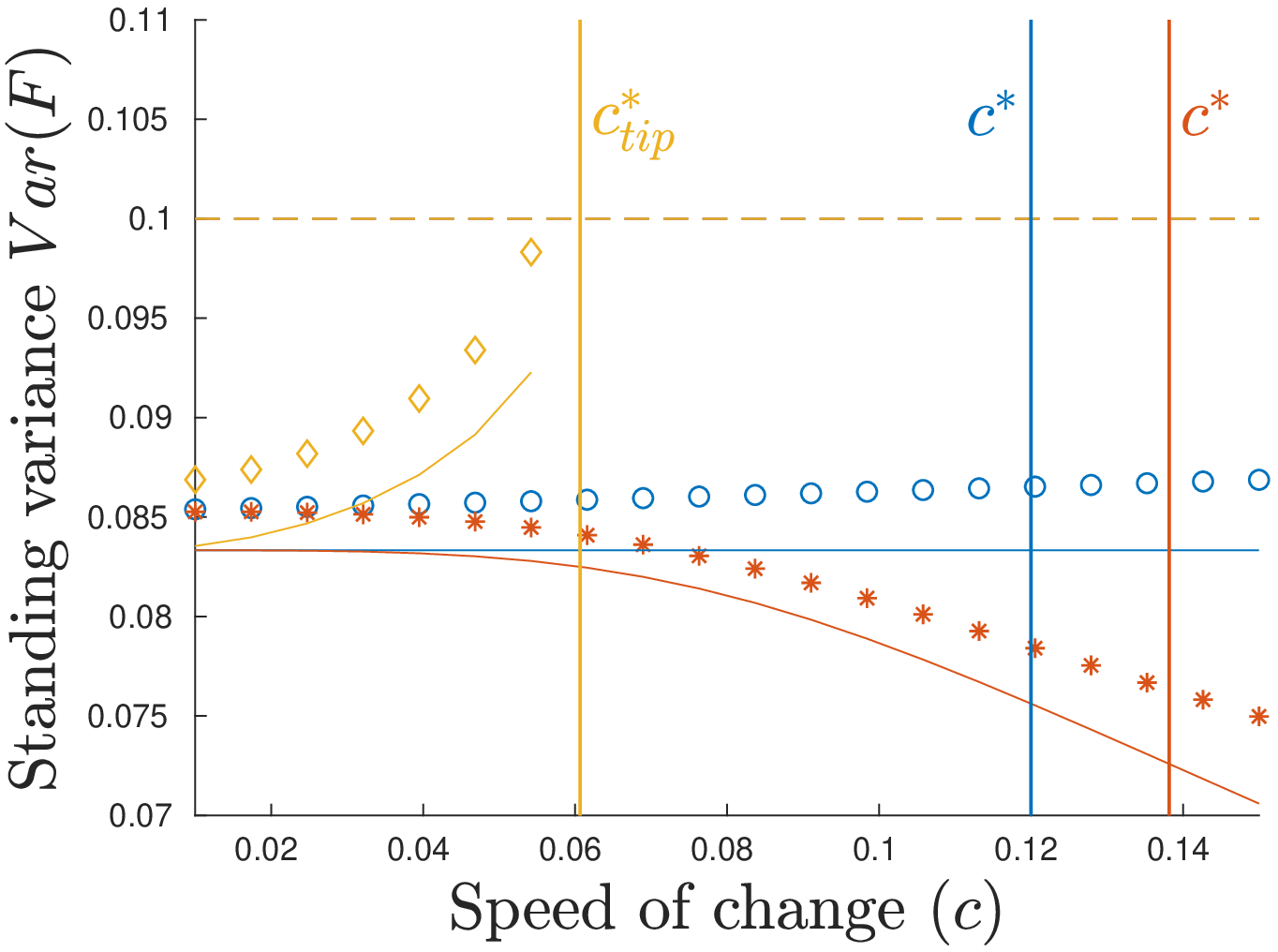}} 
   \caption{Influence of the speed of environmental change $\cb$ for three different selection function $\mb$: quadratic function $\mb(\zb) = \alphab\zb^2/2$ (blue curves), super--quadratic $\mb(\zb) = \alphab\zb^2/2 + z^6/64$ (red curves) or bounded $\mb(\zb) =m_\infty( 1-\exp(-\alphab\zb^2/(2m_\infty))$ (orange curves). Other parameters are: $\alphab = 1$, $\betab = 1$, $\sigmab = 0.1$ and $m_\infty=0.5$ in the asexual model and $m_\infty=1$ in the infinitesimal sexual model. In the asexual model, the mutation kernel is Gaussian.
   We compare our analytical results (first approximation dashed lines and second approximation plain lines) with the numerical simulations of the stationary distribution of \eqref{eq:equilibrium c} (marked symbols) for both asexual and sexual infinitesimal model. 
}\label{fig:m_c_selection}
\end{figure}

\begin{figure}[t!]
   \subfigure[Determination of the lag $\zb_0^*$ in the asexual model]{
\includegraphics[width = 0.45\linewidth]{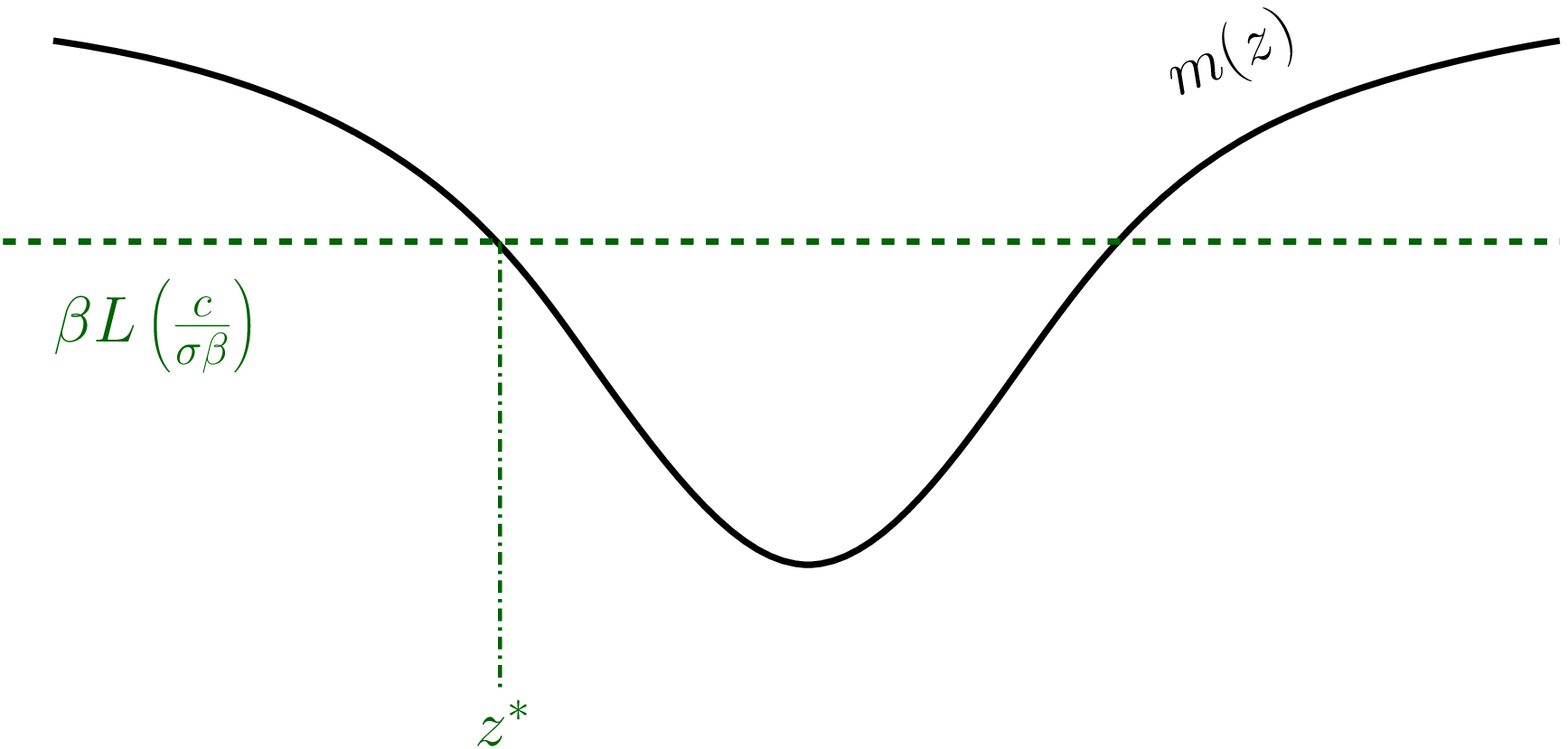}  \label{fig:formula_asex}
}
\hspace{7mm}
\subfigure[Determination  of two  possible lags $\zb^*_s$ and $\zb_u^*$ in the infinitesimal sexual model]{
\includegraphics[width = 0.45\linewidth]{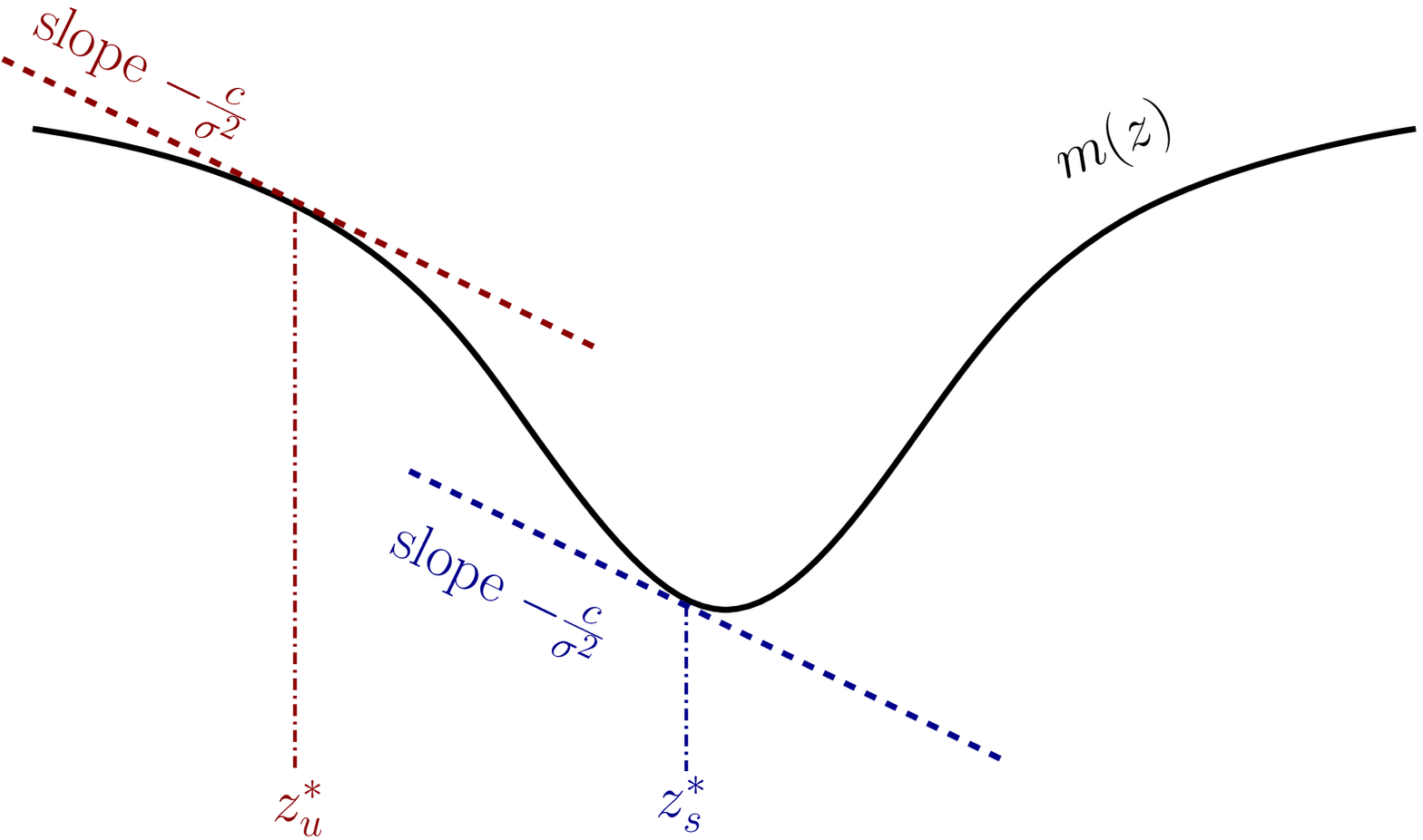} \label{fig:formula_sex}
}
\caption{Graphical illustration of the two ways to characterize the evolutionary lag $\zb^*$ (in original units). 
(a) In the asexual model, the evolutionary lag is found where the mortality rate $\mb$ equals a specific value $\betab L(\frac{\cb}{\sigmab\betab})$. In this case we only have one possible lag $\zb^*$ because $\mb'(\zb^*)$ and $c$ should have opposite signs.
(b) In the sexual infinitesimal model, the evolutionary lag $\zb^*$ is found where the selection gradient $\mb'$ equals a specific value $\frac{-\cb}{\sigmab^2}$. In this case, we may obtain two possible values, a stable point  $\zb^*_s$ in the convex part of $\mb$ and an unstable point $\zb^*_u$ in its concave part.}
\label{fig:formul_m_bounded}
\end{figure}

 \begin{figure}[t!]
   \subfigure[Asexual model]{
\includegraphics[width = 0.5\linewidth]{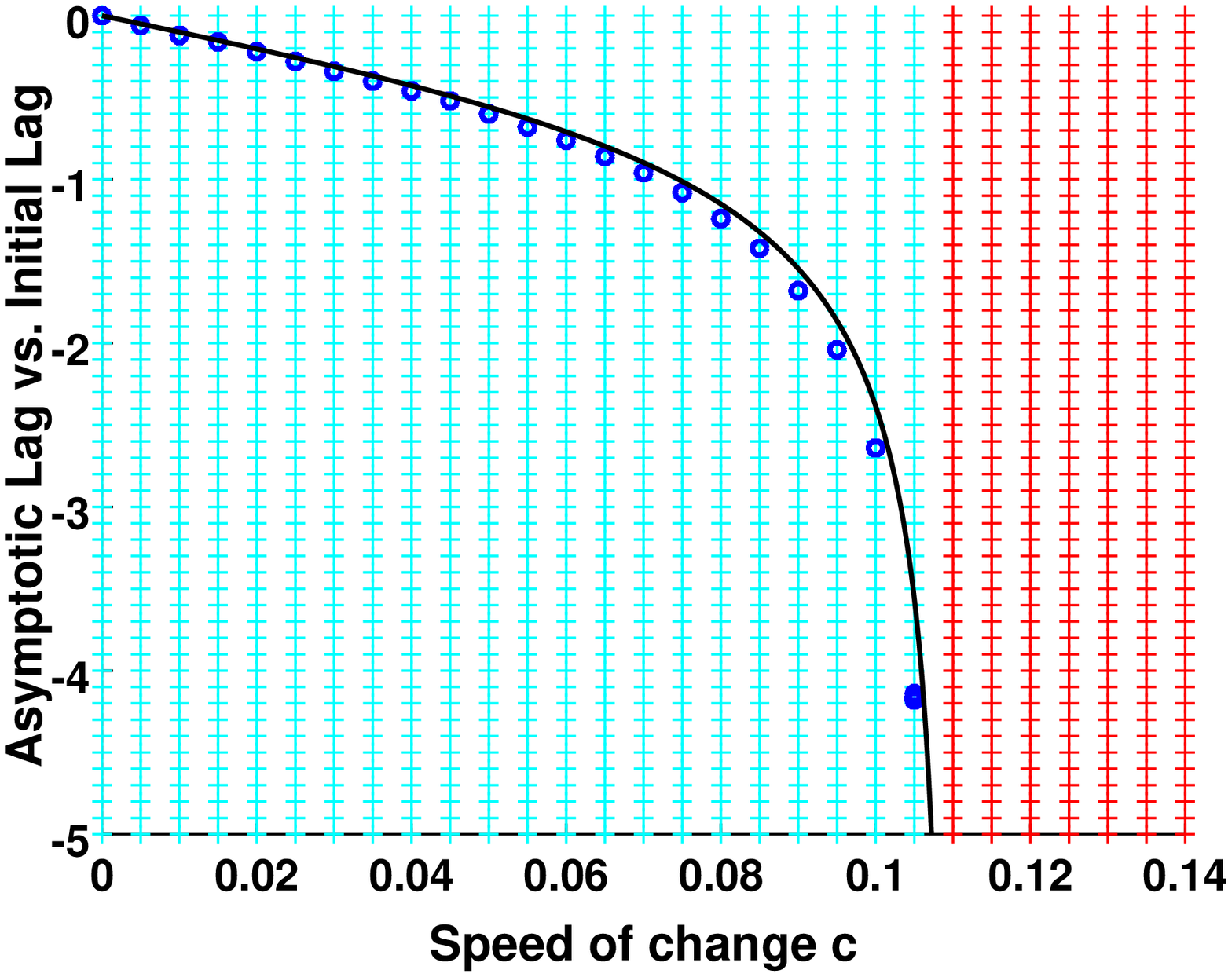}  \label{fig:stability}
}
\hspace{7mm}
\subfigure[Infinitesimal sexual model]{
\includegraphics[width = 0.5\linewidth]{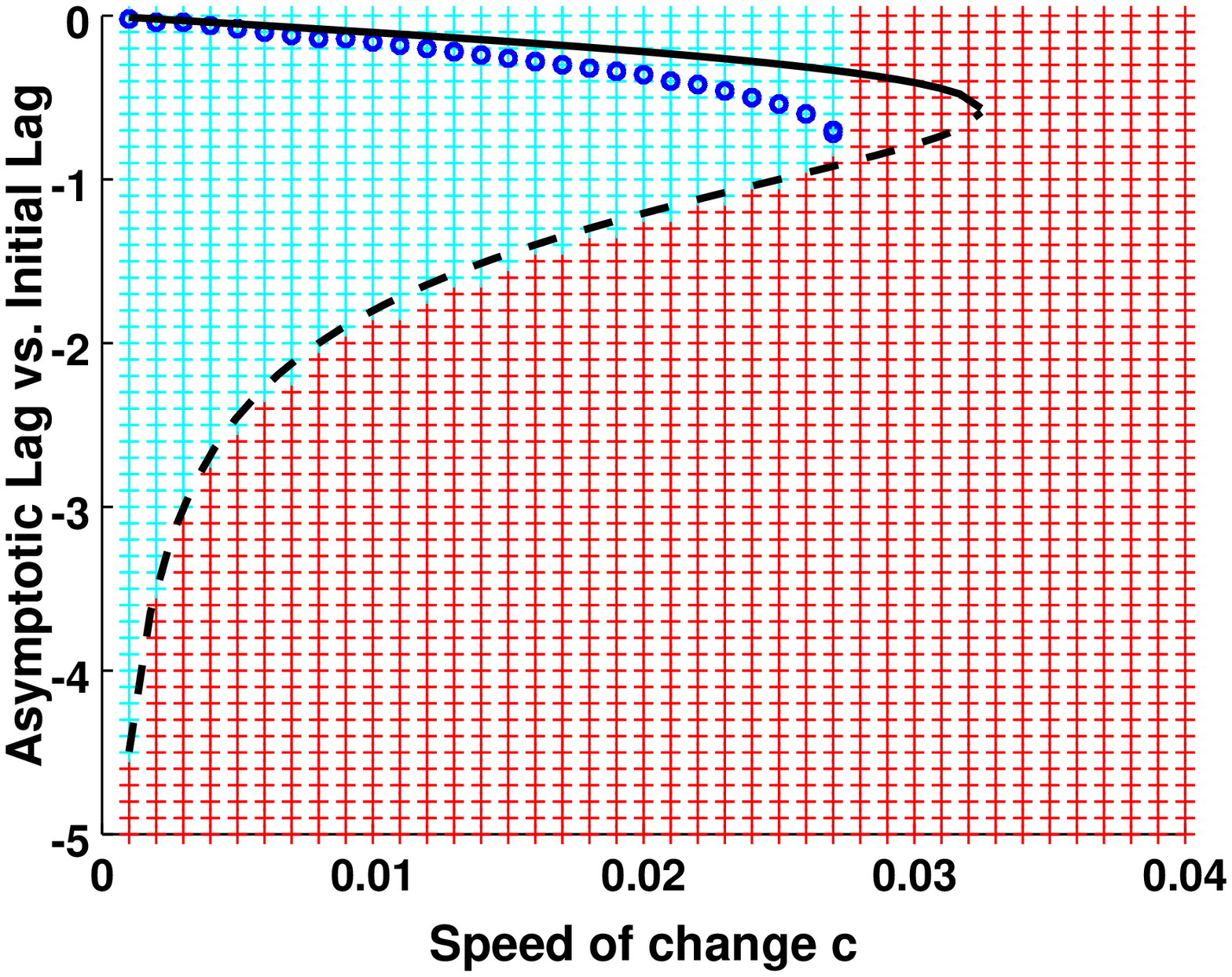} \label{fig:bistability}
}
\caption{Effect of the mean initial lag in the persistence of the population with various rates of environmental change $\cb$. We compute numerically, the solutions of the time--dependent problem~\eqref{eq:f z} with Gaussian initial conditions centered on various evolutionary lags $\zb^*_{init}$, depicted by crosses in the figures. We repeated this exploration for various speeds $\cb$ ranging in $(0,1.5\,\cb^*_{tip})$. For each case, we plot the evolutionary lag $z^*$ at the final time of computations (blue circles) and the analytical evolutionary lags given by the first line of Table~\ref{tab-summary}
(black lines): the plain lines corresponds to the stable trait ($z^*$ in  asexual model and $z^*_s$ in infinitesimal sexual model) while the dashed lines corresponds to the unstable trait $z^*_u$ occurring in the infinitesimal sexual model. The light blue crosses correspond to initial data such that the final evolutionary lag is finite while the red crosses correspond to  initial data such that the final evolutionary lag diverges. In the asexual simulations, the mutation kernel is Gaussian. 
}
\label{fig:m_bounded}
\end{figure}

\paragraph*{The lag increases with the speed of environmental change.} In both the asexual model and infinitesimal model, we recover the classic result that  the lag $|\zb_0^*|$  is an increasing function of $\cb$ (as illustrated by Fig.~\ref{fig:asex_kernel_C} and  Fig.~\ref{fig:m_c_selection}).  

In the asexual model, the evolutionary lag at equilibrium is such that the mortality rates equals $\betab L\left(\frac{\cb}{\sigmab\betab}\right)$ (see Table~\ref{tab-summary}). The latter quantity increases with the rate of environmental change. As the mortality rate $\mb$ increases when we move away from the optimal trait, the lag $|\zb_0^*|$  must also increase with respect to $\cb$. 

In the infinitesimal model of sexual reproduction,the evolutionary lag at equilibrium is found where the gradient of selection $( \mb')$ equals $- \frac{\cb}{\sigmab^2}$, which increases in absolute value with the rate of environmental change $\cb$ (see Table~\ref{tab-summary}).  In the convex neighborhood of the optimal trait, the gradient of selection $(\mb')$ is increasing with deviation from the optimum, hence the lag $|\zb_0^*|$ is increasing with respect to $\cb$. However, if the fitness function has both a convex and a concave part (as in the yellow curves in Fig.~\ref{fig:m_c_selection}),  there may be multiple equilibria fulfilling the condition in Table~\ref{tab-summary} (see Fig.~\ref{fig:formula_sex}). In the concave part of the fitness function, the selection gradient is decreasing when $c$ increases, and so would the lag (see dashed curve in Fig.~\ref{fig:bistability}). However, heuristic argument and numerical simulations suggest that equilibrium points in the concave part of the fitness function are unstable (see Fig.~\ref{fig:bistability} and more detailed discussion of this scenario below).

\paragraph*{The lag increases faster or slower than the speed of environmental change.}
Our analytical predictions suggest that a linear relationship between the rate of environmental change and the evolutionary lag is expected only under special circumstances. We indeed show that the rate of increase of the lag according to the speed of change $c$ crucially depends on the shape of the selection in both the infinitesimal and asexual models (Fig.~\ref{fig:m_c_selection}). 

In addition, in the asexual model, this rate of increase will crucially depend on the shape of the mutation kernel through the Lagrangian function $L$. Indeed, we can show from our formula in Table~\ref{tab-summary} that the lag increases linearly with the speed of change $c$ as soon as the function $c\mapsto m^{-1}(L(c))$ is linear.  Thus, both the shape of selection and that of the mutation kernel interact to determine how the evolutionary lag responds to faster environmental change. If the selection function is quadratic (i.e. $m(z)=z^2/2,$), we can show from the convexity of the Lagrangian function $L$ that the lag increases linearly with $c$  only in the diffusion approximation $L(c)=c^2/2$ (see Table~\ref{tab-summary-quadra} and blue curve in Fig.~\ref{fig:asex_kernel_C}), while it increases sub--linearly for any other mutation kernels (see red, orange, purple and green curves in Fig.~\ref{fig:asex_kernel_C}). We can further show that the lag in this scenario increases more slowly with the rate of environmental change when the kurtosis of the mutation kernel is higher (see Appendix~\ref{app:kernel_kurtosis} for mathematical details). In Fig.~\ref{fig:asex_kernel_C}, we compare four different mutation kernels with increasing kurtosis: uniform distribution kernel (red), Gaussian kernel (orange), double exponential kernel (purple) and Gamma kernel (green). In the asexual model, a fat tail of the mutation kernel thus tends to reduce the lag, even though this effect is most visible when the environment changes fast (Fig.~\ref{fig:asex_kernel_C}) .

To examine the effect of the selection function on how the evolutionary lag increases in faster changing environment, we now focus on the case of diffusion approximation in the asexual model ($L(c)=c^2/2$), for the sake of simplicity, and compare it to the results in the infinitesimal model.
In both cases, we can exhibit a simple criteria to decipher the nature of this increase. Let us first observe that, in those cases, the lag increases  {\em linearly} with  $c$ if the selection function is quadratic (see Table~\ref{tab-summary-quadra} and the blue curves in Fig.~\ref{fig:m_c_selection}). The lag however accelerates with $c$ if $m$ is \emph{sub-quadratic} in the following senses (see orange curves in Fig.~\ref{fig:m_c_selection}):
\begin{equation}\label{eq:weak_selec_zstar0}
\dfrac{m'' m }{(m')^2} < \frac12 \quad \text{(asexual)}\,, \quad m{'''} > 0\quad \text{(infinitesimal sexual)}\, .
\end{equation} 
Conversely, the lag decelerates with $c$ if $m$ is \emph{super-quadratic} in the following senses (see red curves in Fig.~\ref{fig:m_c_selection}): 
\begin{equation}\label{eq:strong_selec_zstar0}
\dfrac{m'' m }{(m')^2} > \frac12 \quad \text{(asexual)}\,, \quad m{'''} < 0\quad \text{(infinitesimal sexual)}\, .
\end{equation}
The criteria are of different nature depending on the model of reproduction (asexual versus infinitesimal). However, they coincide in the case of a homogeneous selection function $m(z) = |z|^p$  $(p>1)$. Indeed,  selection is super-quadratic in both cases if and only if $p>2$.  More generally, the lag is reduced when the selection function has a stronger convexity  in the sense of~\eqref{eq:strong_selec_zstar0}. This behavior is illustrated in Fig.~\ref{fig:m_c_selection}.

\paragraph*{The lag can diverge for a faster speed of environmental change.}
As observed by ~\cite{OsmKla17}, we also find that the lag may diverge, i.e. grow infinite for some finite threshold in the speed of environmental change, when the selection function is too weak away from the optimum. Such "evolutionary tipping point" are predicted, both for the infinitesimal model and in the asexual model, for some shapes of the selection function. The underlying mechanisms are however qualitatively different in the two models, as explained below.

In order to illustrate this phenomenon, we consider a bounded selection function depicted in Fig.~\ref{fig:m_c_selection} (orange curve).
We find the following critical speed $\cb^*_{tip}$, 
\begin{equation}
\begin{array}{ll}
\ds \cb^*_{tip} = \sqrt{2\sigmab^2\betab \left ( \max_{\zb \in(-\infty,0)}  \mb (\zb) \right ) } 
&  \text{(asexual)} \\[3mm]   
\ds \cb_{tip}^* = \sigmab^2\left(  \max _{\zb \in(-\infty,0)} |\mb'(\zb)|\right )
& \text{(infinitesimal sexual)}
\end{array},
\end{equation}
so that the lag is finite if and only if $\cb < \cb^*_{tip}$, while the lag diverges if  $\cb > \cb^*_{tip}$ and the population cannot keep pace with the environmental change. The difference between the two formulas can be understood through graphical arguments (see Fig.~\ref{fig:formul_m_bounded}). In the asexual model, the lag at equilibrium is found where the mortality rate equals a specific value, which increases with the speed of change $\cb$. This point is found where the selection function intersects an horizontal line, of higher elevation as $\cb$ increases  in Fig.~\ref{fig:formul_m_bounded}. With a bounded mortality function, there is thus a finite value of $\cb$ for which this critical quantity equals the maximal mortality rate, the latter being reached for an infinitely large lag. In the infinitesimal model, the equilibrium lag is found where the selection gradient equals a specific value increasing with $\cb$. Graphically, this point is found where the local slope of the selection function equals such a critical value. With a bounded mortality function such as in Fig.~\ref{fig:formul_m_bounded}, there are in general two equilibrium points characterized by such local slope, one stable in the convex part and one unstable in the concave part. As the speed of environmental change increases, so does the local slope at the two equilibria, which gradually converge towards the inflection point of the mortality function with the maximal local slope. This point characterizes the maximal speed of environmental change for which there is a finite evolutionary lag. Above that critical speed of change, the lag grows without limit. 
We illustrate this phenomenon of severe maladaptation in Fig.~\ref{fig:m_c_selection} (see the orange curves).

Despite the existence of tipping points in both cases, the transition from moderate ($\cb<\cb^*_{tip}$) to severe maladaptation ($\cb>\cb^*_{tip}$) have different signatures depending on the reproduction model.
In the asexual model, the lag becomes arbitrarily large as the speed $\cb$ becomes close to the maximal sustainable speed $\cb^*_{tip}$. Conversely, in the infinitesimal model, the lag remains bounded by the value of the inflexion point until it drops to infinity when the speed becomes larger than $\cb^*_{tip}$. 

We can also see a major difference between the two reproduction models when we look at the time dynamics (Fig.~\ref{fig:m_bounded}). We run simulations of the density dynamics described by  the equation~\eqref{eq:f z} and 
starting from various initial data centered at different traits (see crosses in Fig.~\ref{fig:m_bounded}). In the infinitesimal model, when the initial lag is beyond the unstable point $z^*_u$, defined in Fig.~\ref{fig:formul_m_bounded}(b), the lag diverges, whereas it converges to the stable point $z^*_s$, also defined in Fig.~\ref{fig:formul_m_bounded}(b), if the lag is initially moderate. We see that the long term adaptation of the population to a changing environment does not only depend on the speed of change, but also on the initial state of the population. 
In the asexual model, the initial configuration of the population does not play a significant role in the long term dynamics of adaptation: we observe that the population can adapt whatever the initial lag is, if the speed of change is below $\cb^*_{tip}$ (see Fig.~\ref{fig:m_bounded}).
We can expect such difference because the lag at equilibrium is uniquely defined in the asexual model while it can take multiple values in the infinitesimal model if the function has an inflexion point, which is the case for the bounded selection function (see Fig~\ref{fig:formul_m_bounded}).

\subsection{The mean fitness}\label{sec:fitness}
 We now investigate the effect of the changing environment on the mean fitness of the population.

 \paragraph*{The mean fitness decreases with increasing speed of environmental change $\cb$.} In both scenarios the \emph{lag load} $\triangle \lambdab$, defined as the difference between the mean fitness without changing environment ($\lambdab_0=\betab-\mub_0$) and the mean fitness under changing environment, is (unsurprisingly) given by the increment of mortality at the evolutionary lag $\mb(\zb_0^*)$
 \[
    \triangle\lambdab = \betab-\mub_0- \lambdab = \mb(\zb^*_0).
 \]
 Since $\mb$ is symmetric increasing and the lag $|\zb_0^*|$ is increasing with respect to $\cb$, we deduce that the mean fitness decreases with respect to $\cb.$ It is illustrated in Fig.~\ref{fig:m_c_selection} for different selection functions. 
 
 In the asexual model, the lag-load takes the following form 
 \[
    \triangle\lambdab = \betab L\left(\frac{\cb}{\sigmab\betab}\right). 
 \]
 which is exactly the expression~\eqref{eq:L} in the original units with a speed $\cb$.
 Since $L$ increases with the kurtosis of the mutation kernel, we deduce that higher kurtosis of the mutation kernel increases the mean fitness (see Fig.~\ref{fig:asex_kernel_C}). 
 
\begin{figure}[h!]
  \makebox[0.45\linewidth][c]{Asexual model} \hfil
  \makebox[0.45\linewidth][c]{Infinitesimal sexual model} \\[3mm]
  \subfigure[]{\includegraphics[width = 0.45\linewidth]{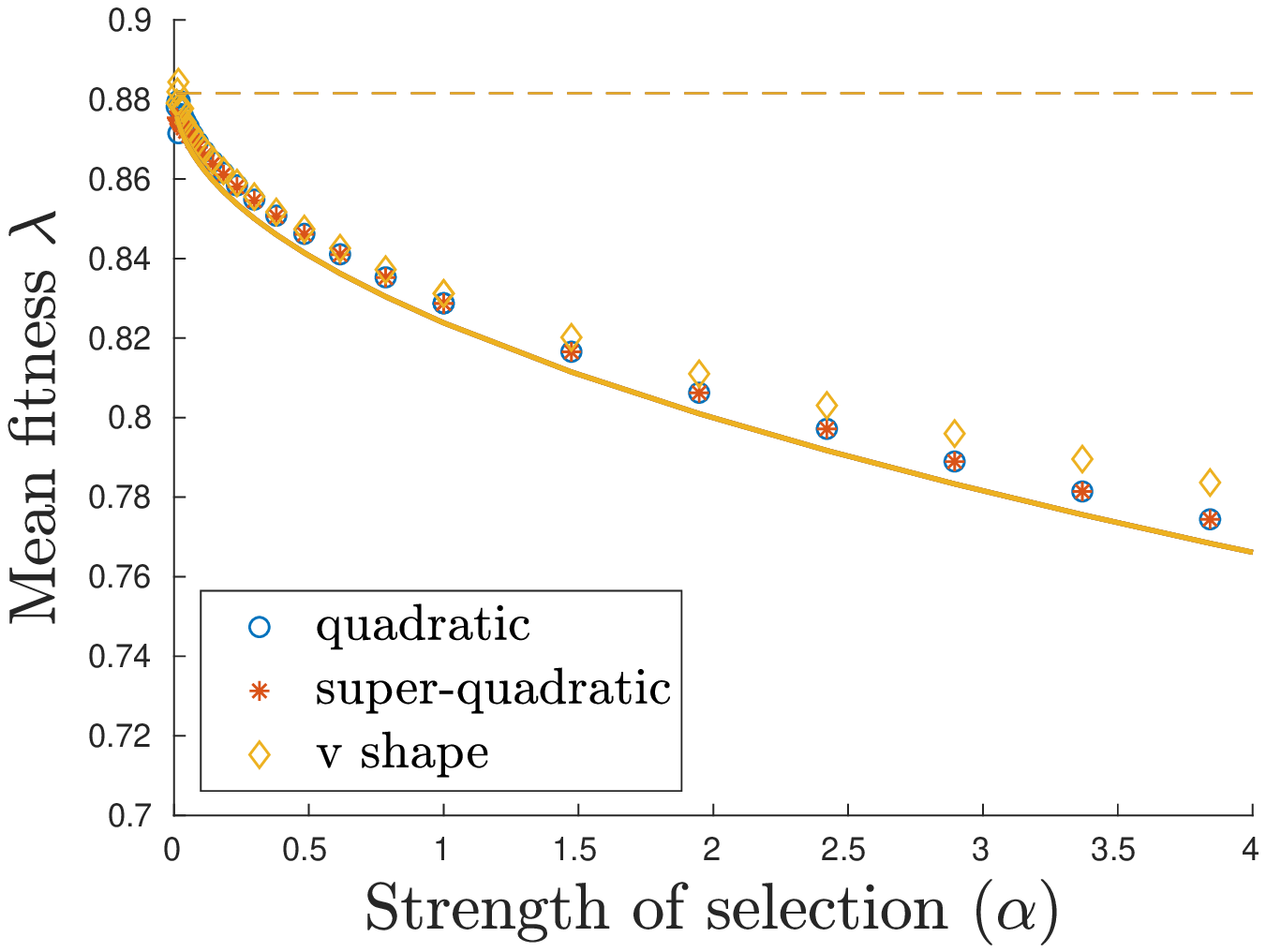}}
   \hspace{7mm}
   \subfigure[]{\includegraphics[width = 0.45\linewidth]{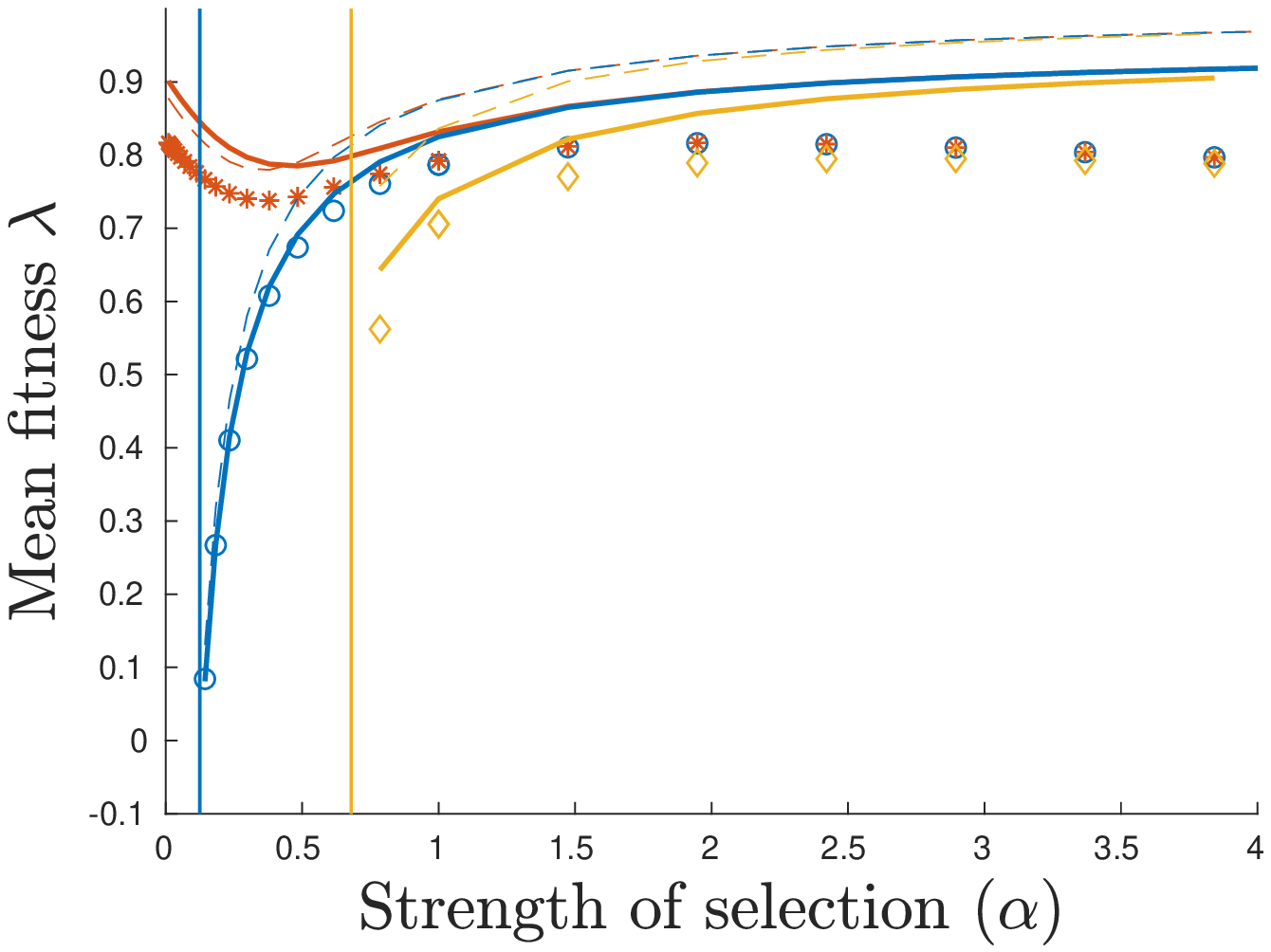}} 
   
   \subfigure[]{\includegraphics[width = 0.45\linewidth]{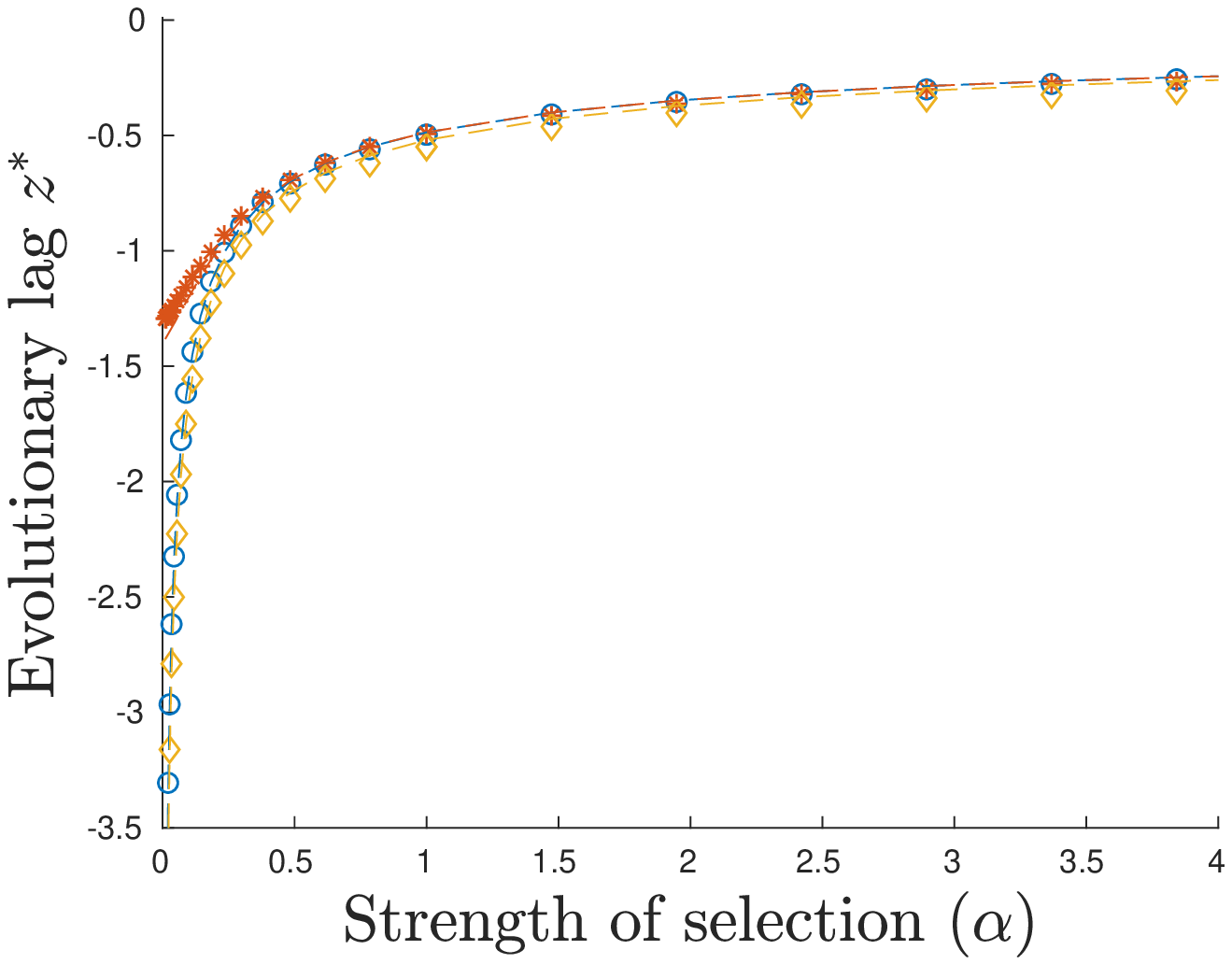}}
   \hspace{7mm}
   \subfigure[]{\includegraphics[width = 0.45\linewidth]{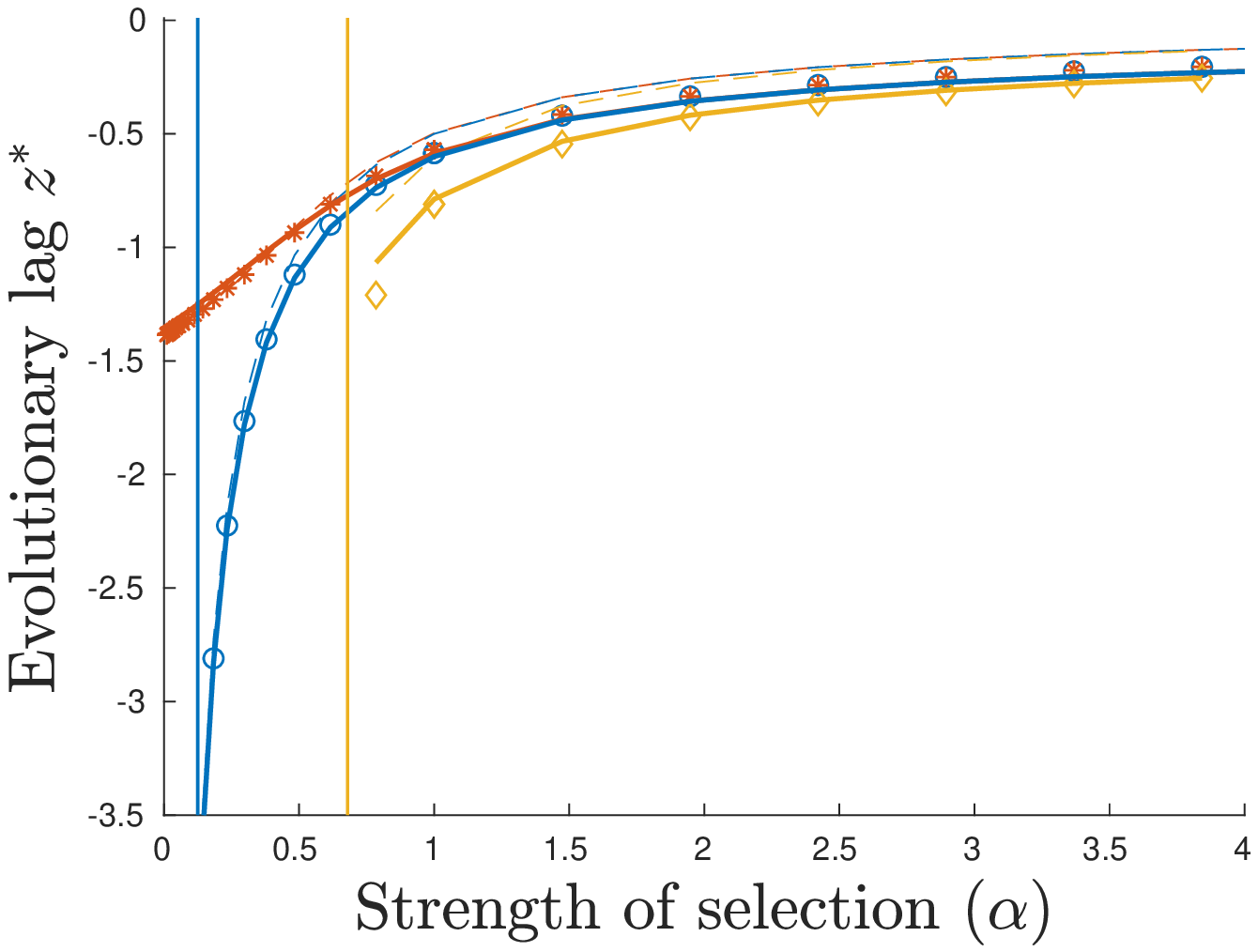}} 
   
   \subfigure[]{\includegraphics[width = 0.45\linewidth]{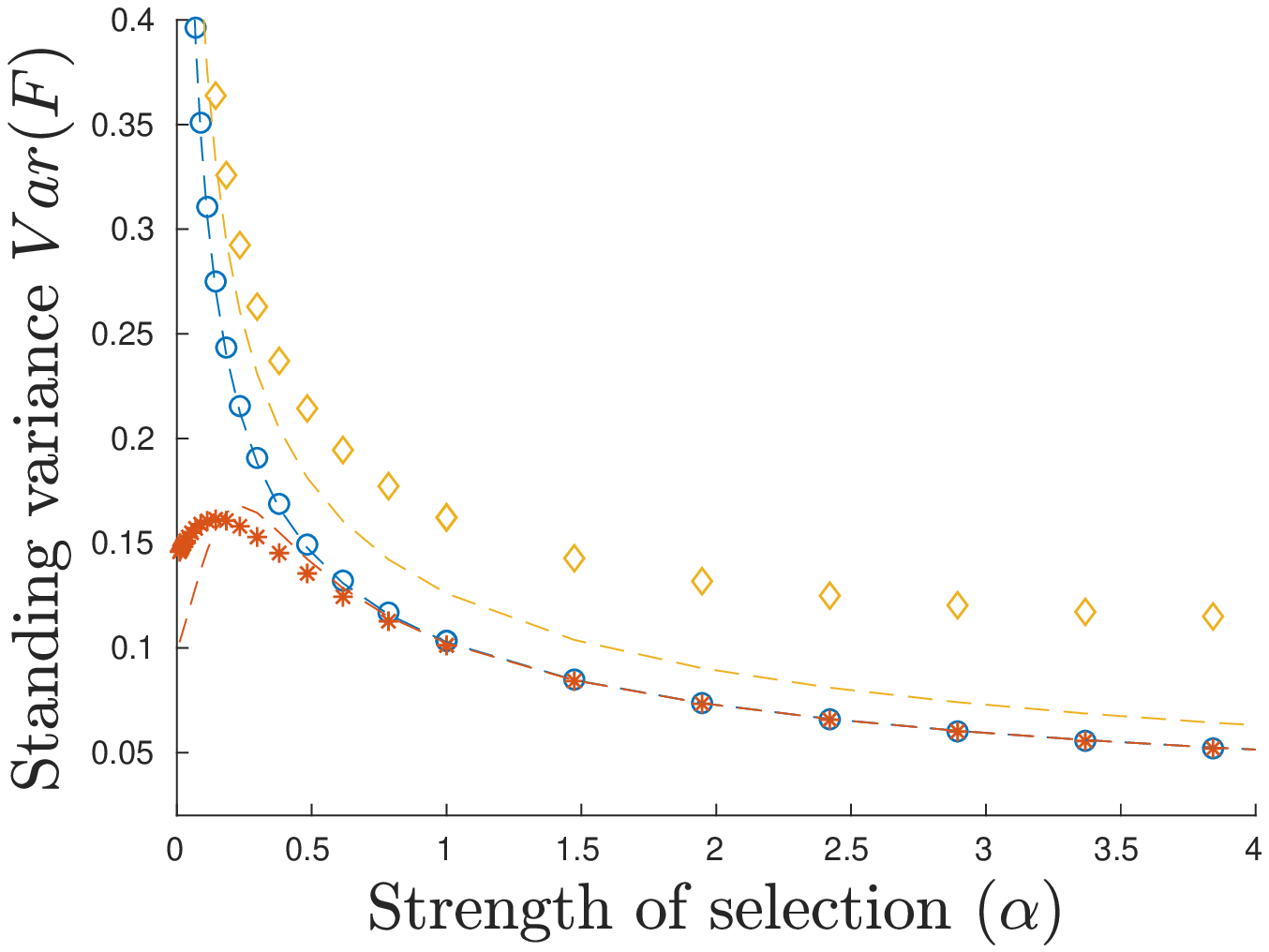}}
   \hspace{7mm}
   \subfigure[]{\includegraphics[width = 0.45\linewidth]{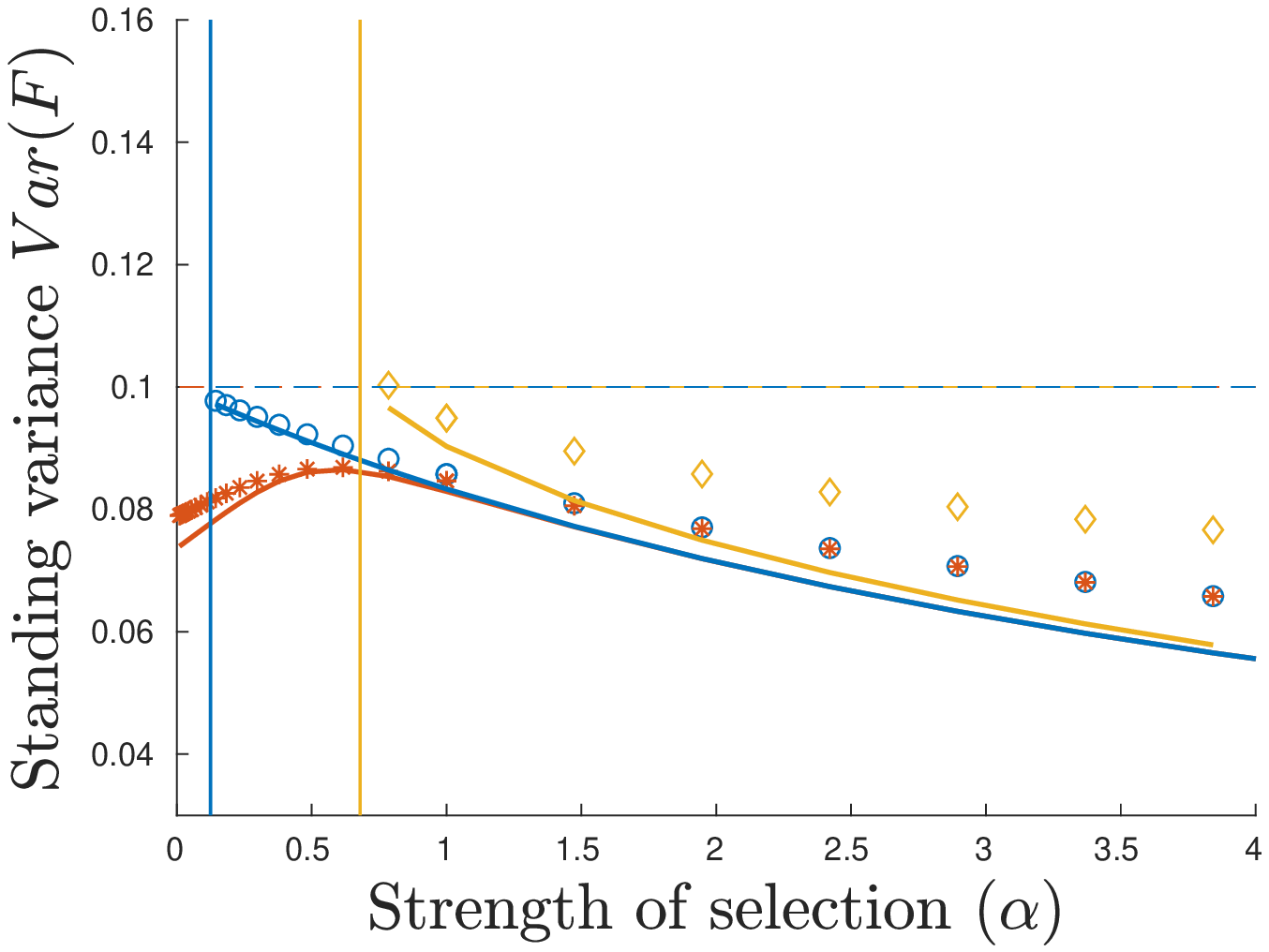}} 
   \caption{Influence of the strength of selection $\alphab$ on the mean fitness $\lambda$, the evolutionary lag $z^*$ and the standing phenotypic variance $\Var(F)$ at equilibrium in an environment changing at rate $\cb=0.05$  and with three different selection patterns: quadratic (blue curves), super--quadratic (red curves) or bounded (orange curves). Other parameters are: $\betab = 1$, $\sigmab = 0.1$ and the  intensity of selection $\alphab$ ranges from $10^{-2}$ to  $4$. 
   We compare our analytical results (first approximation dashed lines and second approximation plain lines) with the numerical simulations of the stationary distribution of \eqref{eq:f z} (marked symbol) for both asexual and sexual infinitesimal model. In the asexual model, we only consider a Gaussian mutation kernel. 
   }\label{fig:alpha_m}
\end{figure}

\paragraph*{The shape of selection affects the lag load in the infinitesimal model, but not in the asexual model.}
In the asexual model, the lag load only depends, at the leading order, on the speed of environmental change and the mutation kernel through the Lagrangian function $L$  (see Table~\ref{tab-summary}). In particular, for a given speed $\cb$, a fertility rate $\betab$ and a given mutation kernel, we predict that the lag load is constant (see dashed line in Fig.~\ref{fig:alpha_m}(a)). At the next order of approximation, the mean fitness however depends on the local shape of the selection function around the optimal trait: $\alphab = \mb''(0)$. The mean fitness is then predicted to decline as the strength of stabilizing selection around the optimum $\alphab$ increases (see Fig.~\ref{fig:alpha_m}(a)), due to increasing standing load.
These predictions are confirmed by our  numerical simulation see Fig.~\ref{fig:m_c_selection}(a) and~\ref{fig:alpha_m}(a).

In contrast, the influence of the selection pattern is more intricate in the case of the infinitesimal model of reproduction. The lag load depends strongly on the global shape of $\mb$ (see Fig.~\ref{fig:m_c_selection}(b) and~\ref{fig:alpha_m}(b)). In particular, we see that for low strength of selection $\alphab$, the mean fitness crucially depends on the shape of selection. Mean fitness is higher in the scenario with super--quadratic selection than with quadratic selection, and lowest when selection is sub-quadratic in Fig.~\ref{fig:m_c_selection}(b) and~\ref{fig:alpha_m}(b)). Moreover,  the mean fitness increases with increasing strength of selection in the quadratic case, while it initially decreases for the super--quadratic case. However, for stronger strength of selection,  the shape of selection has less importance. Our approximation allows us to capture those differences. For instance, in the quadratic case (blue curves in Fig.~\ref{fig:alpha_m} and~\ref{fig:m_c_selection}), we can see from Table~\ref{tab-summary-quadra} that the mean fitness increases with the strength of selection $\alphab$ at the leading order,  which corresponds to small value of $\alphab$. However, when the strength of selection becomes stronger, antagonistic effects occur at the 
next order so that the fitness may decrease due to standing load, defined in~\eqref{eq:lambda_1}~\citep{LynLan93,LanSha96,KopMat14}. This effect is illustrated in Fig.~\ref{fig:alpha_m}(b).

\subsection{The standing variance}\label{sec:variance}
 In both asexual diffusion approximation and the infinitesimal model, the standing variance does not depend on the speed of change $\cb$ when the selection function is quadratic (see blue curves in Fig.~\ref{fig:m_c_selection}).
The standing variance however increases with $\cb$ if the selection function is \emph{sub-quadratic} in the sense of~\eqref{eq:weak_selec_zstar0} (see orange curves in Fig.~\ref{fig:m_c_selection}). Conversely, the standing variance decreases with $\cb$ if the selection function is \emph{super-quadratic} in the sense of~\eqref{eq:strong_selec_zstar0} (see red curves in Fig.~\ref{fig:m_c_selection}) -- see details in appendix~\ref{app:Var_c}.

 The standing variance is less variable in the infinitesimal model than in the asexual model. It was expected from our analysis (see formula of Table~\ref{tab-summary}) because the infinitesimal model tends to constrain the variance of the phenotypic distribution. Indeed, we know from previous analysis~\citep{MirRao13,BarEthVeb17}, that in the absence of selection, the infinitesimal model generates a Gaussian equilibrium distribution with variance $\sigmab^2$. Our analysis shows that under the small variance assumption, the standing variance is close to this variance $\sigmab^2$ and our numerical analysis shows that standing variance slowly deviates from the genetic variance without selection $\sigmab^2$, when either the speed of change increases or the strength of selection increases. This pattern is observed whatever the shape of selection. We can thus conclude that for the infinitesimal model under the small variance hypothesis, the standing variance is not very sensitive to either 
selection (strength of selection or shape of selection) or the speed of environmental change.

Conversely, in the asexual model, the standing variance is quite sensitive to the selection function. This is emphasized in the case of a bounded selection function. The standing variance dramatically increases as the speed of change becomes close to the critical speed $\cb^*_{tip}$ because  the selection gradient becomes flat (see Table~\ref{tab-summary}).

 In the asexual model, the standing variance is moreover sensitive to the shape of the mutation kernel. We see from Fig.~\ref{fig:asex_kernel_C}(c) that the standing variance generally increases with a fatter tail of the mutation kernel. There are however exceptions to this pattern  (see for instance the Gamma mutation kernel at low speed of environmental change, green curves in Fig.~\ref{fig:asex_kernel_C}). This situation, unexpected by our approximation, might be due to the fact that when the speed of change is low, the mutations with large effects are quickly eliminated by selection, which in turn reduces the standing variance.

\subsection{Persistence of the population: the critical speed $\cb^*$ }\label{sec:persitence}
The final outcome of our analysis is to compute the speed $\cb^*$ beyond which the population cannot keep pace with the environmental change ($\lambdab <0$). In the general case, we can obtain the following approximation formula:
\begin{equation}
   \begin{cases}
      \cb^* = \sigmab\betab\, L^{-1}\left( \dfrac{\betab - \mub_0}{\betab} \right)       & \quad \text{(asexual model)}\\[3mm]
      \cb^* = \sigmab^2\mb'\left(\mb^{-1}(\betab-\mub_0)\right) & \quad \text{(infinitesimal model)}
   \end{cases}
\end{equation}
We can first observe that, in the small variance regime, the critical speed in the asexual model does not depend on the shape of the selection $\mb$, but  on the mutation kernel through the Lagrangian $L$. Thus, for any selection function, the critical speed is the same. Conversely, for the infinitesimal model, the critical speed crucially depends on the shape of the selection $\mb$.
Moreover, we can mention that the discussion of the dependency of $\lambdab$ with respect to various parameters also holds naturally for $\cb^{*}$.

When we consider the diffusion approximation for the asexual model ($L(v)=v^2/2$) and the quadratic selection function $\mb(\zb) = \alphab\zb^2/2$, we obtain the following formula:
\begin{equation}
\begin{cases}
\cb^*  = \sqrt{2}  \sigmab \sqrt{\betab}   \left( \betab-\mub_0 \textcolor{gray}{- \dfrac{\sigmab(\alphab\betab)^{1/2}}{2}} \right)^{1/2} & \quad \text{(asexual model)}\\
\cb^*  = \sqrt{2}  \dfrac{\sigmab^2}{ \left(1 \textcolor{gray}{ + 4 \sigmab^2\dfrac{\alphab}{\betab}} \right)^{1/2}}
\sqrt{\alphab} 
         \left ( \betab-\mub_0 \textcolor{gray}{- \dfrac{   \sigmab^2\alphab}2 }
            \right )^{1/2}              
            & \quad \text{(infinitesimal model)}
\end{cases}\label{eq:critical speed}
\end{equation}

The formula \eqref{eq:critical speed} is in agreement with previous results where it was assumed that the phenotypic lag $\zb$ is normally distributed in the population, which corresponds in our framework to assuming that the equilibrium distribution $\Fb$ is Gaussian (see for instance Eq. [A6] in~\citep{KopMat14}). There, the formula is given with the standing phenotypic variance as a parameter by: 
\begin{equation}
\cb^* \approx  \Var(\Fb) \sqrt{2\alphab\lambdab(0)}\, .\label{eq:kopp}
\end{equation}
where $\lambdab(0)$ corresponds to the mean fitness in absence of environmental change ($\cb=0$), and it is given by formula of Table~\ref{tab-summary-quadra} with $\cb=0$.
This is perfectly consistent with the formula for the variance obtained in both scenarios. However, the formulation \eqref{eq:kopp} might be misleading, as it omits some possible compensation, such as the selection strength $\alphab$, which disappears in the case of asexual reproduction because it also affects $\Var(\Fb)$.

\subsection{Numerical predictions for the whole distribution of phenotypes}\label{sec:asex_num}

\paragraph*{Quality of approximation.}
For the asexual model, we only compare the simulation results with our first order approximation stated in Table~\ref{tab-summary} (black colored), except for the variation of the mean fitness with respect to the strength of selection, where we need to take into account the standing load that appears at the second order of approximation (see gray colored formula in Table~\ref{tab-summary}). We can first observe from Fig.~\ref{fig:alpha_m} that our first  approximations are accurate when $\epsilon=\sigmab\sqrt{\alphab/\betab}$ is small (see value of $\alphab<0.5$ in Fig.~\ref{fig:alpha_m}). The scale of Fig.~\ref{fig:alpha_m}(a) is of order $\epsilon$, which is why the first order approximation seems less accurate than the second order approximation. This was expected since the standing load, which increases with the strength of selection, occurs at the second order of approximation. 
The approximation of $\zb^*$ and $\lambdab$ remain efficient even when $\epsilon$ increases (see Fig.~\ref{fig:asex_kernel_C} and~\ref{fig:m_c_selection} for small value of $\cb$). However, we see that the approximations deviate from the simulations when the speed of change increases and reaches the critical value $\cb^*$ (see Fig.~\ref{fig:asex_kernel_C} and~\ref{fig:m_c_selection}) or when the mutation kernel becomes leptokurtic (see green curves of Fig.~\ref{fig:asex_kernel_C}). 
The approximation of the standing variance is more sensitive to the parameter $\epsilon$. When $\cb$ and $\epsilon$ are small it is accurate (see Fig.~\ref{fig:asex_kernel_C}). However, when the speed increases, the approximation diverges from the simulations even if $\epsilon$ is small (see Fig.~\ref{fig:asex_kernel_C} and~\ref{fig:m_c_selection}).

For the infinitesimal model, we have compared our simulations to our first order approximation, as well as the second order approximation stated in Table~\ref{tab-summary} (first order approximation is black colored and second order approximation is gray colored). The first order approximation of $\zb^*$ and $\lambdab$ are efficient only when $\epsilon$ is really small, while the first order approximation of the standing variance may deviate from the simulation value even for small $\epsilon$ (see red curve Fig.~\ref{fig:alpha_m}(f)). However, the second approximations are really precise for small value of $\epsilon$ (see Fig.~\ref{fig:alpha_m}) and they remain efficient when $\epsilon$ increases and $\cb$ increases (see Fig.~\ref{fig:alpha_m} and~\ref{fig:m_c_selection}).

\paragraph*{Comparing simulations to the approximation for the entire distribution.}

\begin{figure}[h!]
     \makebox[0.45\linewidth][c]{Asexual model} \hfil
  \makebox[0.45\linewidth][c]{Infinitesimal sexual model} \\[3mm]
    \subfigure[Quadratic selection]{\includegraphics[width=0.45\linewidth]{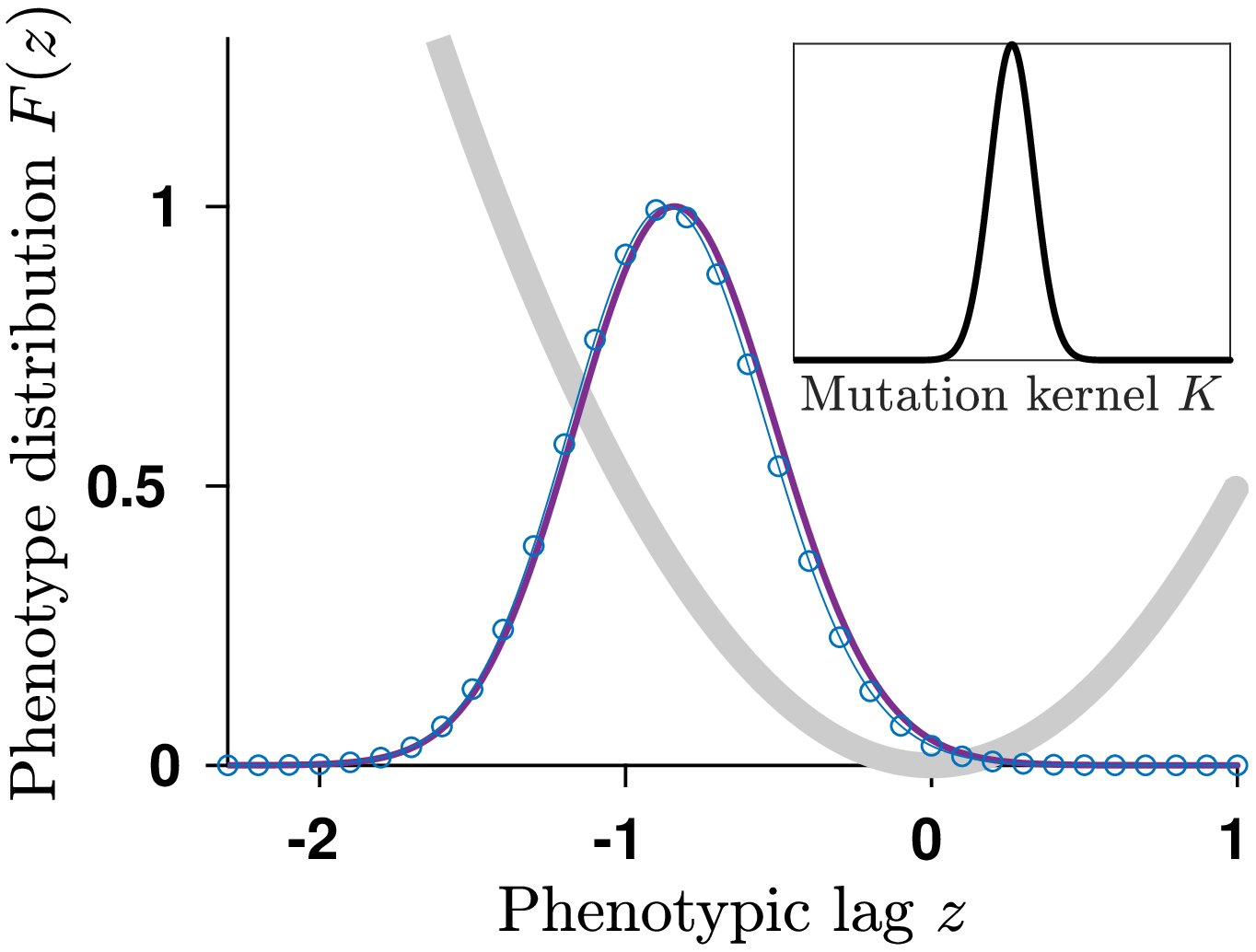}}
    \hspace{7mm}
    \subfigure[Quadratic selection]{\includegraphics[width=0.45\linewidth]{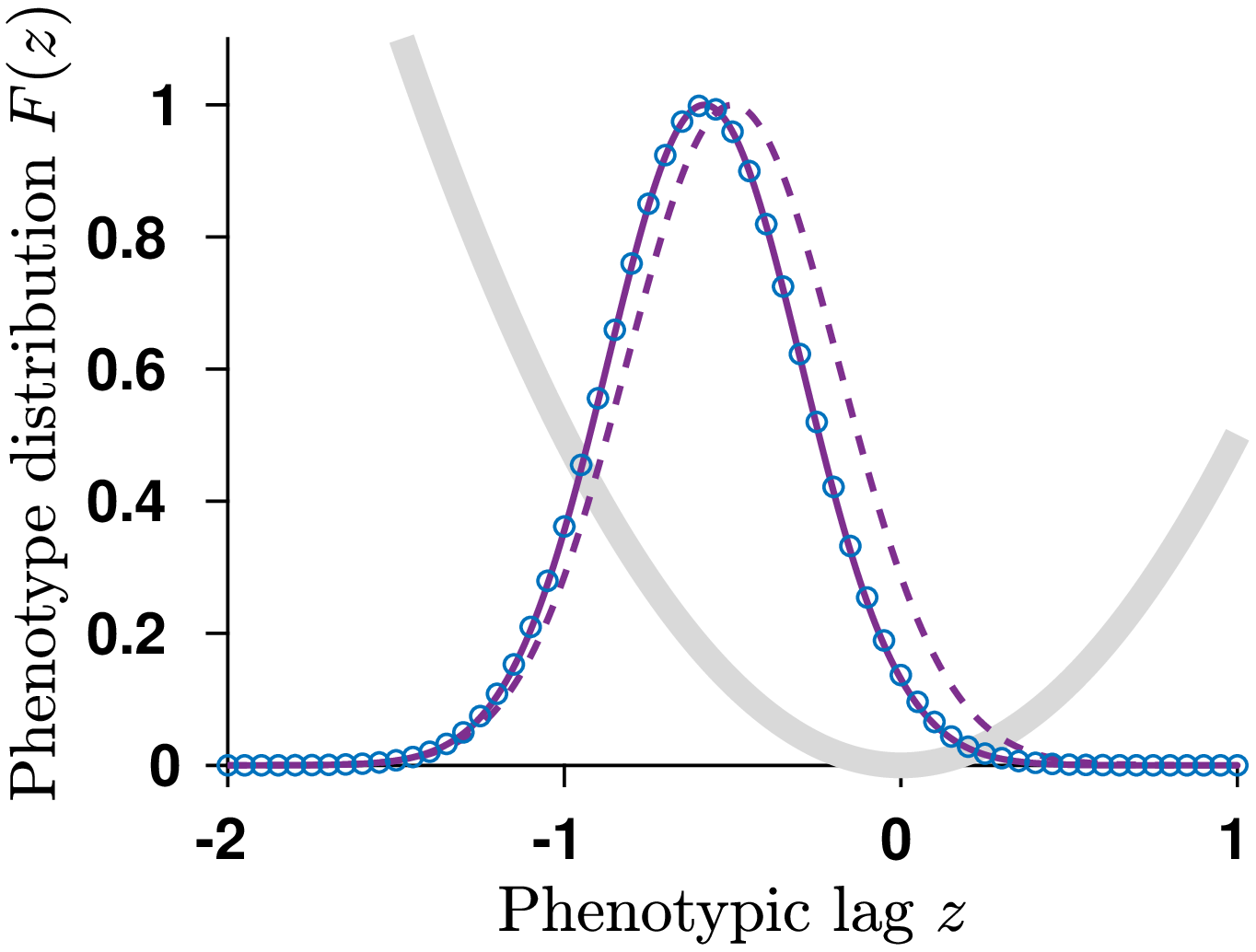}}
    
    \subfigure[Super--quadratic selection]{\includegraphics[width=0.45\linewidth]{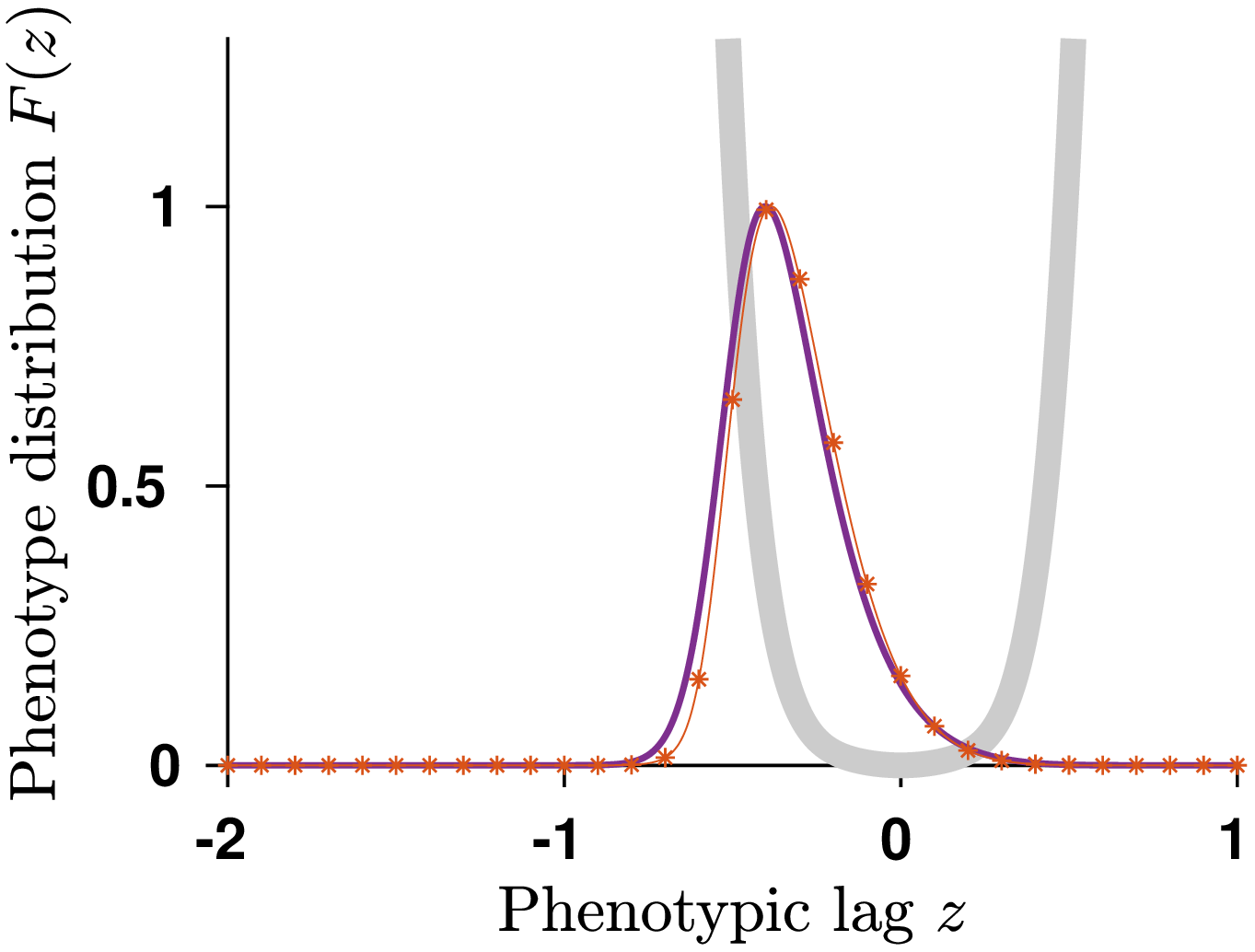}}
    \hspace{7mm}
    \subfigure[Super--quadratic selection]{\includegraphics[width=0.45\linewidth]{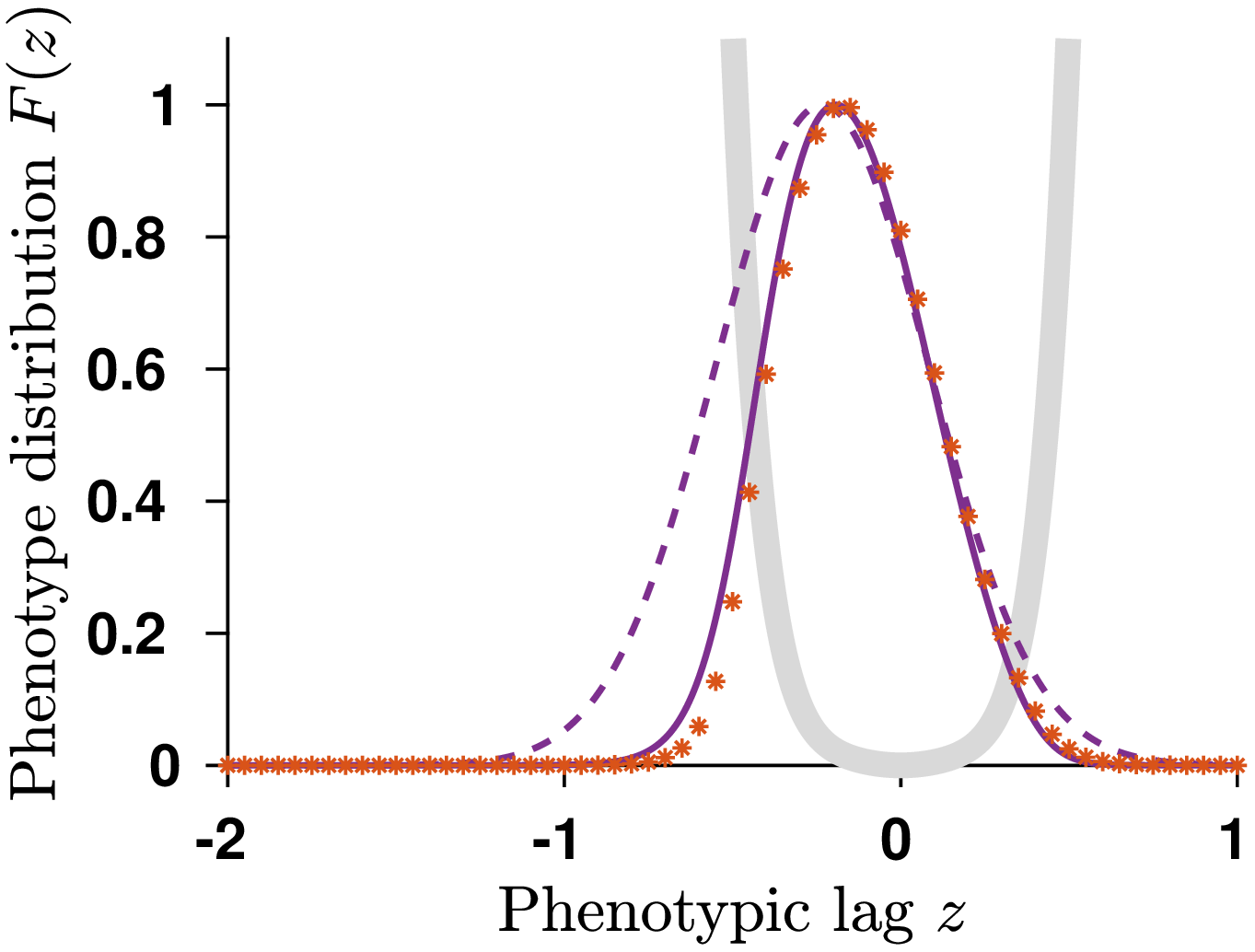}}
    
    \subfigure[Bounded selection]{\includegraphics[width=0.45\linewidth]{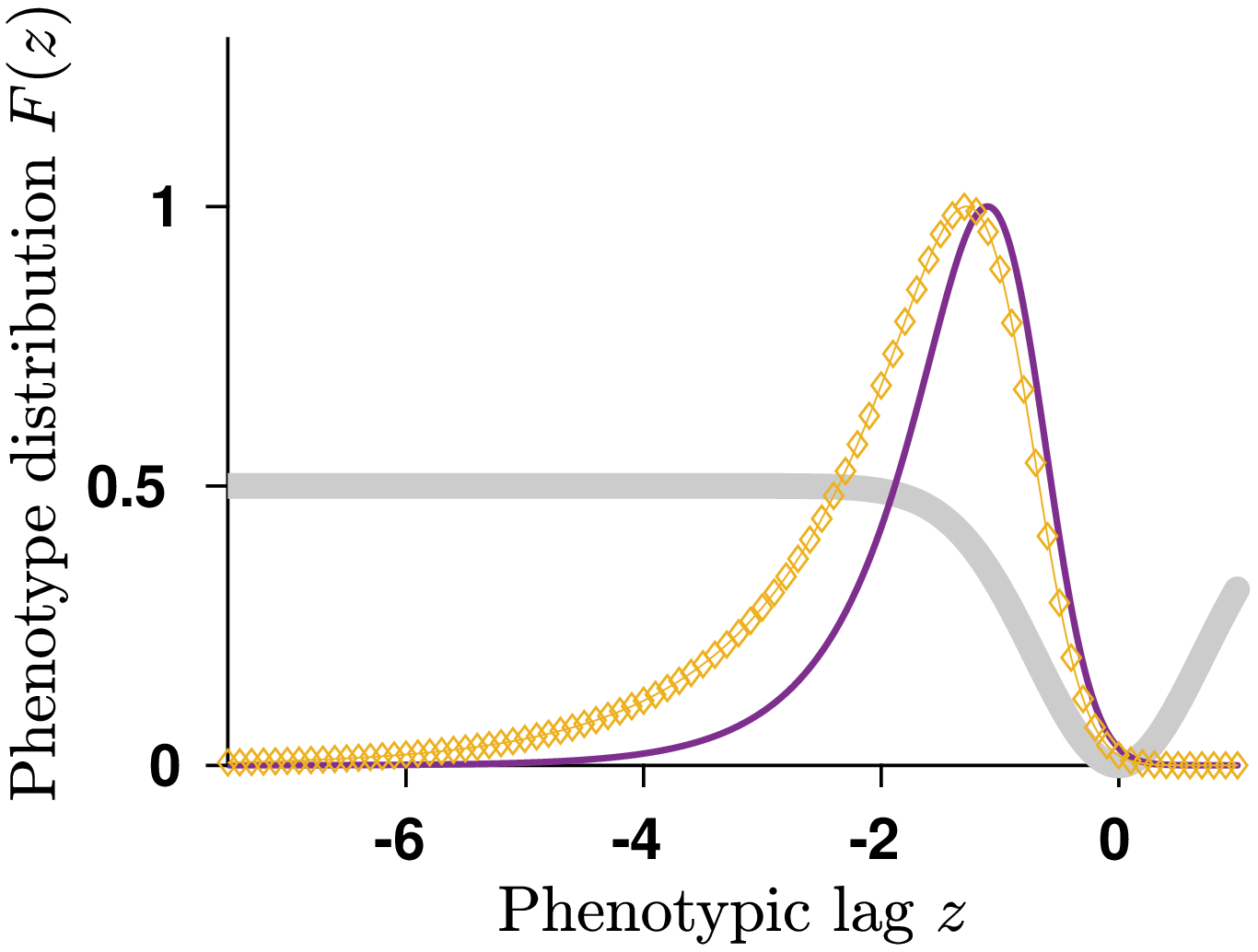}}
    \hspace{7mm}
    \subfigure[Bounded selection]{\includegraphics[width=0.45\linewidth]{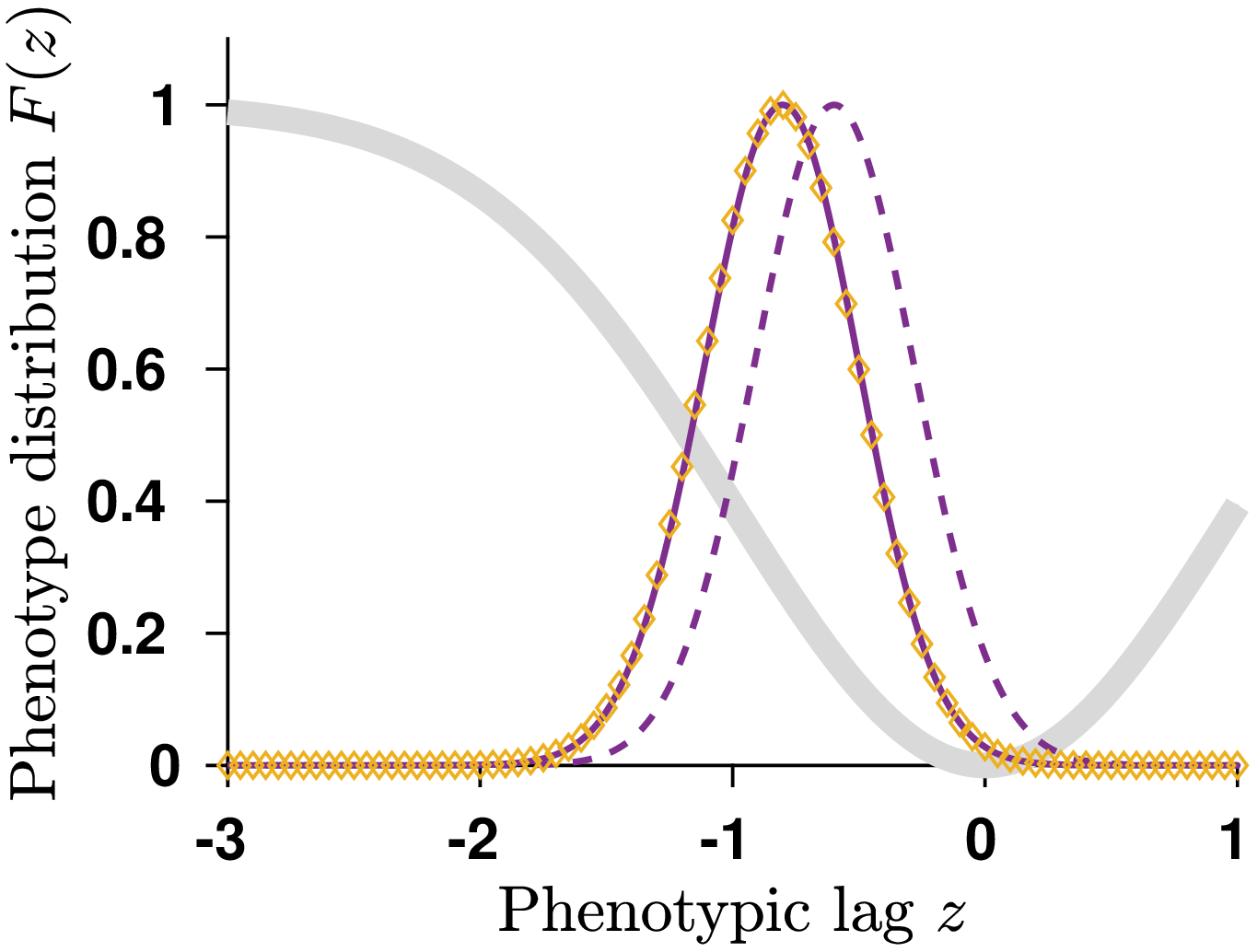}}
    \caption{Mutation-selection equilibria $\Fb$ for a speed of environment change of $\cb=0.09$ in the asexual model and $\cb=0.05$ in the infinitesimal sexual model, with different selection functions: (a)-(b) quadratic selection $\mb(z)=\alphab z^2/2$ (blue circled marked curves); (c)-(d) super-quadratic selection $\mb(z) = \alphab(z^2/2+64z^6)$ (blue star marked curves); (e)-(f) bounded selection function $\mb(\zb) =m_\infty( 1-\exp(-\alphab\zb^2/(2m_\infty))$ (orange dimaond marked curves). Other parameters are: $\alphab = 1$, $\betab = 1$, $\sigmab = 0.1$ and $m_\infty=0.5$ in the asexual model and $m_\infty=1$ in the infinitesimal sexual model. We compare simulated equilibria distribution $\Fb$ (marked curves) with our analytical results (first order results dashed curves and second order results plain curves). For the asexual scenario, we used the Gaussian kernel.}\label{fig:distribution}
 \end{figure}

We compare the simulated equilibrium distribution $\Fb$ with our analytical approximations (Fig.~\ref{fig:distribution}): the first order approximation corresponds to $\Fb_0 = \exp(-U_0/\eps^\gamma)$, where $U_0$ satisfies respectively the differential equation~\eqref{eq:HJ asexual U0} (asexual model) or $U_0(z) = (z-z^*_0)^2/2$ (infinitesimal sexual model), and $\gamma$ is respectively equal to $1$ in the asexual model and $2$ in the infinitesimal case; and the second order approximation $\Fb_1 = \exp(-U_0/\eps^\gamma - U_1)$, where $U_1$ satisfies respectively equation~\eqref{eq_app:asex_U1} (asexual model) or the non--local functional equation~\eqref{eq:U1 main} (infinitesimal model). Our simulations are performed with an $\epsilon^\gamma=0.1$, which is not that small.

In the asexual model, we can observe that the first order analytical approximation is really efficient at tracking the shape of the entire distribution for both super-quadratic and quadratic selection, even if $\epsilon$ is not that small (Fig.~\ref{fig:distribution}). For the bounded selection, our first order approximation fails to fit the left tail of the distribution, mainly because the speed of environmental change is close to the critical speed.

In the infinitesimal model, we can observe that the first order Gaussian approximation is not precise enough to track the entire distribution (Fig.~\ref{fig:distribution}). We need the second order approximation to fit the distribution. This is a direct consequence of our analysis, where we observe that we need the second order approximation to define  the first order approximation of the lag $z^*_0$.

\paragraph*{What me normal?} 
To go further in understanding the effect of a changing environment, we now look at the skewness and the kurtosis of the distributions. Those two indicators allow us to test whether the  distribution $\Fb$ can be well approximated by the Gaussian distribution.

\begin{figure}[h!]
     \makebox[0.45\linewidth][c]{Asexual model} \hfil
  \makebox[0.45\linewidth][c]{Infinitesimal sexual model} \\[3mm]
   \subfigure[]{\includegraphics[width = 0.45\linewidth]{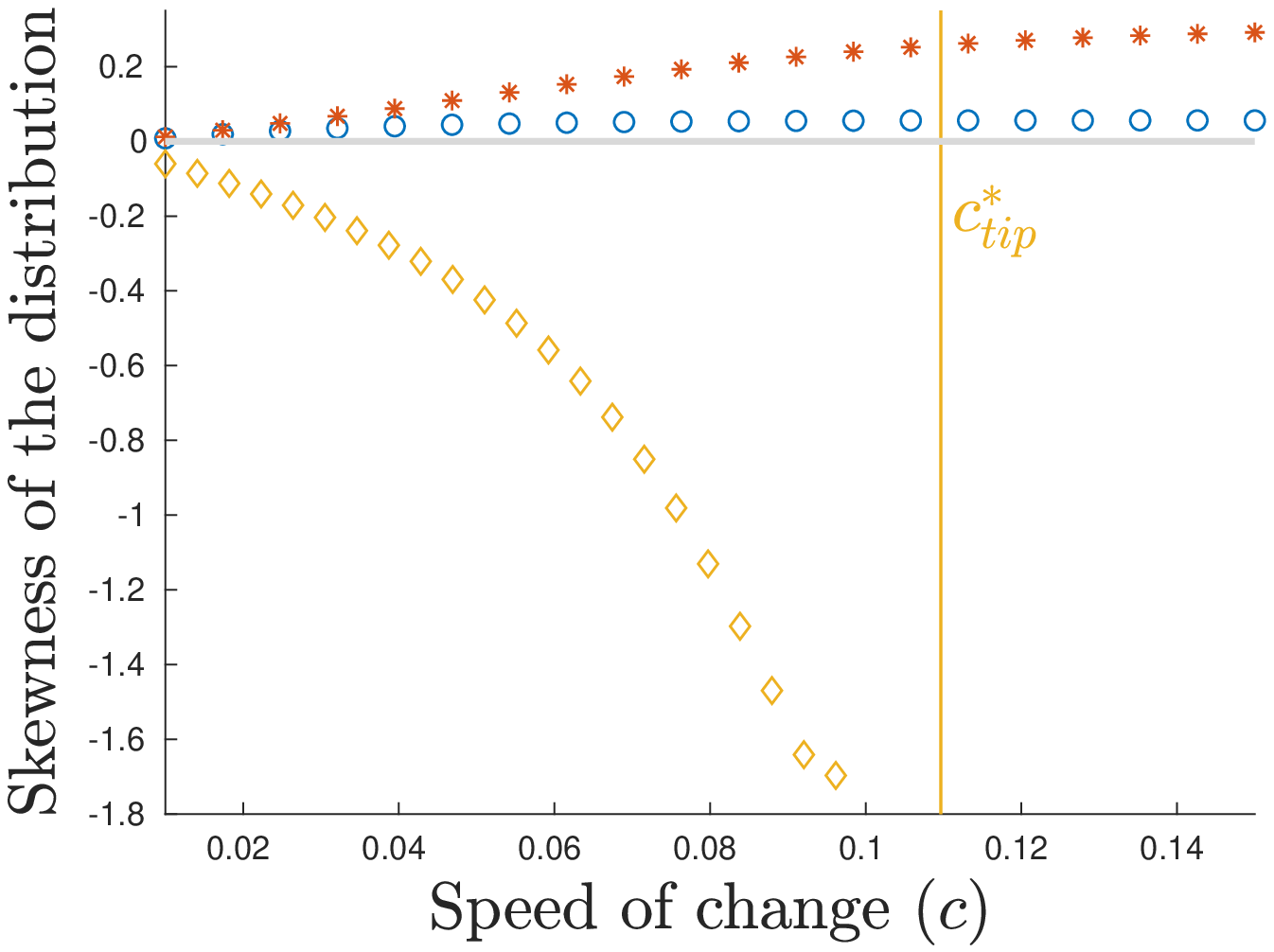}}
   \hspace{7mm}
   \subfigure[]{\includegraphics[width = 0.45\linewidth]{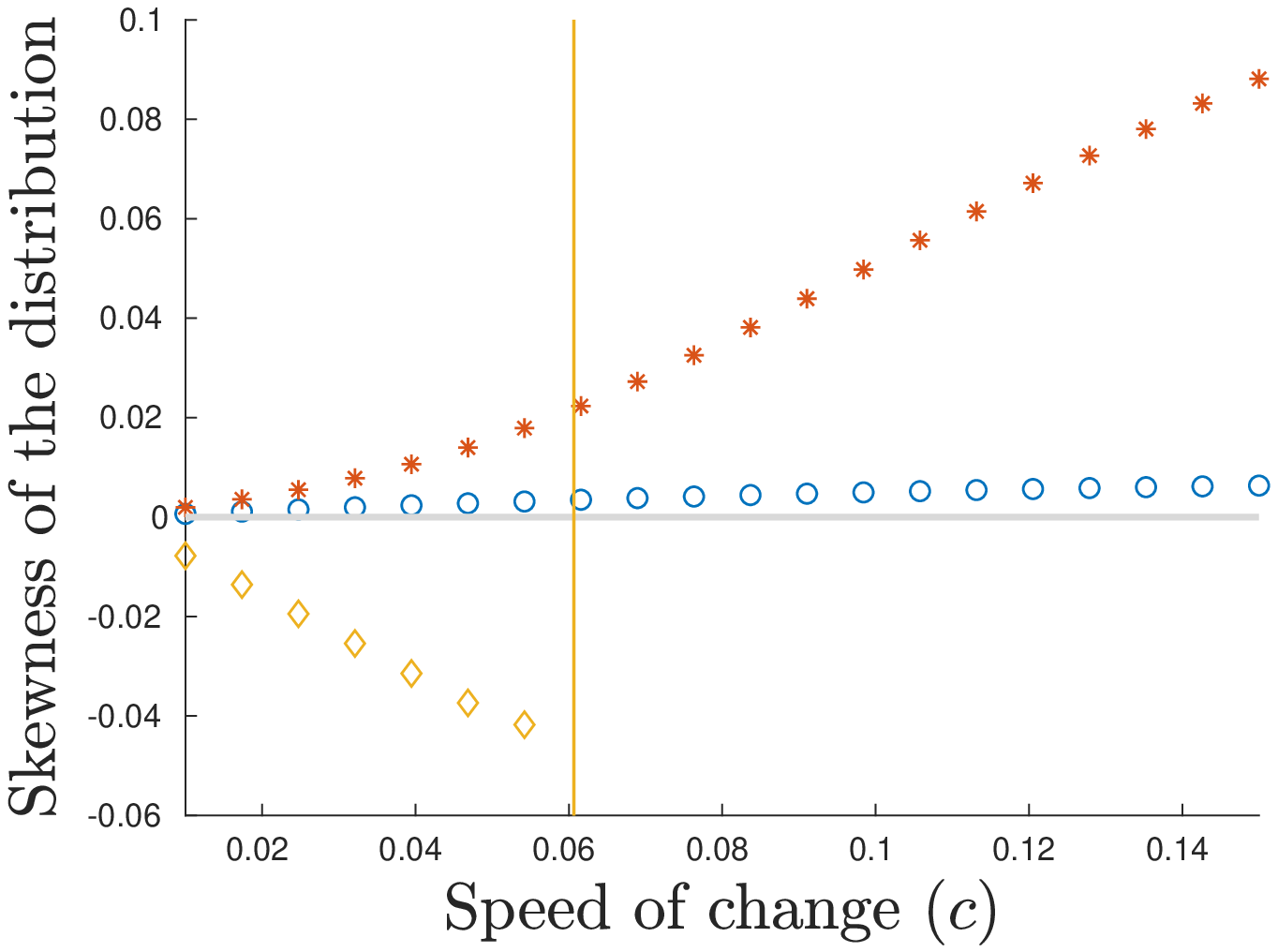}} 
   
   \subfigure[]{\includegraphics[width = 0.45\linewidth]{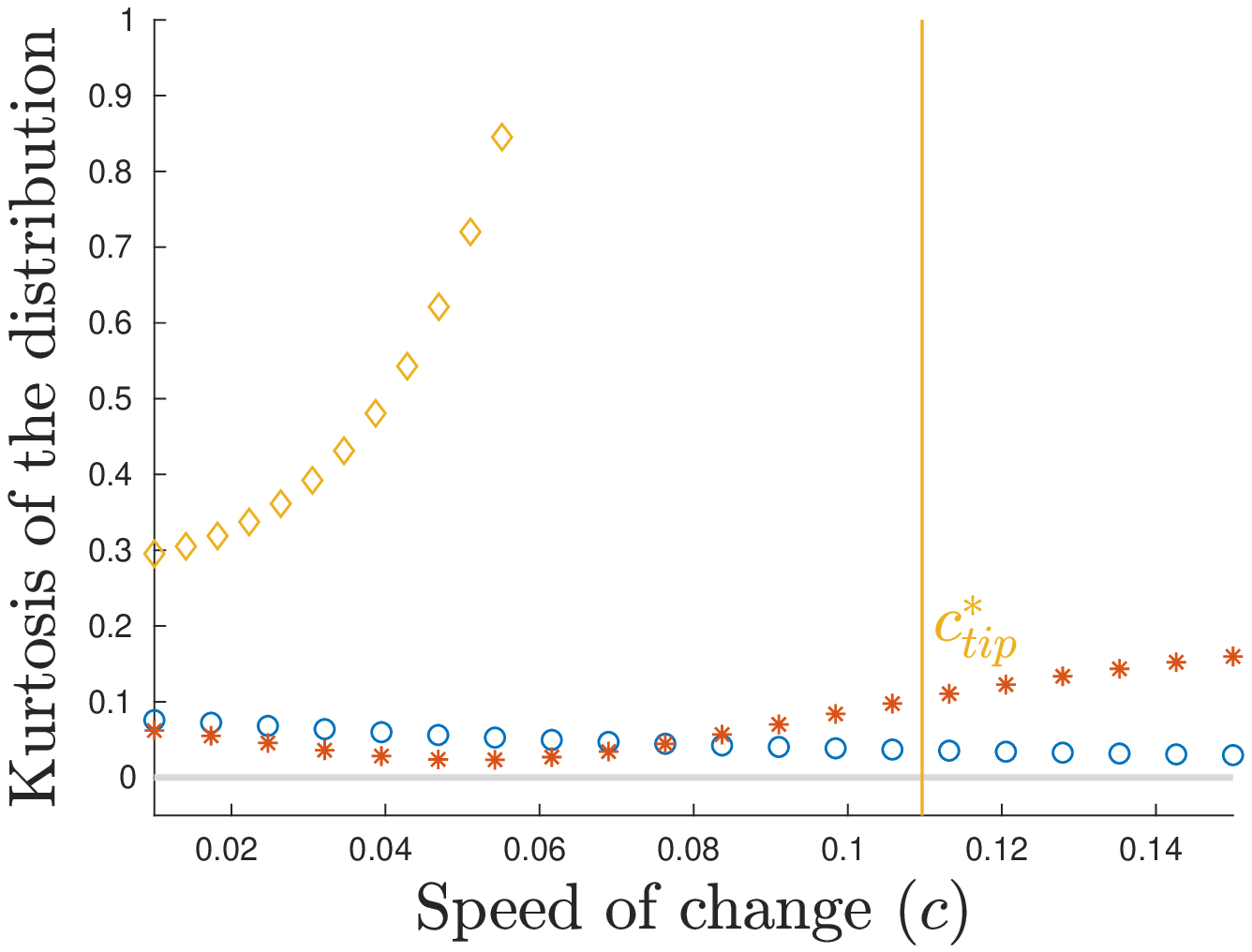}}
   \hspace{7mm}
   \subfigure[]{\includegraphics[width = 0.45\linewidth]{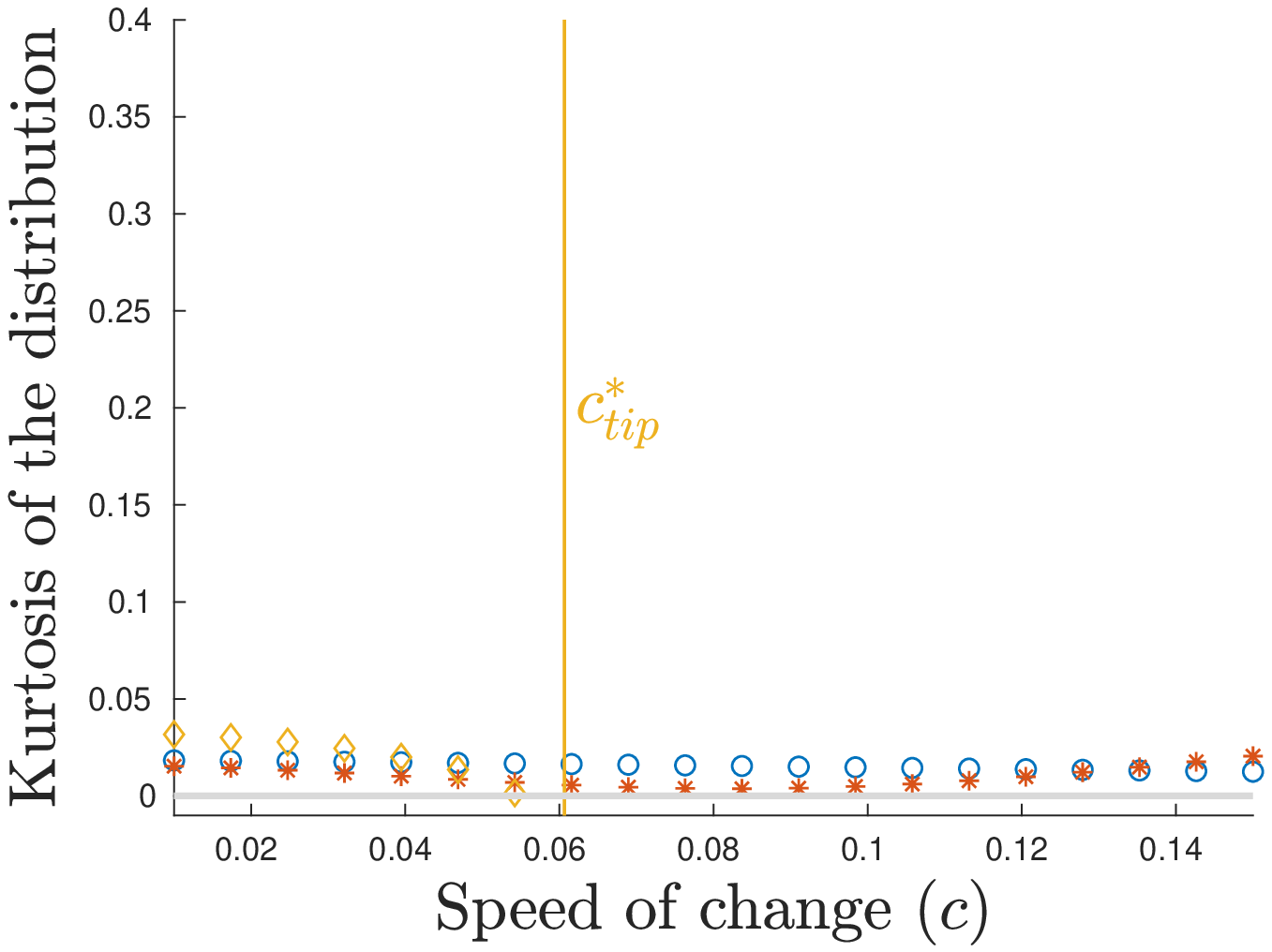}} 
   \caption{Influence of the speed of environmental change $\cb$ for three different selection function $\mb$: quadratic function $\mb(\zb) = \alphab\zb^2/2$ (blue curves), super--quadratic $\mb(\zb) = \alphab\zb^2/2 + z^6/64$ (red curves) or bounded $\mb(\zb) =m_\infty( 1-\exp(-\alphab\zb^2/(2m_\infty))$ (orange curves). Other parameters are: $\alphab = 1$, $\betab = 1$, $\sigmab = 0.1$ and $m_\infty=0.5$ in the asexual model and $m_\infty=1$ in the infinitesimal sexual model.
   In the asexual model, the mutation kernel is Gaussian.
   We compare our analytical results (dashed lines) with the numerical simulations of the stationary distribution of \eqref{eq:f z} (marked symbols) for both asexual and sexual infinitesimal model. It appears that our analytical results are able to catch interesting features even for relatively large speed of change $\cb$. }\label{fig:skewness_kurtosis}
\end{figure}

In the asexual model, we can observe from Fig.~\ref{fig:skewness_kurtosis} that, even for the quadratic selection, the distributions differ from a Gaussian distribution:  they are skewed and leptokurtic, which means that their kurtosis are higher than the kurtosis of the Gaussian distribution with same mean and variance. So the Gaussian distribution fails to track the exact distribution of the trait around the evolutionary lag of the population in a changing environment.   This phenomenon is enhanced when the selection function differs from the quadratic function (see Fig.~\ref{fig:skewness_kurtosis} diamond curves and ~\ref{fig:distribution}).  In addition, we see that, when the selection function is super--quadratic, the distribution has a positive skew, while, for a bounded selection function, it has a negative skew. 

Conversely, in the infinitesimal case, the Gaussian distribution well approximates the equilibrium distribution in general. This was already described by our approximation formula~\eqref{eq:ansatz sexual} in the section~\ref{sec:sex}. We can see that the kurtosis of the equilibrium distribution remains close to zero for any speeds of change and any selection functions.
 However, when the selection function is either super quadratic or bounded, we can observe from Fig.~\ref{fig:distribution} and~\ref{fig:skewness_kurtosis} that the distribution of phenotypes in the infinitesimal model also becomes skewed as the speed increases. The skew of the distribution corresponds to regions where the gradient of selection is low, with the same pattern as in the asexual model.  
 
 \newpage
 
\section{Discussion}\label{sec:discussion}
We have pushed further a recent methodology aimed at describing the dynamics of quantitative genetics models in the regime of small variance, without any {\em a priori} knowledge on the shape of the phenotype distribution. This methodology combines an appropriate rescaling of the equation with Taylor expansions on the logarithmic distribution.

{ Our approach differs from the previous studies based on the cumulant generating function (CGF), which is the logarithm of the Laplace transform of the trait distribution, here $C(t,p) = \log\left( \int e^{pz} f(t,z)\, dz \right)$. In his pioneering work,~\cite{Bur91}  derived equations for the so-called cumulants, which are the coefficients of the Taylor series of the CGF $C(t,p)$ at $p = 0$. However this system of equations is not closed, as the cumulants influence each other in cascade. This analysis was revisited in~\citep{MarRoq16} in the asexual model, using PDE methods. They derived an analytical formula for the CGF itself, but restricted it to a directional selection, when the trait represents the fitness itself. This was further extended to a moving optimum in~\citep{RoqPat20}. However, they made the crucial assumption of the Fisher Geometric Model for selection, which is analogous to our quadratic case, and diffusion for mutations, for which it is known that Gaussian distributions are particular solutions. The common feature with  our present methodology  is the PDE framework. Nevertheless, we focus our analysis on the logarithm of the trait distribution itself, as it is commonly done in theoretical physics to reformulate the wavefunction in terms of its action (see Appendix~\ref{app:conjugacy} for heuristics on this approach). This strategy is well-suited to { provide precise approximations}  with respect to a small parameter, for instance the wavelength in wave propagation (geometric optics) and the Planck constant in quantum mechanics (semi-classical analysis), and the standing variance in our theoretical biology setting.
}





{Here, the small variance regime corresponds either to relatively small effect of mutation, or to weak stabilizing selection. Under this regime, very little variance in fitness is introduced in the population through either
mutation or recombination events during reproduction.
However, the variance in fitness in the population can be relatively large depending on the shape of the selection. This regime can differ from the weak selection approximation classically used in quantitative genetics theory, corresponding to small variation in fitness in the population.}


{Under the small variance regime, we could describe analytically the phenotype distribution (see Table \ref{fig:distribution}), and assess the possible deviation from the Gaussian shape}. We further gave analytical approximations of the three main descriptors of the steady state: the evolutionary lag, the mean fitness, and the standing phenotypic variance (see Table \ref{tab-summary}).

Noticeably,  two  different models of reproduction, assuming either asexual reproduction, or infinitesimal sexual reproduction with an infinite number of freely recombining loci (the infinitesimal model),  could be handled in a unified framework. This allows discussing similarities and discrepancies between the two models, which are frequently used in analytical models of adaptation to changing and/or heterogeneous environments.

\paragraph*{Relaxing the Gaussian distribution assumption. }
Our analytical framework allows us to relax the assumption of a Gaussian distribution of phenotypic values, commonly made by many quantitative genetics models of adaptation to a changing environment with a moving optimum, both in the case of sexually \citep[e.g.][]{BurLyn95,OsmKla17}  and asexually reproducing organisms \citep[e.g.][]{LynGabWoo91}. Consistently with previous simulations and analytical results~\citep{TurBar94,Bur99,Jon12}, our results show that we expect stronger deviations from a Gaussian distribution of phenotypes if the selection function departs from a quadratic shape, if the mutation model departs from a simple diffusion, if reproduction is asexual rather than well described by the infinitesimal model, and/or if the environment changes relatively fast. We in particular recover the observation made by \cite{Jon12}  in their simulations that the skew of the phenotypic distribution is greater in absolute value in faster changing environments, but we further predict that the sign of this skew critically depends on the shape of the selection function away from the optimum, an observation that could not be made by their simulations that only considered quadratic selection.

\paragraph*{Universal relationships. }
Interestingly, despite deviations from the Gaussian distribution, our predictions in the regime of small variance for the evolutionary lag, or the critical rate of environmental change, are consistent with predictions of past quantitative genetics models that have assumed a constant phenotypic variance and a Gaussian distribution of phenotypes. We discuss below the links between the present results and those past predictions and how they provide new insights. 
\label{sec:2relations}
As a direct consequence of the small variance assumption, the two following relationships, linking the three main descriptors of the population (the evolutionary lag, mean fitness and phenotyoic variance), hold true, whatever the model of reproduction (either asexual or infinitesimal): 
\begin{equation}\label{eq:2relations}
\begin{cases}
\lambda \approx 1 - m(z^*)\smallskip\\
\Var(F) \approx  - \dfrac{\eps^\gamma c}{ m'(z^*)}
\end{cases}
\end{equation}
The first relationship corresponds to the demographic equilibrium, when the mean fitness is the balance between (constant) fecundity and mortality at the evolutionary lag. The second one corresponds to the evolutionary equilibrium, when the speed of evolutionary change (as predicted by the product of phenotypic variance and the selection gradient) equals the speed of change in the environment. Note that our model assumes for simplicity that the phenotypic variance is fully heritable. Those relationships are better visualized in adimensional units. They can  be deduced directly from equations~\eqref{eq:asexual no age}-\eqref{eq:sexual no age}. 
Although the reproduction model does not affect the demographic relationship, it  influences the evolution relationship through the scaling exponent $\gamma$ ($\gamma=1$ for asexual reproduction and $\gamma=2$ for infinitesimal sexual reproduction). 
Similar equations appear in quantitative genetics models assuming a Gaussian phenotypic distribution and a constant phenotypic variance.
In particular, with quadratic selection, the second relationship allows us to recover the following results of~\cite{BurLyn95} and ~\cite{KopMat14}:
\begin{equation}
   |\zb^*|\approx \dfrac{\cb}{\alphab \Var(\Fb)}.
\end{equation}

However, the relationships in \ref{eq:2relations} are not enough to compute the three descriptors, if one does not consider the standing phenotypic variance $\Var(F)$ as a fixed parameter, as previous studies often did. 
Our small variance approximations allows us to predict the value of the phenotypic variance in a changing environment in the two models, where previous studies have generally used simulations \cite[e.g.][]{Bur99} to examine how the evolution of the phenotypic variance affects the adaptation of sexual and asexual organisms in a changing environment. Many of our results are ultimately explained by the fact that the evolution of the phenotypic variance is under very different constraints under the asexual model and the infinitesimal model.

In the asexual model, the evolution of the phenotypic variance is not strongly constrained and has in particular no upper bound. 
The mean fitness $\lambda$ does not depend on the shape of the selection function at the leading order (see \eqref{eq:lambda0} and Table \ref{tab-summary}), but only on the speed of environmental change and on the mutation kernel. Once the mean fitness is determined, the evolutionary lag $z^*$ and the standing variance $\Var(F)$ are deduced from respectively the first and the second relationship in \eqref{eq:2relations}. The standing variance then strongly depends on the shape of the selection in the asexual model.
In contrast, in the sexual infinitesimal model, we found that the standing variance  $\Var(F)$ does not depend on the shape of the selection function at the leading order (see \eqref{eq:ansatz sexual} and Table \ref{tab-summary}). {The infinitesimal model sets an upper bound to the phenotypic variance.}
Once the {maximum} variance $\sigma$ is set, the evolutionary lag $z^*$ and the mean fitness $\lambda$ are deduced from respectively the second and the first relationship \eqref{eq:2relations}. Therefore, in the infinitesimal model, constraints on the variance determines the value of the lag, and thus the mean fitness, while in the asexual model, the evolution of the variance allows reaching a demographic equilibrium where the increased loss of fitness due to a changing environment is compensated by the gain of fitness due to beneficial mutations. In the asexual model, it is the mean fitness that determines the value of the lag and in turn the value of the phenotypic variance. Most of our predictions (discussed below) are a consequence of this core discrepancy between the two models. 

\paragraph*{Mean fitness weakly depends on selection in the asexual model, but not in the infinitesimal model.}

 
 In the asexual model, $\lambdab$   depends on $\mb$ only at the second order  through the second derivative around the optimal trait $\alphab =\mb''(0)$ \eqref{eq:scale_alpha}. Hence, up to a reasonable accuracy, the mean fitness depends (weakly) on the local shape of the selection pattern around the optimal trait, even if the population can be localized around an evolutionary lag far from the optimal trait. This happens because, in a gradually moving environment, the asexual population is constantly regenerated by the fittest individuals. This phenomena is apparent when tracing back lineages in the population at steady state: it was proven independently by \cite{PatForGar20} and \cite{CalHenMelTran21} that the typical trajectories of ancestors of individuals sampled uniformly in the population converge to the optimal trait backward in time. 
In contrast, the mean fitness strongly depends on the shape of the selection function in the infinitesimal sexual model. It appears clearly in the quadratic case where $\alphab$ enters into the formula for the mean fitness at the leading order (Table \ref{tab-summary-quadra}). In particular, we recover the previous finding that weak selection represents a "slippery slope" in a changing environment, leading to a lower mean fitness, when effects of selection on the evolution of standing variance are neglected~\citep{KopMat14}. Again, it is interesting to link this finding to the behavior of the typical trajectories of the ancestors in the infinitesimal model, which converge to the evolutionary lag backward in time \cite[Chapter 5]{Patout-PhD}. 

\paragraph*{The shape of selection has strong effects on the evolution of the evolutionary lag and phenotypic variance under both the asexual and infinitesimal models.}
In both models, however, the exact shape of the selection function away from the optimum has noticeable consequences for the evolution of the lag between the  {the mean phenotype} in the population and the optimal moving phenotypic value, and for the evolution of the standing variance, especially in fast changing environments. There is unfortunately very scarce empirical evidence about the exact shape of fitness landscapes and how much they deviate from a quadratic, due to the difficulty to estimate precisely the shape of such fitness functions (see however the predictions of \cite{Gau20} suggesting strong deviations from a quadratic shape in the case of a trait involved in climate adaptation). Most models of adaptation to a moving optimum assume, for mathematical convenience and in the absence of strong empirical support for an alternative, a quadratic selection function. Our analysis allows considering a broad diversity of selection functions and also to draw general conclusions about how their shape may affect the evolution of the phenotypic distribution. In both asexual and infinitesimal models, we found, consistently with previous predictions~\cite[reviewed in][]{KopMat14}, that the lag increases with the speed of environmental change: however there is a linear relationship between the two only when assuming a quadratic selection function. When the selection function is super-quadratic (and selection much stronger away from the optimum), this puts a brake on maladaptation and the evolutionary lag does not increase as fast when the environment changes more rapidly. For the same reason, the phenotypic variance then declines when the environment changes faster in the super-quadratic selection scenarios. Conversely, with a sub-quadratic selection function, the weakening of selection away from the optimum results in larger lags, accelerating maladaptation with increasing speed of environmental change and increasing phenotypic variance. There has been little discussion yet in the theoretical literature of the consequences of the exact shape of selection in changing environments (see however \citep{OsmKla17,Kla20} and discussion of tipping-points below). {In a constant or stationary environment with weak fluctuations, the mean phenotype value is never very far from the optimum and the quadratic selection is an adequate approximation. However, the present results suggest that further empirical investigation of the shape of the fitness landscape far from the optimum is critically needed to understand how much populations may depart from the optimal phenotypic value.}

\paragraph*{Evolutionary tipping points.}
The case of sub-quadratic selection functions has recently attracted some interest, since it was discovered that the weakening of selection away from the optimum could lead to evolutionary tipping points: above some critical speed of environmental change, the evolutionary lag grows without limit and the population abruptly collapses without much warning signal \citep{OsmKla17,Kla20}. This behaviour is very different from the dynamics of the lag under classic models of quadratic selection on moving optimum. \cite{OsmKla17} assumed a Gaussian distribution of phenotypes and a constant phenotypic variance and compared their analytical results to simulations of a sexually reproducting population. \cite{Kla20} went on to show that non quadratic fitness function with inflection points, leading to such tipping points, could emerge from various realistic ecological feedbacks involving density-dependence or interactions with other species. Our analytical results in the infinitesimal model allow us to recover very similar patterns to these previous studies and to predict the critical speed at which such evolutionary tipping points occur. We furthermore show that evolutionary tipping points also emerge in the asexual model, but with a different signature.
 In the asexual model, there is only one possible equilibrium for each value of the speed of environmental change. Again, ultimately, this unique equilibrium is due to the fact that the variance evolves more freely in the asexual model. As the speed increases towards the critical value $c^*_{tip}$, the lag diverges (Figure \ref{fig:formul_m_bounded}(a)-\ref{fig:m_bounded}(a)). As a result, the variance gets arbitrarily large and {the skewness becomes  negative, which shows that more individuals lag behind the evolutionary lag.} Conversely, in the infinitesimal model, the variance is constrained to remain nearly constant, resulting in multiple equilibria, which determine several basins of stability, up to the critical value $c^*_{tip}$. The lag remains bounded in the vicinity of the tipping point, determining a characteristic range for the basin of attraction of the origin (Figure \ref{fig:formul_m_bounded}(b)-\ref{fig:m_bounded}(b)). The lag can diverge, even if $c<c^*_{tip}$, for maladapted initial distributions concentrated far from the origin. This corresponds to a population that cannot keep pace with the environmental change. 
 
 \paragraph*{Effect of the mutation kernel.}
In the asexual model, our results also give analytical insights on the effect of the shape of the mutation kernel on the adaptation to a changing environment. Empirical data on the exact distribution of mutational effects on phenotypic traits are hard to get (even though there is more data on the fitness effects of mutations)~\cite[see e.g.][]{HalKei09,Nei14}. 
Most models therefore assume for mathematical convenience a Gaussian distribution of mutational effects. A few simulation studies have however explored marginally the consequences of a different, leptokurtic, mutation kernel \citep{KeiHil88,Bur99,WaxPec99}
: they found that a fatter tail for the distribution of mutational effects led to higher phenotypic variance, smaller evolutionary lag and greater fitness. The present analytical results are consistent with these past simulation results and show that we may expect in general distributions of mutations with higher kurtosis to reduce maladaptation and improve fitness, especially in fast changing environments.

\paragraph*{The advantage of sex in changing environments.}
Previous studies~\citep{Cha93,Bur99,WaxPec99} have used the Gaussian assumption and/or simulations to compare the dynamics of adaptation to a changing environment in sexual and asexual organisms. They all reached the conclusion that sex should provide a net advantage in a directionally changing environment, with a lower lag and greater fitness, which was ultimately due to the greater standing variance evolving in a sexually reproducing populations. More precisely, \cite{Bur99} and \cite{WaxPec99} found that the standing variance in sexual organisms would increase significantly with the speed of environmental change, while it would have only moderate effects on the variance in the asexual population. These findings seem to contrast with our comparison of the asexual model and sexual infinitesimal model, with more constraints on the evolution of the phenotypic variance for the latter. However, we would warn against interpreting our comparison of the infinitesimal and asexual model as informing about the advantage of sex in a changing environment. We rather see this comparison as informing us about the consequences of some modeling choices, with various constraints on the evolution of the phenotypic variance. First, for the ease of comparison between models, we set the parameter sigma to determine the amount of new variation introduced through reproduction in the progeny of parents in both models: in the asexual model it describes the amount of variance introduced by mutation, while it describes variation due to segregation in the infinitesimal model. It is unclear whether these quantities would be comparable with an explicit genetic model, including mutation and segregation at a finite set of loci. Second, we note that both \cite{Bur99} and ~\cite{WaxPec99}  used in their simulations parameter values for mutation and selection corresponding well to the regime of the House of Card approximation~\citep{Tur84}, with rare mutations of large effects on fitness. {These approximation regime is in sharp contrast with our assumption of small variance. In particular in the asexual model, we assume that the mean effect of mutation $\sigma^2$ is small compared with the frequency of mutations, which may be captured by $\betab$.}


\paragraph{Conclusions and perspectives.}
{ One of the main conclusion of our study is that the genetic standing variance at equilibrium truly depends on the modelling choice of the mode of reproduction. To understand this relationship, the approximation of the phenotype distribution appeared necessary. This approach is indeed robust, as shown by several studies following the same methodology in spatial structured population models: discrete patches~(\citep{Mirrahimi-2017} in asexual model and~\citep{Dekens-2020} in the infinitesimal sexual model); dispersal evolution (\citep{Perthame-Souganidis-2016,Lam-Lou-2017,Lam-2017,Hao-Lam-Lou-2021,Cal-Hen-Mir-Tur-2018} in the asexual case and \citep{Dek-Lav-2021} in the infinitesimal sexual case). Moreover, this methodology is expected to be efficient to investigate other structured population models. Our next step will be to study the adaptation of an age--structured population to a changing environment, following~\citep{CotRon14}. Other modes of reproduction with a more complicated genetic underlying architecture are also under investigation, \citep[see for instance][]{Dek-Mir-2021,Dek-Ott-Cal-2021}. 
}

%
%
%
%

\bibliographystyle{plainnat}
\bibliography{biblio_jimmy_single_brace_francais}

\newpage
\appendix

\section*{Supplementary material/Appendix}

The following subsections gather mathematical analysis supporting the adimensional scaling, numerical methods, Taylor expansions and formula derived in the main text. Although some parts are standard methods (rescaling, numerics), some parts are original contributions (dedicated Taylor expansions and formula involving the Lagrangian function), extending the literature in multiple ways. Hence, this supplementary material can be read  as the companion mathematical paper of the main text. 

Before we enter into the technical details, let us highlight some important observations about the Taylor expansions: 
\begin{itemize}
\item These expansions are more than moment closure methods, where one usually tries to guess the higher moments of the distribution in order to derive a close system of equations on some scalar quantities (first moments of the distribution, {\em e.g.} population size, evolutionary lag value, etc). Here, the {\em whole distribution} is approximated, then scalar quantities are deduced without any {\em a priori} assumptions on the shape of the distribution. 
\item In contrast to classical expansions of the distribution $F$ which are {\em linear}, {\em e.g.} $F = F_0 + \eps F_1 + \dots$, we perform here a {\em multiplicative} Taylor expansion, meaning a linear expansion of the logarithm of the density: $U = U_0 + \eps U_1 + \dots$. We claim this is the natural expansion in the regime of small variance in order to discard the variance from the asymptotic calculations. Nonetheless, intermediate computations may appear heavy because of the nonlinear nature of the multiplicative expansion. 
\item We believe all these approximations can be theoretically justified, and error terms can be controlled quantitatively up to some extent. Results in the literature so far cover the case without environmental change (c = 0), see~\citep{PerBar08,MirRao13} (Barles, Perthame, Mirrahimi et al) for the asexual model, and the more recent~\citep{CalGarPat19,Pat20} for the infinitesimal sexual model.  
\end{itemize}

\section{Derivation of generic formula}
Let us consider the equilibrium of our model:
\begin{equation}
\lambda F(z) - \eps^\gamma c  \partial_z F(z)  + m(z) F(z) = \mathcal{B}(F)(z)  \, , \quad \gamma\in \{1,2\}
\end{equation}
By integration over $\R$, we find: 
\begin{equation}
\lambda \rho + \int_\R m(z) F(z)\, dz = \rho\, , \quad \rho = \int_\R  F(z)\, dz\,.
\end{equation}
In the regime of small variance, we expect $F$ to concentrate around the evolutionary lag $z^*$, so as to get the following relationship 
\begin{equation}\label{eq:demographic eq}
\lambda \approx 1 - m(z^*)\,, 
\end{equation}
which corresponds to the demographic equilibrium. Next, we multiply by $(z-z^*)$, where $z^*$ is the mean value of the distribution $F$. Then, we integrate over $\R$ to find: 
\begin{equation}
   \eps^\gamma c \rho + \int_\R (z-z^*) m(z) F(z)\, dz = \int_\R (z - z^*) \mathcal{B}(F)(z)\, dz  \,.
\end{equation}
For any operator $\mathcal{B}$ defined by~\eqref{eq:asex sex}, we find that the right-hand-side vanishes by definition of $z^*$. 
The concentration of the distribution $F$ motivates the Taylor expansion of the selection function: $m(z) \approx m(z^*) + (z- z^*)\partial_z m(z^*)$ which implies the following:
\begin{equation}\label{eq:evol eq}
 \eps^\gamma c  \approx - \partial_z m(z^*) (\Var(F)) \, .
\end{equation}

\section{Adimensional scaling}\label{app:scaled_pb}
We present in this  section the details of the scaling procedure which leads to equations~\eqref{eq:asexual no age} and~\eqref{eq:sexual no age} in adimensional form. By convention, the variables and parameters in original units are written in bold, whereas adimensional quantity are in normal font. 

The stationary state $(\lambdab,\Fb)$ satisfies
\begin{equation*}\lambdab \Fb(\zb) - \cb \partial_{\zb} \Fb(\zb)  + \mub(\zb) \Fb(\zb) = \betab \mathcal{B}( \Fb )(\zb)\, .
\end{equation*}
where the mortality rate $\mu(z) = \mub_0 + \mb(\zb)$ is decomposed as a basal rate $\mub_0$ (miminum value of the mortality rate), and a mortality increase $\mb(\zb)\geq 0$ which is trait-dependent. 
Dividing by the fertility rate $\betab,$ (trait-independent) it becomes 
\begin{equation}\label{eq_app:F}
\dfrac{\lambdab+\mub_0}\betab  \Fb(\zb) - \dfrac{\cb}{\betab} \partial_{\zb} \Fb(\zb)  + \dfrac{1}{\betab}\mb(\zb) \Fb(\zb) =   \mathcal{B}( \Fb )(\zb)\, .
\end{equation}
Around the optimum trait $\zb=0,$ the mortality per individual per generation $\mb/\betab$ is equivalent to 
$$\bf \dfrac1\betab\mb(z) = \dfrac12 \dfrac{\mb''(0)}{\betab} \zb^2 + o(\zb^2)= \dfrac12 \left(\zb\sqrt{\dfrac\alphab\betab}\right)^2 +o(\zb^2)$$
It is natural to measure traits with this selection scale: 
$$\Zc^2 =  \dfrac{\betab}{\alphab}\,.$$
The { mean fitness} and the phenotypic distribution becomes in the scaled trait variable $z = \zb/\Zc$:
$$ \lambda =  \dfrac{\lambdab+\mub_0}\betab \quad \hbox{ and } \quad    F( z) = \Fb(\Zc z)\,.$$
The mortality rate per individual becomes 
\begin{equation*}
 m(z) = \dfrac{1}{\betab}\mb\left ( \Zc z \right ),
\end{equation*}
so that the selection strength $\alphab$ around the optimum is scaled to a unit value:
\begin{equation*}
m''(0) = 1\,.
\end{equation*}
Our main assumption is that there is a small variability with respect to the selection scale $\Zc$. Denoting by $\Zv$ the 
standard deviation of offspring traits from the parental traits, we define $\epsilon$ the scaling ratio:
$$\epsilon = \frac{\Zv}{\Zc}. $$
Then, our main assumption can be summarized as $\epsilon\ll 1$, paving the way to suitable Taylor expansions. 
In both models, the standard deviation $\Zv$ is denoted by the common parameter $\sigmab$ in the original units. However, we emphasize that it corresponds to mechanisms of variability associated with very different genetical background. 

The reproduction operators $\mathcal{B}$ are transformed   as follows:
\paragraph*{Asexual reproduction operator in scaled variables.}
\begin{equation*}
 \mathcal{B}(\Fb)(\Zc z) = \dfrac1\sigma \int_\R K\left(\dfrac{\Zc z-\zb'}{\sigma}\right)\Fb(\zb')\,d\zb' 
                              = \dfrac1\Zv \int_\R K\left(\dfrac{\Zc}{\Zv} \left(z-\dfrac{\zb'}\Zc\right)\right)\Fb(\zb')\,d\zb'. 
\end{equation*}
Using the change of variable $z' = \zb'/\Zc$ in the integral and the definitions of $\epsilon=\Zv/\Zc$ and $F,$ we obtain
\begin{equation*}
 \mathcal{B}(\Fb)(\Zc z) = \dfrac\Zc\Zv \int_\R K\left(\dfrac{\Zc}{\Zv} (z-z')\right)\Fb(\Zc z')\,dz'
                              = \dfrac1\epsilon \int_\R K\left(\dfrac{z-z'}{\epsilon}\right) F(z')\,dz'.
\end{equation*}

\paragraph*{Sexual reproduction operator in scaled trait.}
\begin{equation*}
\begin{array}{rl}
   \mathcal{B}(\Fb)(\Zc z)  &\ds  = \dfrac{1}{\sqrt{\pi\sigmab^2}} \iint_{\mathbb{R}^2}  \exp\left ( - \left(\Zc z - \dfrac{ \zb_1 +  \zb_2}2 \right)^2 \right ) \Fb(\zb_1) \ds\dfrac{\Fb(\zb_2)}{ \int_\R \Fb( \zb_2')\, d \zb_2'}\, d\zb_1 d \zb_2 \\[4mm]
                                 &\ds  =  \dfrac{1}{\sqrt{\pi}} \dfrac{1}{\Zv}   \iint_{\mathbb{R}^2}  \exp\left ( - \left (\dfrac{\Zc}{\Zv}\right)^2 \left(z - \dfrac12\left(\dfrac{\zb_1}\Zc + \dfrac{\zb_2}\Zc\right)\right )^2 \right ) \Fb( \zb_1)\dfrac{ \Fb( \zb_2)}{ \int_\R  \Fb( \zb_2')\, d \zb_2'}\, d \zb_1 d \zb_2 
\, .
\end{array}
\end{equation*}
Using the change of variable $z_1 = \zb_1/\Zc,$ $z_2 = \zb_2/\Zc,$ and $z_2' = \zb_2'/\Zc,$ in the integrals and the definitions of $\epsilon=\Zv/\Zc$ and $F,$ we obtain
\begin{equation*}
\begin{array}{rl}
   \mathcal{B}(\Fb)(\Zc z) &\ds  =  \dfrac{1}{\sqrt{\pi}} \dfrac{1}{\Zv}   \iint_{\mathbb{R}^2}  \exp\left ( - \left (\dfrac{\Zc}{\Zv}\right)^2 \left(z - \dfrac{{z_1} + {z_2}}2\right )^2 \right ) \Fb(\Zc z_1)\dfrac{  \Fb(\Zc z_2)}{ \Zc\int_\R  \Fb(  \Zc z_2')\, dz_2'}\, \Zc^2 dz_1 dz_2 \\[4mm]
                                &\ds  =  \dfrac{1}{\sqrt{\pi}} \dfrac{\Zc}{\Zv}   \iint_{\mathbb{R}^2}  \exp\left ( - \left (\dfrac{\Zc}{\Zv}\right)^2 \left(z - \dfrac{{z_1} + {z_2}}2\right )^2 \right ) F(z_1)\dfrac{ F(z_2)}{ \int_\R   F(z_2')\, dz_2'}\, dz_1 dz_2 \\[4mm]
                                &\ds  =  \dfrac{1}{\epsilon\sqrt{\pi}}  \iint_{\mathbb{R}^2}  \exp\left( - \dfrac{1}{\epsilon^2} \left(z - \dfrac{{z_1} + {z_2}}2\right)^2 \right) F(z_1)\dfrac{ F(z_2)}{ \int_\R   F(z_2')\, dz_2'}\, dz_1 dz_2 
\, .
\end{array}
\end{equation*}

\paragraph*{The adimensional speed.}
It remains to express the adimensional speed $c = \cb/C$ with different choices of the typical speed $C$. This choice depends on the reproduction mode as follows:
\begin{equation}\label{eq:Cada}
  C = \begin{cases}
\sigmab\betab & \text{(asexual model)}\medskip\\
 \sigmab^2\sqrt{\alphab\betab}& \text{(infinitesimal sexual model)}
\end{cases}\, . 
\end{equation}
We thus deduce the adimensional the following expression of the advection term:
\begin{equation}\label{eq:speed adim}
- \dfrac{\cb}{\beta} \partial_{\zb} \Fb(\zb) = - c \dfrac{C}{\beta \Zc}  \partial_z F(z) = \begin{cases}
- c \dfrac{\sigma}{\Zc}  \partial_z F(z) = - \eps c  \partial_z F(z) & \text{(asexual model)}\medskip\\
- c \dfrac{\sigmab^2 \alpha^{1/2}}{\Zc \beta^{1/2}}   \partial_z F(z) = - \eps^2 c  \partial_z F(z)& \text{(infinitesimal sexual model)}
\end{cases}\, . 
\end{equation}
We obtain eventually the two rescaled problems as shown in~\eqref{eq:asexual no age} and~\eqref{eq:sexual no age}. To conclude, let us mention that the discrepancy between the two values of $C$ \eqref{eq:Cada} is due to the very last step \eqref{eq:speed adim}, where the adimensional speed must be of order $\epsilon$ in the asexual model, resp. of order $\eps^2$ in the infinitesimal sexual model, in order to balance the other contributions. A  mismatch at this step ({\em e.g.} any other power of $\eps$) would result in a severe unbalance between the contributions, namely dramatic collapse of the population if the effective speed is too large, or no clear effect of the change if the effective speed is too small. 

%

\section{Derivation of the variance}\label{sec:var_derivation}

We compute below the formula of the standing variance $\Var (F)$ in terms of $U = - \eps^\gamma \log F$,
\begin{equation}
\Var(F) = \left( \int_\R( \left (z - z^*\right )^2 \exp\left (-\dfrac{U(z)}{\eps^\gamma}\right )\, dz\right ) \Big/\left (\int _\R   \exp\left (-\dfrac{U(z)}{\eps^\gamma}\right )\, dz\right ) 
\end{equation}
We assume that $U$ reaches a non-degenerate minimum point at a unique $z^*$, such that $U(z) = U(z^*) + \frac12 (z - z^*)^2 \partial_z^2 U(z^*) + o((z - z^*)^2)$ as $z\to z^*$.
The denominator is equivalent to 
\begin{equation}
\dfrac{\eps^{\gamma/2}}{\sqrt{2\pi} \sqrt{\partial_z^2 U(z^*)}}\exp\left ( -\dfrac{U(z^*)}{\eps^\gamma} \right )
\end{equation}
whereas the numerator is equivalent to 
\begin{equation}
\dfrac{\eps^\gamma}{\partial_z^2 U(z^*)}\dfrac{\eps^{\gamma/2}}{\sqrt{2\pi} \sqrt{\partial_z^2 U(z^*)}}\exp\left ( -\dfrac{U(z^*)}{\eps^\gamma} \right )\, .
\end{equation}
Thus, the ratio is equivalent to \eqref{eq:variance}:
\begin{equation}
\Var(F) \sim \dfrac{\eps^\gamma}{\partial_z^2 U(z^*)}\, . 
\end{equation}

\section{Asexual type of reproduction (Details of Section \ref{sec:asex})}\label{app:asexual}
This long section is devoted to the details of the Taylor expansion of $U$ defined by~\eqref{eq:log transform}. The equations verified by the successive terms $U_0$ and $U_1$ are derived. The meaningful  formula are computed.

We can formally expand the pair $(\lambda,U)$ with respect to $\epsilon$ as follows,
\begin{equation}\label{eq_app:ansatz asexual} \begin{cases} 
U(z) = U_0(z) \textcolor{gray}{+ \eps U_1(z) + o(\eps)} \\
\lambda  = \lambda_0 \textcolor{gray}{+ \eps \lambda_1  + o(\eps)}
\end{cases}
\end{equation}
where $(\lambda_0,U_0)$ gives the limit shape  as $\epsilon\to0$, and $(\lambda_1,U_1)$ is the correction for small $\epsilon>0$. We focus on the   leading order contribution in this work. The corrector is required to refine our approximation in some part of the discussion. 


\subsection{Equations for $(\lambda,U),$ $(\lambda_0,U_0)$ and $(\lambda_1,U_1)$}\label{app:Mastre_eq_asex}

We begin with the diffusion approximation for the sake of simplicity. This enables to present the main ingredient, namely {\em the completion of the square} in the equation, that will be generalized next for a general mutation kernel.

\subsubsection{The diffusion approximation}\label{app:diffusion_approx}

The equation for $F$ \eqref{eq:asexual no age}, together with the logarithmic transformation $F(z) = \exp(-U(z)/\eps)$, is equivalent to the following one:
\begin{equation}
\lambda  +  c  \partial_z U(z)  + m(z)   =  1+  \frac12 \left (\partial_z  U(z) \right)^2 + \frac\eps2 \partial_z^2U(z) \, .
\end{equation} 
Clearly, the limiting problem for $(\lambda_0,U_0)$ is
\begin{equation}
\lambda_0 + c \partial_z U_0(z) + \m(z) = 1 + \frac12 \left (\partial_z  U_0(z) \right)^2  \, .\label{eq:HJ asexual U0 apporx diff}
\end{equation} 
It is instructive to gather all the $\partial_z U_0$ in the right hand side, then to complete the square:
\begin{equation}
\label{eq:U0 approx diff}
m(z) + \left [ \lambda_0 - 1 + \dfrac{c^2}{2}\right ] = 
\frac12\left ( \partial_z U_0(z) - c\right )^2   \,.\end{equation} 
The key point is that there exist admissible solutions of this ODE if, and only if, the value between brackets vanishes, {\em i.e.}  $\lambda_0 = 1  - \frac{c^2}{2}$. The argument is as follows.

\paragraph*{Completion of the square.}
On the one hand, evaluating \eqref{eq:U0 approx diff} at $z = 0$, we find that $\lambda_0 - 1 + \frac{c^2}{2}\geq 0$ since $m(0) = 0$. On the other hand, if $\lambda_0 - 1  + \frac{c^2}{2}$ is positive, then  $\partial_z U_0  - c$ does not change sign. Assuming without loss of generality that it is everywhere positive, we find that $U_0(z) \geq c z + U_0(0)$ for $z \geq 0$ and $U_0(z) \leq c z + U_0(0)$ for $z \leq 0$. In particular, we have $U_0(z) \to -\infty$ as $z\to -\infty$, and $U_0(z) \to +\infty$ as $z\to +\infty$, which is clearly not admissible because $F$ is a population density. Therefore, $\lambda_0 - 1 + \frac{c^2}{2} = 0$.

Next, we can deduce the lag by evaluating \eqref{eq:HJ asexual U0 apporx diff} at $z_0^*$ such that $\partial_z U_0(z_0^*) = 0$, 
\begin{equation}
\label{eq:lambda and co approx diff}
m(z_0^*) = \dfrac{c^2}{2}, 
\end{equation}
and also the value of the second derivative by differentiating once and evaluating at $z_0^*$:
\begin{equation}
 c \partial^2_z U_0(z_0^*) + \partial_z m(z_0^*) = 0.
\end{equation}
Finally, we deduce the variance from \eqref{eq:variance}
\begin{equation}
\Var(F) = -\dfrac{\eps c}{\partial_z\m(z^*_0)} + o(\eps)\, 
\end{equation}
consistently with \eqref{eq:2relations}. 

We can even provide a formula for the profile $U_0$ by solving the ODE \eqref{eq:U0 approx diff}:
\begin{equation}\label{eq:diff profile}
U_0(z) = c z + \left| \int_0^z \left( 2 m(z') \right)^{1/2}\,dz'\right|\, .
\end{equation}
Notice that the environmental change acts here as a linear correction of the equilibrium profile obtained in the case $c = 0$. However, this is a peculiarity of the diffusion approximation.

It is another peculiarity that a quadratic selection function $m(z) = \frac{z^2}{2}$ results in a quadratic profile $U_0(z) = cz + \frac{z^2}{2}$ \eqref{eq:diff profile}, which corresponds to a Gaussian distribution function  $F$ with variance~$\eps$.  



\subsubsection{The case of a general mutation kernel}
Again, we can reformulate the problem~\eqref{eq:asexual no age} in an equivalent form:
\begin{equation}\label{eq:eq U eps}
\left ( \lambda + c \partial_z U(z) + m(z) \right) \exp\left(-\dfrac{U(z)}\eps\right) = \dfrac1{\eps} \int_\R  K\left(\dfrac{z-z'}\eps\right) \exp\left(-\dfrac{U(z')}\eps\right)\, dz'  
\end{equation}
After the change of variables $z' = z - \eps y$ in the integral term, we obtain: 
\begin{align*}
\lambda + c\partial_z U(z) + m(z)
& =      \int_\R   K\left(y\right) \exp\left( \dfrac{U(z)-U(z-\eps y)}\eps\right)\, dy 
\\
 &=   \int_\R   K\left(y\right) \exp\left( y \partial_z U(z) - \dfrac\eps2 y^2  \partial^2_z U(z) + o(\eps)\right)\, dy\, .
\end{align*}

Injecting \eqref{eq_app:ansatz asexual} into \eqref{eq:asexual no age}, but dropping terms of order higher than $\epsilon$, we get
\begin{align}
   \lambda_0 + \eps \lambda_1 + c\partial_z\left( U_0(z) + \eps U_1(z) \right) + \m(z) 
   & \ds =   \int_\R   K\left(y\right) \exp\left( y \partial_z \left( U_0(z) + \eps U_1(z) \right)  - \dfrac\eps2 y^2  \partial^2_z U_0(z)  + o(\eps) \right)\, dy \nonumber \\
   & \ds =   \int_\R   K\left(y\right) \exp\left( y \partial_z  U_0(z) \right) \left( 1 + \eps y \partial_zU_1(z)  - \dfrac\eps2 y^2  
   \partial^2_z U_0(z)  \right)\, dy +o(\eps)\, .
\label{eq_app:asex_dev}
\end{align}
By identification of the contributions having  the same order in $\epsilon$ in equation~\eqref{eq_app:asex_dev}, we obtain the following  equations for the pairs $(\lambda_0,U_0)$ and $(\lambda_1,U_1)$
\begin{align}
\label{eq_app:U0}
      \hbox{\bf Limit problem:}   & \quad
      \lambda_0 + c \partial_z U_0(z) + \m(z) = 1 + H  \left(\partial_z  U_0(z) \right)  \, , \\[3mm]
      \label{eq_app:asex_U1}
      \hbox{\bf First correction problem:}  & \quad
      \lambda_1 + \left( c - \partial_p H(\partial_z  U_0(z))\right) \partial_z  U_1(z) =    - \dfrac12 \partial^2_p H(\partial_z  U_0(z)) \partial^2_z U_0(z)\,,
\end{align}
where the Hamiltonian function $H$ is the two-sided Laplace transform of $K$ up to an additive constant:
\begin{equation*}
H(p) =  \int_\R   K\left(y\right) \exp\left( y p \right)\, dy - 1  \, , \quad \partial_p H(p) = \int_\R   y K\left(y\right) \exp\left( y p \right)\, dy \, , \quad \partial^2_p H(p) = \int_\R   y^2 K\left(y\right) \exp\left( y p \right)\, dy\, .
\end{equation*}

\subsubsection{Computation of the mean fitness}\label{app:mean_fitness_asex}

The argument of Section \ref{app:diffusion_approx} for computing $\lambda_0$ can be extended to the general case. Quadratic functions are replaced by convex ones, but the argument is essentially the same.

Again, let us reorganize \eqref{eq:HJ asexual U0} as follows, gathering  the $\partial_z U_0$ in the right hand side,
\begin{equation}
m(z) +  \lambda_0 - 1  = - \left (  c \partial_z U_0(z) -  H(\partial_z U_0(z)) \right )  \, .  \label{eq:HJ asex general}
\end{equation}
The function $p\mapsto cp - H(p) $ reaches a maximum value, denoted as $L(c)$ by definition \eqref{eq:L}. Adding this value on each side, we find
 \begin{equation}
m(z) + \left [ \lambda_0 - 1 + L(c) \right] = H(\partial_z U_0(z)) - c \partial_z U_0(z) + L(c) \, . 
\end{equation}

\begin{figure}
\begin{center}
\includegraphics[width = 0.6\linewidth]{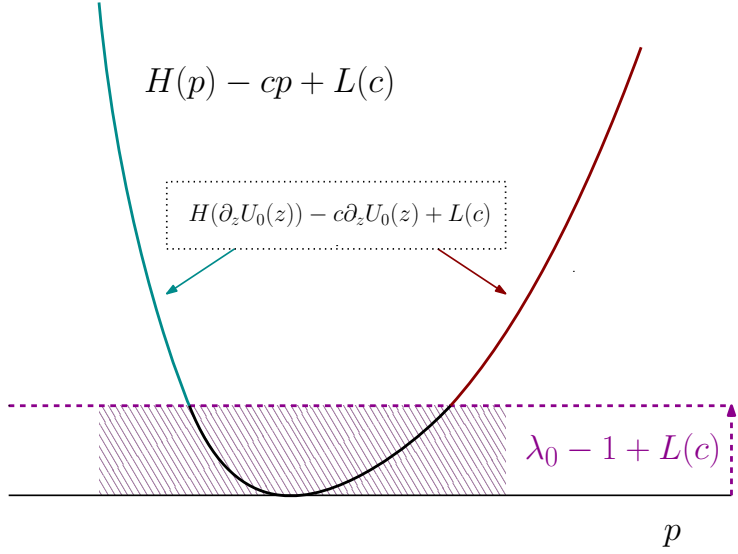}
\end{center}
\caption{Completion of the generalized square $H(p) - cp + L(c)\geq 0$ (with zero minimum value). It is a true quadratic expression in the case of the diffusion approximation. It is a convex function in the general case.}
\label{fig:square}
\end{figure}

\paragraph*{Completion of the generalized square.}
As in \eqref{eq:U0 approx diff}, the function $p\mapsto H(p) - cp + L(c)$ in the right-hand-side is convex, nonnegative and touches zero. This is the analogous computation of {\em the completion of the square} by means of adding $L(c)$.
The same reasoning as above implies that  the constant between brackets must vanish, {\em i.e} $\lambda_0 = 1  - L(c)$. Otherwise, the quantity $H(\partial_z U_0(z)) - c \partial_z U_0(z) + L(c)$ would take positive values for $z\in \R$, hence the function $\partial_z U_0(z)$ could  take values only on one of the two branches of the function $p\mapsto H(p) - cp + L(c)$, as depicted in Fig~\ref{fig:square}.
As the function $p\mapsto H(p) - cp + L(c)$ is invertible on each separate branch, we could determine unambiguously the value of $\partial_z U_0(z)$ for $z\in \R$. In particular, it would have the same limiting value (possibly infinite) as $z\to -\infty$ and $z\to +\infty$ since $m(-\infty) = m(+\infty)$. 
This would preclude the asymptotic behavior $U_0(\pm \infty) = +\infty$ which is equivalent to vanishing population density at infinity. Hence, $\lambda_0 = 1  - L(c)$ is the only possible value. 

\subsection{Summary}\label{sec:app summ}

So far we have obtained an analytical formula for the mean fitness, 
\begin{equation}
\lambda_0 = 1- L(c)\,, \label{eq:app mean fitness}
\end{equation}
by means of the Lagrangian function which is the Legendre transform of the Hamiltonian function,
\begin{equation}
L(c) = \max_p \left( p c - H(p)\right)\,,
\end{equation}
where $H$ is the Laplace transform of the mutation kernel $K$. 

The knowledge of the mean fitness enables deriving the lag load, which equilibrates birth and death in the population concentrated at trait $z_0^*$: $\lambda_0 = 1 - m(z_0^*)$, or equivalently
\begin{equation}\label{eq:app mean trait}
m(z_0^*) = L(c)\,.
\end{equation}
Note that the latter is equivalent to setting $\partial_z  U_0(z_0^*) = 0$ in \eqref{eq:HJ asex general} (critical point of the density), which is another characterization of the lag load.  

The variance can be completed subsequently by differentiating \eqref{eq:HJ asex general} with respect to $z$ and evaluating at $z = z_0^*$. It is found that the variance equilibrates the fitness gradient and the speed of environmental change ({\em i.e.} the variations in the trait value in the moving frame):
\begin{equation}\label{eq:app d2U0}
\partial^2_z U(z_0^*) =  -\dfrac{\partial_z\m(z^*_0)}{c} \, .
\end{equation}

\subsection{Conjugacy: Enlightening heuristics}\label{app:conjugacy}
 There exists an alternative way to get some of the previous formula. The idea is to twist the unknown distribution $F$ by a well chosen exponential function, in order to remove the transport part $-c \partial_z F$ due to the environmental change. An enlightening example is the case of the diffusive approximation. Suppose the model is 
 \begin{equation}
 \lambda F(z) - \eps c \partial_z F(z) - \dfrac{\eps^2}{2} \partial_z^2 F(z) = (1 - m(z)) F(z)\, .
 \end{equation}
 Then, the twisted distribution $\mathfrak{F}(z) = F(z) e^{c z/\eps}$ satisfies the following equation:
 \begin{equation}
 \lambda \mathfrak{F}(z)   - \dfrac{\eps^2}{2} \partial_z^2 \mathfrak{F}(z) = \left (1 - \dfrac{c^2}{2} - m(z)\right ) \mathfrak{F}(z)\, .
 \end{equation}
 Therefore, we are reduced to a simpler problem without environmental change, at the expense of a global increase of mortality of value $c^2/2$, consistently with the result of Section \ref{app:diffusion_approx}. 

 However, the general case is based on heuristics rather than formal arguments. Starting from equation \eqref{eq:asexual no age}, or equivalently:
 \begin{equation}
 \lambda F(z) - \eps c \partial_z F(z) - \int_\R K_\eps(z-z')\left ( F(z') - F(z) \right )\, dz' = (1 - m(z)) F(z)\,,
 \end{equation}
 the density $F$ is replaced with $\mathfrak{F}(z) = F(z) e^{p_0z/\eps}$, for some $p_0\in \R$ to be characterized later on. The equation for $\mathfrak{F}$ is:
 \begin{equation}
 \lambda \mathfrak{F}(z) +   c p_0 \mathfrak{F}(z) - \eps c  \partial_z \mathfrak{F}(z) - \int_\R K_\eps(z-z')\left ( e^{p_0(z-z')/\eps} \mathfrak{F}(z') - \mathfrak{F}(z) \right )\, dz' = (1 - m(z)) \mathfrak{F}(z)\,,
 \end{equation}
 It is useful to rearrange the terms as follows:
 \begin{equation}\label{eq:conjugacy}
 \lambda \mathfrak{F}(z)  - \eps c  \partial_z \mathfrak{F}(z) - \int_\R K_\eps(z-z') e^{p_0(z-z')/\eps}\left (  \mathfrak{F}(z') - \mathfrak{F}(z) \right )\, dz' = \left (1 -  c p_0 + \left ( \int_\R K_\eps(z')e^{p_0 z'/\eps}\,dz' - 1\right ) - m(z)\right ) \mathfrak{F}(z)\,,
 \end{equation}
 A natural way to choose $p_0$ is  to guarantee that the combination of transport and mutations preserves the center of mass of the distribution. This is a way to remove artificially the asymmetrical transport part. Thus, we propose the following characterization of $p_0$: for any distribution $\mathfrak{F}$,
 \[ \int_\R z \left ( - \eps c   \partial_z \mathfrak{F}(z)  - \int_\R K_\eps(z-z') e^{p_0(z-z')/\eps}\left (  \mathfrak{F}(z') - \mathfrak{F}(z) \right )\, dz'\right )\, dz = 0\, . \]
 This is equivalent to:
 \begin{align*}
 \eps c \int_\R \mathfrak{F}(z)\, dz &= \iint z K_\eps(z-z') e^{p_0(z-z')/\eps}   \mathfrak{F}(z')\, dz' dz - \iint z K_\eps(z-z') e^{p_0(z-z')/\eps} \mathfrak{F}(z)  \, dz' dz\\
 &= \iint z K_\eps(z-z') e^{p_0(z-z')/\eps}   \mathfrak{F}(z')\, dz' dz - \iint z' K_\eps(z-z') e^{p_0(z-z')/\eps} \mathfrak{F}(z')  \, dz' dz\\
 &= \iint (z - z') K_\eps(z-z') e^{p_0(z-z')/\eps}   \mathfrak{F}(z')\, dz' dz \\
 &= \left ( \int z K_\eps(z) e^{p_0z/\eps}\, dz \right ) \left ( \int_\R \mathfrak{F}(z)\, dz\right ) 
 \, .
\end{align*}
 Finally, the required  condition is equivalent to the following one, which appears to be independent of $\eps>0$:
 \begin{equation}
 c = \int y K(y)e^{p_0 y} \, dy\, .
 \end{equation}
 With the notations of Section \ref{app:Mastre_eq_asex}, this is also $c = \partial_p H(p_0)$. The right hand side of \eqref{eq:conjugacy} becomes:
 \begin{equation}
 \left ( 1 - c p_0 + H(p_0) - m(z)\right ) \mathfrak{F}(z) = \left ( 1 - L(c) - m(z)\right ) \mathfrak{F}(z)\, . 
 \end{equation} 
 As a conclusion, we have shown that the combination of transport and mutations is equivalent to an operator which preserves the center of mass, up to a global increase of mortality of value $L(c)$.

\subsection{Some properties of the Hamiltonian and Lagrangian functions}\label{app:kernel_kurtosis}

We gather below some classical properties of the special functions that appeared useful in the analysis above. 

\paragraph*{Diffusion approximation as an extremal case of the convolution case.}
By symmetry of the kernel $K$, and its properties, the Hamiltonian function can be bounded below:
\begin{equation}
H(p) =  \int_\R K(y) \left ( \dfrac{\exp(yp) + \exp(-yp)}2 - 1\right )dy   \geq \dfrac{|p|^2}2 \int_\R K(y) y^2dy = \dfrac{|p|^2}2\,.
\end{equation}
The latter expression is realized by the so-called diffusion approximation, see Section \ref{app:diffusion_approx}. Indeed, the Hamiltonian function there was simply the square of the gradient \eqref{eq:HJ asexual U0 apporx diff}. It is a direct consequence of the formula $L(c) = \max_p p c- H(p)$ (completion of the generalized square) that the Lagrangian function is bounded above:
\begin{equation}
L(c) \leq \dfrac{c^2}{2}\, .
\end{equation}
Hence, the maximum of lag load is realized for the diffusion approximation.

\paragraph*{The Hamiltonian function contains all the moments of the mutation kernel.}

By definition of the exponential function we have:
\begin{align*}
H(p) &= \int_\R K(y) \left ( \sum_{k = 0}^{\infty} \frac{(py)^k}{k!} \right )\, dy - 1\\
& = \sum_{k = 1}^{\infty} \left (  \int_\R K(y) y^k\, dy\right ) \dfrac{p^k}{k!}\, .
\end{align*}
Hence, the moments of $K$ are successive derivatives of $H$ at the origin. 

\paragraph*{Influence of the kurtosis of the mutation kernel.}

As an immediate consequence, we see that the mean fitness $\lambda_0 = 1-L(c)$ crucially depends on the full shape of the mutation kernel $K$. Indeed, the Lagrangian function $L$ is related to the Laplace transform of the mutation kernel $K$ \eqref{eq:H} via the Legendre transform \eqref{eq:L}. To investigate this relationship, we investigate five kernels having the same variance, but different shapes, see Table \ref{tab:shapes}. We can show from the Taylor expansions that the Hamiltonian functions are ordered from top to bottom as follows:
\begin{equation}
H_{\rm diff} \leq H_{\rm unif} \leq H_{\rm gauss} \leq H_{\rm exp} \leq H_{\rm gamma}\, .
\end{equation}
Accordingly, the Lagrangian functions are ordered in the opposite way, and the resulting mean fitnesses are ordered as follows:
\begin{equation}
\lambda_{\rm diff} \leq \lambda_{\rm unif} \leq \lambda_{\rm gauss} \leq \lambda_{\rm exp} \leq \lambda_{\rm gamma}\, .
\end{equation}
Hence, the lag load is ordered with respect to the kurtosis of the kernel.


  \begin{table}
\begin{center}
\renewcommand{\arraystretch}{2.5}
\begin{tabular}{l l l}
{\sffamily  } & \sffamily Mutation kernel $K(y)$ &\sffamily Hamiltonian function $H(p)$\\\hlineB{3}
\rowcolor{gray!30}
Diffusion approximation       & $\dfrac{1}{2}\partial_z^2 $ & $\dfrac12 p^2$  \\
Uniform distribution   & $\dfrac{1}{2\sqrt{3}} {\bf 1}_{(-\sqrt{3}, \sqrt{3})}$ & $ \dfrac{\sinh(\sqrt{3} p )}{\sqrt{3}p} - 1$\\
\rowcolor{gray!30}
Gaussian distribution    & $\dfrac{1}{\sqrt{2\pi}}\exp\left (-\dfrac{y^2}{2}\right )$ & $\exp\left (\dfrac{p^2}{2}\right ) - 1$ \\
Exponential distribution   &  $\dfrac{1}{\sqrt{2}} \exp\left (-\sqrt{2}|y|\right ) $ & $\dfrac{1}{1 - \dfrac{p^2}{2}} - 1$  \\
\rowcolor{gray!30}
Gamma distribution & $ |y|^{\gamma - 1 } \exp\left (- \sqrt{\gamma(\gamma+1)} |z|\right ) $ &  $ \dfrac12 \left ( (1 - \theta p )^{-\gamma} + (1 + \theta p )^{-\gamma} \right ) - 1$\\

\end{tabular}
\end{center}
\caption{(Left) Five examples of mutation kernels with same (unit) variance, ordered by increasing kurtosis (from top to bottom). (Right) The associated Hamiltonian functions, with analytical formula. The corresponding Lagrangian functions cannot be expressed with classical functions, but the first one, up to our knowledge.}
\label{tab:shapes}
\end{table}

\subsection{Consistency of the formula for $\partial_z^2U_0(z_0^*)$ at $c=0$}\label{app:asex_c0}
Here, we justify Remark \ref{rem:local shape}, meaning that the formula obtained for $\partial_z^2U_0(z_0^*)$ at $c>0$ \eqref{eq:app d2U0} coincides with the formula at $c=0$, namely $\partial^2_z U_0(0) = 1$. The latter is derived as follows. Firstly, the  { mean fitness}~\eqref{eq:app mean fitness} is $\lambda_0 = 1$, as $L(0) = 0$,   and the evolutionary lag~\eqref{eq:app mean trait} is naturally $z^*_0=0$ at $c=0$ by definition of the mortality rate, optimum at the origin. Secondly, the expression of $\partial^2_z U_0(0)$ can be obtained by two alternative ways.

By differentiating twice~\eqref{eq_app:U0} with respect to $z,$ yields
\[  \partial^2_z\m(z) =  \partial^2_p H  \left(\partial_z  U_0(z) \right) \left( \partial^2_z U_0(z)\right)^2 +  \partial_p H  \left(\partial_z  U_0(z) \right)   \partial^3_z U_0(z)  \, .  \]
By evaluating this expression at $z = 0$, the last contribution vanishes because $\partial_p H  \left(\partial_z  U_0(0) \right) = \partial_p H  \left(0\right)=0$. Hence, we get that
\[
   \partial_z^2U_0(0) = \left(\dfrac{\partial_z^2m(0)}{\partial^2_p H  \left(\partial_z  U_0(0) \right)}\right)^{1/2} =1 \,,
\]
since $\partial^2_z\m(0) = \partial^2_p H(0) = 1$. 

Alternatively, performing suitable Taylor expansions in expressions of, respectively, $z^*_0$ \eqref{eq: def z0star} and $\partial_z^2U_0(z_0^*)$~\eqref{eq:d2U_0}, as $c\to 0$, yields: 
\begin{align*}
&z^*_0 = \dfrac{\partial z^*_0}{\partial c} c + o(c) \,,\quad  \hbox{ and } \quad  
\dfrac12\partial^2_z\m(0) \left( \dfrac{\partial z^*_0}{\partial c} c \right)^2  = \dfrac12 \partial^2_v L(0) c^2\,, \\
& \partial^2_z U_0(0) = - \dfrac{\partial^2_z\m(0)}c \left(\dfrac{\partial z^*_0}{\partial c} c \right)=\partial^2_z\m(0) \left( \dfrac{\partial^2_v L(0)}{\partial^2_z\m(0)} \right)^{1/2} = \left(\partial^2_z\m(0) \partial^2_v L(0)  \right)^{1/2} = 1\, .
\end{align*}
By reciprocity of the derivatives of $H$ and $L$, we have $\partial^2_v L(0)  = 1/\left(\partial^2_p H(0)\right) = 1$. Both calculations coincide.


\subsection{Quantitative description of the first correction $(\lambda_1,U_1)$}
\label{sec:asex_U1}

We derive useful informations from the equation \eqref{eq_app:asex_U1} about the pair $(\lambda_1,U_1)$. 
The methodology goes as in Section \ref{app:Mastre_eq_asex}. 

We give the formula for the correctors $\lambda_1$, $z_1^*$, and the local shape around the minimal value: $\partial^2_z (U_0 + \eps U_1)(z^*_0 + \eps z^*_1)$. However, only the former one  ($\lambda_1$) is meant to be used in the main text, as it contains useful information about the mutation load in the population. 

The formula are summarized in the following list, which completes those obtained in Section \eqref{sec:app summ} at the leading order:
\begin{equation}
\begin{array}{ll}
\hbox{\bf  { Mean fitness} } & 
\lambda = 1  - L(c) - \dfrac\eps2   \left (\dfrac{ 1}{\partial^2_v L(c)}\right )^{1/2}    + o(\eps)   \medskip\\
\hbox{\bf  Evolutionary lag  } & z^* = z^*_0  + \dfrac\eps2 \left( \dfrac1{\partial_z \m(z^*_0)} \left( \dfrac{1}{\partial^2_v L(c)}   \right)^{1/2}  + \dfrac1{c} \right)  + o(\eps)   
\medskip\\
\hbox{\bf  Local shape} & \partial^2_z U(z^*) =  -\dfrac{\partial_z\m(z^*_0)}{c}  - \dfrac \eps2  \left(  \dfrac1{c} \dfrac{\partial^2_z \m(z^*_0)}{\partial_z \m(z^*_0)} \left( \dfrac{ 1}{\partial^2_v L(c)}   \right)^{1/2} + \left ( \dfrac{\partial_z m(z^*_0) }{c^2}\right )^2 \right)  + o(\eps) 
\end{array}
\end{equation}


\paragraph*{Description of the { Mean fitness} $\lambda_1$. }
The  equation~\eqref{eq_app:asex_U1} evaluated at the optimal trait $z=0$ yields
$
\lambda_1 = -  \partial^2_p H(p_0) \partial^2_z U_0(0)/2,
$ 
where $p_0 = \partial_z  U_0(0)$. To compute $\partial^2_z U_0(0),$ we  differentiate~\eqref{eq_app:U0} twice, and evaluate the expression at $z=0$: 
\begin{equation}\label{eq_app:d2U0_0}
   1 = \partial^2_p H  \left(p_0 \right) \left(\partial^2_z U_0(0)\right)^2\,.
\end{equation}
Recall that  $p_0 = \partial_v L(c)$. Moreover, since $\partial_p H$ and $\partial_v L$ are reciprocal functions, then the second derivatives are inverse from each other. Therefore $\partial^2_p H(p_0)  = \left( \partial^2_vL(c) \right)^{-1}$. Thus,  $\lambda_1$ is given by the following expression:
\begin{equation} \label{eq:lambda1}
\lambda_1 = - \dfrac12 \left( \dfrac{ 1}{\partial^2_v L(c)}   \right)^{1/2} \, .
\end{equation}

\paragraph*{Description of the evolutionary lag $z^*_1$. }
By pushing the computations further, it is also possible to derive the first order correction of the lag $z^*_1.$ It is defined such that $z^*_0 + \eps z^*_1$ is the critical point of $U_0 + \eps U_1$, that is $\partial_z (U_0 + \eps U_1)(z^*_0 + \eps z^*_1) = 0\, .$
By expanding this relation, but keeping only the first order terms, we obtain $ z^*_1 = -{\partial_z U_1(z^*_0)}/{ \partial^2_z U_0(z^*_0) } .$ On the other hand, evaluating the  equation~\eqref{eq_app:asex_U1} at $z=z^*_0$ yields $-{\partial_z U_1(z^*_0)}/{ \partial^2_z U_0(z^*_0) } = \lambda_1/(c\partial^2_zU_0(z^*_0))+1/(2c).$ Using the expression~\eqref{eq:app d2U0} of $\partial_z^2U_0(z^*_0)$, we obtain: 
\begin{align}\label{eq:z1star}
z^*_1 & = \dfrac1{2\partial_z \m(z^*_0)} \left( \dfrac{ 1}{\partial^2_v L(c)}   \right)^{1/2}  + \dfrac1{2c} \, .
\end{align}

\paragraph*{Description of the local shape. }\label{app:local_shape_U1}
We expand the second derivative of $U_0 + \eps U_1$ at the lag point $z^*_0+\epsilon z^*_1$ with respect to $\epsilon$ and we obtain
\begin{equation}\label{eq_app:expansion variance}
\partial^2_z (U_0 + \eps U_1)(z^*_0 + \eps z^*_1) = \partial^2_z U_0(z^*_0) + \eps \left( \partial^3_z U_0(z^*_0) z^*_1 +   \partial^2_z U_1(z^*_0)\right) + o(\eps)\, . 
\end{equation}
We aim at characterizing the term of order $\epsilon$ in this expansion. The first additional contribution $\partial^3_z U_0(z^*_0)$ can be easily deduced from the equation~\eqref{eq_app:U0} by differentiating it twice,  and evaluating at $z=z^*_0$:
\[ c \partial^3_z U_0(z^*_0) + \partial^2_z\m(z^*_0) =  \partial^2_p H  \left(0 \right) \left(\partial^2_z U_0(z_0^*)\right)^2 =  \left(\partial^2_z U_0(z^*_0)\right)^2 \, .
  \] 
The second additional contribution $ \partial^2_z U_1(z^*_0)$ is deduced from the   equation~\eqref{eq_app:asex_U1} by differentiating once and evaluating at $z = z^*_0$:
\[
c \partial^2_zU_1(z^*_0) =   \partial^2_z U_0(z^*_0) \partial_z U_1(z^*_0) - \dfrac12 \partial^3_z U_0(z^*_0) \, .
\]
Combining these two expressions with the expression~\eqref{eq:z1star} of $z_1^*$, and $\partial_z U_1(z_0^*)$, we get
\begin{align*}
      &\partial^3_z U_0(z^*_0) z^*_1 +   \partial^2_z U_1(z^*_0) 
     \\ & = \partial^3_z U_0(z^*_0)\left( z^*_1 -  \dfrac1{2c}\right)  + \dfrac1{c}\partial^2_z U_0(z^*_0) \partial_z U_1(z^*_0)  \\
      & = \frac1c \left ( \left(\partial^2_z U_0(z^*_0)\right)^2 -   \partial^2_z\m(z^*_0) \right ) \left(  \dfrac1{2\partial_z \m(z^*_0)} \left( \dfrac{ 1}{\partial^2_v L(c)}   \right)^{1/2}   \right)  - \dfrac1{c}\partial^2_z U_0(z^*_0)   \left (  \dfrac{\lambda_ 1}c  + \dfrac{\partial^2_z U_0(z^*_0)}{2c}  \right ) \\
      & = \frac1c \left ( \left(\dfrac{\partial_z \m(z^*_0)}c\right)^2 -   \partial^2_z\m(z^*_0) \right ) \left(  \dfrac1{2\partial_z \m(z^*_0)} \left( \dfrac{ 1}{\partial^2_v L(c)}   \right)^{1/2}   \right)  + \dfrac{\partial_z \m(z^*_0)}{c^2}    \left (  - \dfrac{1}{2c}\left( \dfrac{ 1}{\partial^2_v L(c)}   \right)^{1/2}   - \dfrac{\partial_z \m(z^*_0)}{2c^2}  \right ) \\
      & = - \dfrac{\partial^2_z\m(z^*_0)}{2c\partial_z \m(z^*_0)} \left( \dfrac{ 1}{\partial^2_v L(c)}   \right)^{1/2} -\frac12 \left ( \dfrac{\partial_z m(z^*_0) }{c^2}\right )^2
\end{align*}
This concludes the analysis of the corrector problem at first order.


\subsection{Numerical computation of the distributions $U_0$ and $U_1$ in the asexual model}\label{app:num_sol_asex_U0_U1}
The equation for $U_0$ \eqref{eq_app:U0} is a non linear Ordinary Differential Equation (ODE). It has a singular point at $z = 0$, where the function $p\mapsto cp - H(p)$ cannot be inverted. It was solved numerically in the following way: after differentiation with respect to $z$, equation~\eqref{eq_app:U0} becomes
\begin{equation*}
\left (\partial_p H(\partial_z U_0(z)) - c\right )  \partial_z^2 U_0(z) = \partial_z m(z) \quad \Leftrightarrow\quad  \dfrac{d}{dz}\left ( U_0'(z)\right ) = \dfrac{m'(z)}{\partial_p H(U_0'(z)) - c}\, .
\end{equation*}
This ODE on $U_0'(z)$ was solved using a classical solver (RK45), separately on the two branches $z>0$ and $z<0$. The issue is to initialize appropriately the solver for $z = 0^+$, and $z = 0^-$. The correct initialization was deduced from the analytical expressions of $U_0'(0)=p_0 = \partial_v L(c)$.

Next, the linear ODE for $U_1$ \eqref{eq_app:asex_U1} was computed along characteristic lines:
\begin{align*}
\dot \bz(s) = \partial_p H(\partial_z U_0(\bz(s))) - c\quad \Longrightarrow\quad  \dfrac{d}{ds} \left ( U_1(\bz(s))\right )&  = \lambda_1 + \dfrac12 \partial^2_p H(\partial_z  U_0(\bz(s))) \partial^2_z U_0(\bz(s)) \\
& = \lambda_1 + \dfrac12 \left (\dfrac{d}{dz} \partial_p H(\partial_z  U_0)  \right )(\bz(s)) \, . 
\end{align*} 
Integrating this formula with respect to  $s$ yields 
\begin{align*}
U_1(\bz(s)) - U_1(\bz(0)) & = \lambda_1 s + \frac12 \int_{0}^s \left (\dfrac{d}{dz} \partial_p H(\partial_z  U_0)  \right )(\bz(s'))\, ds' \\
& =  \lambda_1 s + \frac12 \int_{0}^{\bz(s)} \left (\dfrac{d}{dz} \partial_p H(\partial_z  U_0)  \right )(\bz) \left ( \dfrac{1}{\partial_p H(\partial_z U_0(\bz)) - c}\right )\, d\bz \\
& = \lambda_1 s + \frac12 \log\left|\dfrac{\partial_p H(\partial_z U_0(\bz(s))) - c}{\partial_p H(\partial_z U_0(\bz(0))) - c} \right|\, . 
\end{align*}
Again, the delicate issue is to evaluate appropriately the value $U_1(\bz(0))$ for a starting point $\bz(0)$ close to $0$ (notice that $0$ is an equilibrium point for the ODE: $\dot \bz(s) = \partial_p H(\partial_z U_0(\bz(s))) - c$). The correct approximation is given by the analytical expression of $\partial_z U_1(0)$ obtained by differentiating equation~\eqref{eq_app:asex_U1} with respect to $z$ and evaluating it at $z=0.$ 

\section{Qualitative properties of the standing variance at equilibrium $\Var(\Fb)$}\label{app:Var_c}
In this appendix, we discuss in detail the behavior of the standing variance at equilibrium with respect to the speed of change $\cb$ in the scenario of asexual reproduction.
Let us remind that in this case the standing variance at equilibrium is well approximated by the following expression at the leading order:
\[
   \Var(\Fb) \approx - \dfrac{\cb}{\partial_{\zb} \mb( \zb_0^*)}\,. 
\]
It is convenient to introduce the positive lag $|\zb^*_0|$, which is the distance to the optimal trait located at $\zb=0$, so that 
\[
   \Var(\Fb) \approx  \dfrac{\cb}{\partial_{\zb} \mb( |\zb_0^*|)} \, .
\]
Recall that the lag is deduced from the inversion of the increment of mortality $\mb$:
\begin{equation}\label{eq_app:lag}
   |\zb_0^*| = \mb^{-1}\left(\betab L\left(\dfrac{\cb}{\sigmab\betab}\right) \right),
\end{equation}
where $\mb^{-1}$ is the inverse of the function $\mb$ on $(0,\infty)$. 
The differentiation of the lag $|\zb_0^*|$ with respect to $\cb$ goes as follows:
\begin{equation}
\dfrac{d |\zb_0^*|}{d\cb}(\cb) = \dfrac1\sigmab \partial_v L\left(\dfrac{\cb}{\sigmab\betab}\right) \partial_\zb(\mb^{-1})\left(\betab L\left(\dfrac{\cb}{\sigmab\betab}\right) \right),
\end{equation}
Since $\partial_\zb(\mb^{-1}) = 1/\partial_\zb\mb(\mb^{-1})$, the previous expression becomes
\begin{equation}
\dfrac{d |\zb_0^*|}{d\cb}(\cb) = \dfrac1\sigmab \partial_v L\left(\dfrac{\cb}{\sigmab\betab}\right) \dfrac1{\partial_\zb\mb\left(\mb^{-1}  \left(\betab L\left(\dfrac{\cb}{\sigmab\betab}\right) \right)\right)} = \dfrac1\sigmab \partial_v L\left(\dfrac{\cb}{\sigmab\betab}\right) \dfrac1{\partial_\zb\mb(|\zb_0^*|)},
\end{equation}
Reformulating this expression, we get an alternative expression for the variance:
\begin{equation}
    \Var(\Fb) \approx  \dfrac{\cb}{\partial_{\zb} \mb( |\zb_0^*|)} = \dfrac{d |\zb_0^*|}{d\cb}(\cb) \times  \sigmab\cb \left (\partial_v L\left ( \frac{\cb}{\sigmab\betab}\right )\right )^{-1}
\end{equation}
Now let differentiate the latter expression with respect to $\cb$:
\begin{equation}
   \dfrac{d}{d\cb} \left (\dfrac{\cb}{\partial_{\zb} \mb( |\zb_0^*|)}\right) = 
   \dfrac{d^2 |\zb_0^*|}{d\cb^2}(\cb) \times  \sigmab\cb \left (\partial_v L\left ( \frac{\cb}{\sigmab\betab}\right )\right )^{-1} + \dfrac{d |\zb_0^*|}{d\cb}(\cb) \times   \sigmab\left (\partial_vL\left ( \frac{\cb}{\sigmab\betab}\right )\right )^{-1} \left( 1 - \dfrac{\cb}{\sigmab\betab} \dfrac{\partial^2_vL\left ( \frac{\cb}{\sigmab\betab}\right )}{\partial_vL\left ( \frac{\cb}{\sigmab\betab}\right )}  \right)
\end{equation}
We shall establish that for all $\cb>0$, the following inequality holds true:
\[
   \left( 1 - \dfrac{\cb}{\sigmab\betab} \dfrac{\partial^2_vL\left ( \frac{\cb}{\sigmab\betab}\right )}{\partial_vL\left ( \frac{\cb}{\sigmab\betab}\right )}  \right) \geq0\,.
\]
Indeed, it can be reformulated by means of $p$ such that $p = \partial_vL\left ( {\cb}/{\sigmab\betab}\right )$, as follows:
\begin{equation}
     1 - \dfrac{\cb}{\sigmab\betab} \dfrac{\partial^2_vL\left ( \frac{\cb}{\sigmab\betab}\right )}{\partial_vL\left ( \frac{\cb}{\sigmab\betab}\right )}   = 1 - \dfrac{\partial_p H(p)}{p \partial^2_p H(p)} = 1 - \dfrac{\displaystyle\int_{\R} y K(y) e^{py}\, dy}{\displaystyle p \int_{\R} y^2 K(y) e^{py}\, dy} =
     1 - \dfrac{\displaystyle\int_{\R_+} y K(y) \sinh(py)\, dy}{\displaystyle p \int_{\R_+} y^2 K(y) \cosh(py)\, dy}\,. 
\end{equation}
The conclusion follows from the pointwise inequality $\tanh(py) \leq py$ for $p,y\geq 0$, which is equivalent to $\sinh(py) \leq py \cosh(py)$. 

On the other hand, we have shown that the lag  increases with respect to the speed of change $\cb,$ thus $d |\zb_0^*|/d\cb \geq 0$. Then, if the lag is convex with respect to the speed of change $\cb$, that is $d^2 |\zb_0^*|/d\cb^2 \geq 0$, then the standing variance at equilibrium  increases with respect to the speed $\cb$. 

However, the convexity of the lag depends on the convexity of the function $c \mapsto\mb^{-1}(\betab L(c))$. If the selection is quadratic $\mb(\zb)=\alphab\zb^2$, this function is concave for any mutation kernel. However, if the selection function is more than quadratic, we can find mutation kernels such that the lag becomes convex. 

In the diffusion approximation $L(c)=c^2/2$, we can go further. In this case, we know from equation~\eqref{eq:weak_selec_zstar0} that the lag accelerates with $\cb$ if $\mb$ is sub-quadratic. Whereas it lag is concave if $\mb$ is super-quadratic in the sense of~\eqref{eq:strong_selec_zstar0}. 

As a result, we have shown that the variance $\Var(\Fb)$ increases with $\cb$ if the function $c \mapsto\mb^{-1}(\betab L(c))$ is convex. More precisely, in the diffusion approximation, the variance increases with $\cb$ if  $\mb$ is sub-quadratic in the sense of~\eqref{eq:weak_selec_zstar0}.


\section{Sexual type of reproduction (details of Section \ref{sec:sex})}\label{app:sexual}
In this section we develop the computations required to describe  $U$ up to order $\eps^2$, as in \eqref{eq:ansatz sexual}. We present arguments from convex analysis to characterize $U_0$. We provide an explicit formula for   the first order correction $U_1$ as an infinite series. Meanwhile, we present tedious computations needed to identify the linear part of  $U_1$, and we derive the first order correction of the { mean fitness} $\lambda_1$ as a by-product. 

Our starting point is the following relationship which is equivalent to finding a stationary density in the moving frame, expanded at first order in  $\eps^ 2$:
\begin{multline}\label{eq:log transform sexue}
\lambda_0 + c\partial_z U_0(z) + \m(z) =\\
\f{  \displaystyle 
\dfrac{1}{\eps^2 \sqrt{2}\pi } \iint_{\R^2} \exp\left(-\dfrac1{\eps^2} \left[   \left(z-\dfrac{z_1+z_2}{2}\right)^2 + U_0(z_1) + U_0(z_2)  - U_0(z) \right] - U_1(z_1) - U_1(z_2) + U_1(z) \right) d z_1 d z_2}{\displaystyle \dfrac{1}{\eps \sqrt{2\pi} }  
  \int_\R \exp\left(-\dfrac{1}{\eps^2}U_0(z') - U_1(z') \right) d z'}\, . 
\end{multline}
Note that the prefactors (involving $\eps, \pi$ have been arranged for the sake of normalizing singular integrals). 

The arguments below are formal computations. We refer to \citep{CalGarPat19} for a rigorous analysis of this asymptotic analysis in the case $c=0$, and to \citep{Pat20} for the time marching problem.

\subsection{The characterization of $U_0$ by convex analysis}\label{app:sex_U0}

Recall that the identity satisfied by $U_0$ is the following one, ensuring that the right hand side of \eqref{eq:log transform sexue} does not get trivial as $\eps \to 0$: 
%
\begin{align}
&\forall z\in \R \quad \min_{(z_1,z_2) \in \R^2} \left[  \left(z-\dfrac{z_1+z_2}{2}\right)^2 + U_0(z_1) + U_0(z_2)  - U_0(z) - \min U_0 \right] = 0\nonumber\\
&\Longleftrightarrow\quad U_0(z) + \min U_0 = \min_{(z_1,z_2) \in \R^2}  \left(  \left(z-\dfrac{z_1+z_2}{2}\right)^2 + U_0(z_1) + U_0(z_2)   \right) \label{eq_app:function ineq U0}
\, .
\end{align}

The goal of this section is to prove that any solution of the functional equation \eqref{eq_app:function ineq U0} is given by a member of the three parameters family 
\begin{equation}\label{eq_app:3 param family}
U_0(z) = C+ \frac{(z-a)_-^2}{2} + \frac{(z-b)_+^2}{2}\, ,   
\end{equation} 
where the parameters $a,b$ are such that $a\leq b$ and $C$ is an arbitrary constant. We denote by  $z_0^*$ a minimum point of $U_0$. We can restrict to $\min U_0 = 0$ without loss of generality (so that the additive constant $C$ is set to $0$).  The characterization of $U_0$ is done in several steps. 

\paragraph*{Regularity and $\lambda-$concavity.}
Firstly, notice that $U_0(z) - z^2$ is a concave function, as it can be written as
\begin{align*}
U_0(z) - z^2 & = \min_{(z_1,z_2) \in \R^2}  \left(  - z(z_1 + z_2) + \left(\dfrac{z_1+z_2}{2}\right)^2 + U_0(z_1) + U_0(z_2)   \right)\\
& = \min \left \{ \text{affine functions with respect to $z$} \right\}\, . 
\end{align*}
We deduce that $U_0$ is continuous, and that it admits left and right derivatives everywhere. 

\paragraph*{The  convex conjugate.}
The trick is to introduce the convex conjugate $\hat U_0$ (also called the Legendre transform of $U_0$): 
\begin{equation*}
\hat U_0(y) = \max_{z\in \R} ( (z-z_0^*) y - U_0(z))\,,
\end{equation*}
where $z_0^*$ is a minimum point of $U_0$. 
The basic properties of $\hat U_0$ are listed below:
\begin{itemize}
\item $\hat U_0$ is convex, so it is continuous, and it admits left and right derivatives everywhere, 
\item $\hat U_0(0) = \max \left (-U_0\right ) = - \min \left (U_0\right ) = 0$,
\item for all $y$, $\hat U_0(y) \geq -U_0(z_0^*) = 0$, thus $\min \hat U_0 = 0$. 
\end{itemize}
We deduce from the functional identity \eqref{eq_app:function ineq U0}, that
\begin{align*}
\hat{U}_0(y) &= \max_{z\in \R} \left(  (z-z_0^*)y -  \min_{\substack{(z_1,z_2)\\ \in \R^2}}  \left(  \left(z-\dfrac{z_1+z_2}{2}\right)^2 + U_0(z_1) + U_0(z_2)   \right)\right) \\
& = \max_{ (z,z_1,z_2)  \in \R^3} \left( (z-z_0^*)y -  \left(z-\dfrac{z_1+z_2}{2}\right)^2 - U_0(z_1)- U_0(z_2)   \right) \\
& = \max_{ (z_1,z_2)  \in \R^2}  \left( \max_{z\in\R} \left( (z-z_0^*)y -  \left(z-\dfrac{z_1+z_2}{2}\right)^2\right) - U_0(z_1)- U_0(z_2)   \right) \\
& = \max_{ (z_1,z_2)  \in \R^2}  \left(  \frac{y^2}{4} + \frac{1}{2}\left(z_1+z_2\right) y -z_0^*y - U_0(z_1)- U_0(z_2)   \right) \\
& = \frac{y^2}{4} + \max_{z_1 \in \R} \left(\frac{1}{2} (z_1-z_0^*) y   - U_0(z_1)  \right)
+  \max_{z_2 \in \R} \left(\frac{1}{2} (z_2-z_0^*) y   - U_0(z_2)  \right)\, .
\end{align*}
Finally, we end up with the following functional identity,
\begin{equation}\label{eq_app:hat U0}
\hat{U}_0(y) = \frac{y^2}{4}  + 2\hat{U}_0\left(\frac{y}{2}\right)\, .
\end{equation}
We observe that $\hat{U}_0(y) = y^2/2$ is a solution to the latter identity. However, it is not the only one. More generally, let $a = \hat{U}_0'(0^-) $ and $b = \hat{U}_0'(0^+) $ denote the left and the right derivative at $y = 0$, respectively. By convexity, and optimality at the origin $y=0$ (namely, $\min \hat U_0 = \hat U_0(0) = 0$), we have $a\leq 0\leq b$.  
We deduce recursively from \eqref{eq_app:hat U0} the series expansion
\begin{align}
\hat{U}_0(y) &= \dfrac{y^2}4 + \dfrac{y^2}8 + \dfrac{y^2}{16} + \dots + 2^n \dfrac{(2^{-n}y)^2}4 + 2^{n+1} \hat{U}_0\left(2^{-(n+1)}y\right)\,,\nonumber\\
\Longrightarrow \hat{U}_0(y) & = \dfrac{y^2}{2} +  \hat{U}_0'(0^{\pm})y\, . \label{eq:app56}
\end{align}
Obviously, the choice of the left or right derivative depends on the sign of $y$. 

\paragraph*{The  convex bi-conjugate.} Next, we define the convex bi-conjugate 
\begin{equation*}
 \breve{U}_0(z) =  \max_{y\in \R} \left  ((z - z_0^*) y - \hat{U}_0(y)\right )\,.
\end{equation*}
Standard results in convex analysis states that $\breve{U}_0$ and $U_0$ coincide if $U_0$ is convex. More generally, $\breve{U}_0$ is the (lower) convex envelope of $U_0$ \citep{Roc70}. This is quite useful, because the characterization \eqref{eq:app56} enables to compute the convex bi-conjugate:
\begin{equation}\label{eq_app:3 param family app}
\breve{U}_0(z) = \frac{(z - z_0^*-a)_-^2}{2} + \frac{(z - z_0^*-b)_+^2}{2}\, . 
\end{equation}  
We deduce that the latter function is the (lower) convex envelope of $U_0$. The last (delicate) step consists in proving that it coincides with $U_0$. 

\paragraph*{From the convex envelope to the function.}
The idea is to use the functional identity \eqref{eq_app:function ineq U0} iteratively. As $z = z_0^* + a$ is an extremal point of the graph of $\breve{U}_0$, the values of $U_0$ and $\breve{U}_0$ must coincide at this point. Hence, we have $U_0(z_0^*+ a) = 0$, and similarly $U_0(z_0^*+b) = 0$. Recall that $U_0(z_0^*) = 0$ by definition. As a consequence, we have for $z_1 = z_0^* + a$, $z_2 = z_0^*$, and $z = z_0^* + a/2$ in \eqref{eq_app:function ineq U0}:
\begin{equation*}
U_0\left ( z_0^* + \dfrac{a}{2} \right ) \leq 0\, ,
\end{equation*} 
from which we deduce that $U_0$ vanishes at $z = z_0^* + a/2$ as well, and similarly at $z = z_0^* + b/2$. The same argument shows that $U_0$ vanishes at each middle point between two vanishing points. So, it vanishes on a dense set of points in $z_0^* + (a,b)$. By continuity of $U_0$, it vanishes everywhere on $z_0^* + [a,b]$. Finally, it coincides with its (lower) convex envelope \eqref{eq_app:3 param family app} because the latter is strictly convex outside the interval $[a,b]$. 

Finally, it is necessary that $a = b = 0$ in the present context. Otherwise $F$ would not correspond to a population density uniformly with respect to vanishing $\eps$. 

We have proved that $U_0$ is necessary of the form 
\begin{equation}\label{eq_app:U0 polynom}
U_0(z) = \dfrac{(z-z_0^*)^2}{2} \, .
\end{equation} 
However, we are not able at this point to characterize the evolutionary lag $z^*_0.$ We need to push the analysis beyond the first order and compute the profile $U_1$, as done in the following sections.

\paragraph*{Discussion.}
There is an immediate interpretation of this result: we found that the equation is dominated by the reproduction term in the regime of small variance. Hence, the stationary distribution at the leading order equilibrium is the Gaussian distribution with prescribed variance (here, renormalized to a unit value), meaning a quadratic polynomial after taking the logarithm. In fact, Gaussian distributions are known to be stationary distributions of the Infinitesimal model in the absence of selection. As selection does not act on reproduction, there is no way to find the evolutionary lag at equilibrium, and so $z_0^*$ must be unknown at this point of analysis. The situation is quite different from the case of asexual reproduction, where no stationary distribution can be achieved without selection, and the evolutionary lag is deduced from the knowledge of $U_0$, accordingly. 

\subsection{Description of the corrector $U_1$}\label{app:sex_U1}


Next, we can rearrange the right hand side in \eqref{eq:log transform sexue} using the characterization of $U_0$ \eqref{eq_app:U0 polynom}. 
It is instructive to begin with the denominator integral, which is a classical computation:
\begin{align*}
\dfrac{1}{\eps \sqrt{2\pi}}  \int_\R \exp\left ( -\frac{(z' - z_0^*)^2}{2\eps^2} \right )\exp\big(- U_1(z')\big) d z'
&= \dfrac{1}{
\sqrt{2\pi}} \int_\R \exp\left (-\frac{y'^2}{2}\right ) \exp\left( - U_1(z_0^* + \eps y') \right) d y'\\ 
&\underset{\eps\to 0}{\longrightarrow}    \exp(- U_1(z_0^*) )\, .
\end{align*}
Indeed, the function $(\eps \sqrt{2\pi})^{-1}\exp\left (- (z' - z_0^*)^2/(2\eps^2)\right )$ is the approximation of a Dirac mass as $\eps\to 0$. Hence the integral concentrates on the evolutionary lag $z_0^*$: this yields the convergence of the integral towards $\exp(- U_1(z_0^*) )$. An alternative way to say  is that, in the integral $\int F(z')\,dz'$, most of the contribution comes from those $z'$ which are close to~$z_0^*$.

\subsubsection{What are the most representative parental traits?}\label{eq:seeking U1}


The same kind of computation allows  handling the numerator in \eqref{eq:log transform sexue}. 
The key point is to understand how the term inside the integral gets concentrated as $\eps \to 0$. In other words, we shall identify what are the most representative traits $(z_1,z_2)$ of parents giving birth to an offspring of trait $z$. Those will contribute mainly to the integral in the right hand side. They will enable to derive the equation for $U_1$.

A preliminary computation is required:  the double integral gets concentrated at the minimum points  (with respect to variables $(z_1,z_2)$) of the  quadratic form under brackets:
\begin{equation} 
\left(z-\dfrac{z_1+z_2}{2}\right)^2 + U_0(z_1) + U_0(z_2)  - U_0(z)\quad \hbox{ where } \quad U_0(z) = \dfrac{(z-z_0^*)^2}{2} \,.\label{eq:trade off}
\end{equation}
We know already that the minimum value is zero thanks to the characterization \eqref{eq_app:function ineq U0}. The values above the minimum will contribute very little to the integral as they will have size of order $\exp(-\delta/\eps^2)$, for $\delta>0$. Indeed, this decays to zero very fast as $\eps\to 0$. 

\begin{figure}
\begin{center}
   \includegraphics[width=0.8\linewidth]{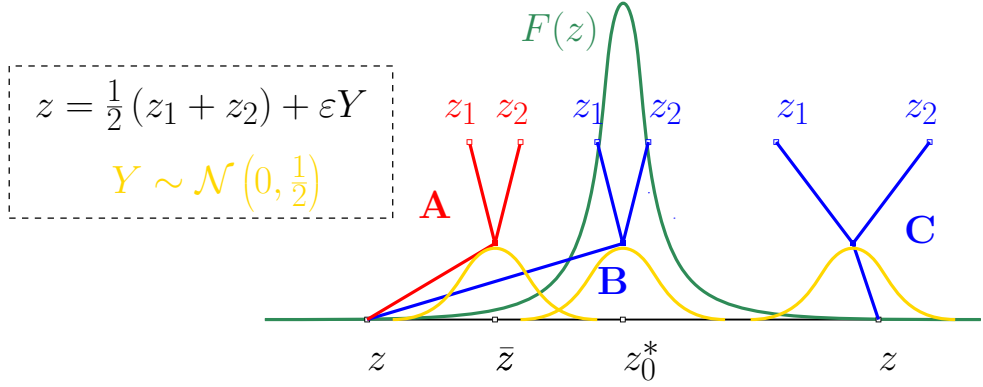}
   \caption{Sketch of the argument that underpins the estimation of the double integral in \eqref{eq:log transform sexue}. Recall that the infinitesimal model assigns to an offspring the trait $z$ which is the mean value of the parental traits plus a normal random variable with standard deviation $1/\sqrt{2}$ (in adimensional variables). Among the three scenarios $\mathbf{A},\mathbf{B},\mathbf{C}$, the first one is by far the most likely in the regime of small variance $\eps^2 \ll 1$. In scenario $\mathbf{B}$, the parental traits $(z_1,z_2)$ are close to the evolutionary lag $z_0^*$: this is a likely event from the point of view of the parental trait distribution. However, it is very unlikely to draw a random number $Y$ so large resulting in $z$ at the next generation. In scenario $\mathbf{C}$, the deviation is small, so that the mean parental trait is close to $z$: this is a likely event from the point of view of the "choice" of the offspring trait. However, it is very unlikely to draw a parent with trait $z_2$ 
from the phenotypic distribution $F$: that one is too far from the evolutionary lag in the tail of the distribution. Scenario $\mathbf{A}$ is the compromise between these two antagonistic effects.}\label{fig:sketch mid point}
\end{center}
\end{figure}

Direct computation provides the unique minimum $(z_1,z_2) = (\bar z, \bar z)$, with $\bar z = (z+z^*_0)/2.$
This means that an offspring of trait $z$ is very likely to be the combination of equal parental traits $z_1 = z_2$, equal to the mid-value between $z$ and the evolutionary lag $z_0^*$. This is the result of an interesting trade-off: parents with phenotype close to the evolutionary lag value $z_0^*$ are more frequent but the chance of producing an offspring with phenotype $z$ decreases when their own phenotype departs from the latter value. As a compromise, the most likely configuration  is when both parents have the mid-point trait $\bar z$, see Figure \ref{fig:sketch mid point}. 

We thus define the following change of variable centered around this minimum point:
\begin{equation} \label{eq:change var z}
\begin{cases}
z_1 = \bar z + \eps  y_1  \\
z_2 = \bar z + \eps  y_2
\end{cases}
\end{equation}
The quadratic form between brackets $[\cdots ]$ in the numerator of \eqref{eq:log transform sexue} is transformed into an expression which does not depend on $\eps$:
\begin{equation}
\frac{1}{\eps^2}\left[\left(z-\dfrac{z_1+z_2}{2}\right)^2 + U_0(z_1) + U_0(z_2)  - U_0(z)\right] 
=  \frac{1}{2} y_1 y_2 + \frac34\left(y_1^2 + y_2^2\right)\,.
\end{equation}
And the numerator finally writes
\begin{align*}
&\dfrac{1}{
\sqrt{2}\pi}  \iint_{\R^2} \exp\left(-  \left[   \frac{1}{2} y_1 y_2 + \frac34\left(y_1^2 + y_2^2\right) \right] - U_1(\bar z + \eps y_1) - U_1(\bar z + \eps y_2) + U_1(z) \right) d y_1 d y_2 \\
&  \underset{\eps\to 0}{\longrightarrow} \dfrac{1}{
\sqrt{2}\pi} \left(  \iint_{\R^2} \exp\left(-  \left[   \frac{1}{2} y_1 y_2 + \frac34\left(y_1^2 + y_2^2\right) \right] \right)  d y_1 d y_2\right) \exp\left( - U_1(\bar z ) - U_1(\bar z ) + U_1(z) \right)\\
& =  \exp\left(  -2 U_1(\bar z ) + U_1(z) \right)
\end{align*}
Note that the prefactor $(\sqrt{2}\pi)^{-1}$ is such that the integral in $(y_1,y_2)$ has unit value.  

\subsubsection{Equation for the corrector $U_1$}

We conclude that equation \eqref{eq:log transform sexue} converges as $\epsilon\to0$ to the following equation on the corrector $U_1$:
\begin{equation}\label{eq:U1}
\lambda_0 + c (z - z_0^*)  + \m(z) = \exp\left( U_1(z_0^*) -2 U_1(\bar z ) + U_1(z) \right)\, , \quad \text{with}\quad \bar z = \frac{z + z_0^*}{2}\, .
\end{equation}
This equation is simple enough to admit an explicit solution as an infinite series, as shown below. 

Note that the values of $\lambda_0$ and $z_0^*$ can be deduced readily from \eqref{eq_app:calcul U1 1} as explained in the main text \eqref{eq:sex_lambda_z}. 

\subsubsection{Analytical expression of $U_1$}
It is convenient to reformulate equation  \eqref{eq:U1} as follows,
by using the formula~\eqref{eq:sex_lambda_z} for $\lambda_0$ and $z^*_0,$
\begin{equation}\label{eq_app:calcul U1 1}
\log\big(1 + \G(z)\big) = U_1(z_0^*) -2 U_1\left( \frac{z + z_0^*}{2} \right) + U_1(z)\, 
\end{equation} 
where $\G(z) = m(z) - \partial_z m(z_0^*)(z-z_0^*)$ is such that $\G(0) = \partial_z \G(0) = 0$.
Differentiating this equation with respect to $z$, we obtain 
\begin{equation*}
\dfrac{\partial_z \G(z)}{1 + \G(z)} = \partial_z U_1(z) -\partial_z U_1\left( \frac{z + z_0^*}{2} \right)  \, .
\end{equation*}
After the change of variable $z = z_0^* + h$, we get eventually the recursive relation where the value at some $z_0^* + h$ can be computed from the value at $z_0^* + h/2$,
\begin{equation*}
\partial_zU_1(z_0^* + h) = \partial_zU_1\left(z_0^* + \frac h2\right) + \dfrac{\partial_z \G(z_0^* + h)}{1 + \G(z_0^* + h)}\, .
\end{equation*}
We deduce the following series expansion,
\begin{equation}\label{eq_app:dzU1}
\partial_zU_1(z_0^* + h) = \partial_zU_1(z_0^*) + \sum_{n=0}^{\infty}  \dfrac{\partial_z \G(z_0^* + 2^{-n}h)}{1+ \G(z_0^* + 2^{-n}h)}\, .
\end{equation}
This provides an expression for $U_1$ after integration with respect to $h$,
\begin{equation}\label{eq_app: U1 series}
 \quad U_1(z_0^*+h) = U_1(z_0^*) + h \partial_z U_1(z_0^*) + \sum_{n=0}^{\infty} 2^n \log \left( 1+ \G( z_0^*+ 2^{-n}h )\right)\, .
\end{equation}
There are two degrees of freedom in the above expression of $U_1$. First, the constant part $U_1(z_0^*)$ cannot be determined, because $U$ is defined up to an additive constant. Thus, we are free to choose any value for $U_1(z_0^*)$, say $U_1(z_0^*)=0$ for instance. On the other hand, the value $p^* = \partial_z U_1(z_0^*)$ plays a key role in the shape of the distribution, related to the expansion of the evolutionary lag, see \eqref{dzU1 1} below, but its value cannot be elucidated at this stage. We need to push the expansion up to order $\eps^4$  to get the following formula for $p^*$:
\begin{equation} \label{eq_app:p*}
p^* = \dfrac{\partial^3_z \m(z_0^*)}{2 \partial^2_z \m(z_0^*) } + 2c \, ,
\end{equation}
see next section for the complete computation (see also \citep{CalGarPat19} for an alternative path with limited expansions to the next order in the case $c=0$).

We deduce the following expression for $U_1$,
\begin{equation}\label{eq_app:U1 complete}
 \quad U_1(z_0^*+h) =   \left(  \dfrac{ \partial^3_z \m(z_0^*)}{ 2 \partial^2_z \m(z_0^*)}  + 2c\right) h  + \sum_{n=0}^{\infty} 2^n \log \left( 1 + \G( z_0^*+ 2^{-n}h )\right)\,.
\end{equation}

\subsubsection{The missing linear part: calculation of $\partial_zU_1(z^*_0)$}\label{app:p*}

Starting with the equation satisfied by $U$ \eqref{eq:U sexual}, and plugging the ansatz 
\begin{equation}\label{eq_app:ansatz sexual} 
\begin{cases} 
U(z) = U_0(z) + \eps^2 U_1(z) +\epsilon^4 U_2(z)+ o(\eps^4) \\
\lambda = \lambda_0 + \eps^2 \lambda_1 +\epsilon^4\lambda_2  + o(\eps^4)
\end{cases}
\end{equation}
we obtain the following equation up to order $\eps^2$:
\begin{align*}
&\hspace*{-50pt}\lambda_0 + \eps^2 \lambda_1 + c\partial_z U_0(z) + \eps^2c\partial_z U_1(z) + \m(z) \\
&\hspace*{-50pt} = \dfrac{ \dfrac{1}{\sqrt{2}\pi} \displaystyle \iint_{\R^2} \exp\left(-  \left[   \frac{1}{2} y_1 y_2 + \frac34\left(y_1^2 + y_2^2\right) \right] - U_1(\bar z + \eps y_1) - U_1(\bar z + \eps y_2) - \eps^2 U_2(\bar z + \eps y_1) - \eps^2 U_2(\bar z + \eps y_2)  + U_1(z) + \eps^2 U_2(z) \right) d y_1 d y_2  }{\displaystyle\dfrac{1}{\sqrt{2\pi}}\int_\R \exp\left(-\dfrac{y'^2}{2} - U_1(z_0^* + \eps y') - \eps^2 U_2( z_0^* + \eps y') \right) d y'  }
\end{align*}
The integrals were subject to the same change of variables as in \eqref{eq:change var z}.
After elimination of higher order contributions, we obtain for the denominator, up to order $\eps^2$:
\begin{align*}
&\dfrac{1}{\sqrt{2\pi}}\int_\R \exp\left(-\dfrac{y'^2}{2} - U_1(z_0^* + \eps y') - \eps^2 U_2( z_0^* + \eps y') \right) d y' \\
&=
\dfrac{1}{\sqrt{2\pi}}\int_\R \exp\left(-\dfrac{y'^2}{2} - U_1(z_0^*) - \eps y' \partial_z U_1(z_0^*) - \eps^2\frac{y'^2}{2} \partial^2_z U_1(z_0^*)  - \eps^2 U_2( z_0^* ) \right) d y' \\
& = \dfrac{1}{\sqrt{2\pi}}\int_\R \exp\left(-\dfrac{y'^2}{2} - U_1(z_0^*) \right)\left( 1 - \eps y' \partial_z U_1(z_0^*) + \frac{\eps^2}{2}y'^2\left|\partial_z U_1(z_0^*)\right|^2 - \frac{\eps^2}{2}y'^2 \partial^2_z U_1(z_0^*)  - \eps^2 U_2( z_0^* ) \right) d y' \\
& =   \exp\left( - U_1(z_0^*) \right)  \left( 1 + \frac{\eps^2}{2} \left|\partial_z U_1(z_0^*)\right|^2 -  \frac{\eps^2}{2} \partial^2_z U_1(z_0^*)  - \eps^2 U_2( z_0^* ) \right) \, .
\end{align*}
In an analogous way, we obtain for the numerator,
\begin{tiny}
\begin{align*}
&\hspace*{-50pt} \dfrac{1}{\sqrt{2}\pi} \iint_{\R^2} \exp\left(-  \left[   \frac{1}{2} y_1 y_2 + \frac34\left(y_1^2 + y_2^2\right) \right] - U_1(\bar z + \eps y_1) - U_1(\bar z + \eps y_2) - \eps^2 U_2(\bar z + \eps y_1) - \eps^2 U_2(\bar z + \eps y_2)  + U_1(z) + \eps^2 U_2(z) \right) d y_1 d y_2 \\
&\hspace*{-50pt}=
\dfrac{1}{\sqrt{2}\pi}
 \iint_{\R^2} \exp\left(-  \left[   \frac{1}{2} y_1 y_2 + \frac34\left(y_1^2 + y_2^2\right) \right] - 2 U_1(\bar z) + U_1(z)\right)\left( 1  -\eps \left[ y_1 + y_2 \right]\partial_z U_1(\bar z) + \dfrac{\eps^2}{2} \left[ y_1 + y_2 \right]^2\left|\partial_z U_1(\bar z)\right|^2 + \dfrac{\eps^2}2\left[ y_1^2 + y_2^2 \right]\partial^2_z U_1(\bar z) - 2 \eps^2 U_2(\bar z) + \eps^2 U_2(z) \right) d y_1 d y_2\\
&\hspace*{-50pt}=  \exp\left(- 2 U_1(\bar z) + U_1(z)  \right) \left(  1 + \dfrac{ \eps^2}2 \left|\partial_z U_1(\bar z)\right|^2  -  \dfrac{3\eps^2}4 \partial^2_z U_1(\bar z) - 2 \eps^2 U_2(\bar z) + \eps^2 U_2(z)  \right)
\end{align*}
\end{tiny}
Combining all these expansions, we obtain up to order $\eps^2$:
\begin{align*}
&\hspace*{-50pt}\lambda_0 + \eps^2 \lambda_1 + c\partial_z U_0(z) + \eps^2c\partial_z U_1(z) + \m(z) \\
&\hspace*{-50pt}= 
 \exp\left( U_1(z_0^*) -2 U_1(\bar z ) + U_1(z) \right) \dfrac{ 1 + \dfrac{ \eps^2}2 \left|\partial_z U_1(\bar z)\right|^2 -  \dfrac{3\eps^2}4 \partial^2_z U_1(\bar z) - 2 \eps^2 U_2(\bar z) + \eps^2 U_2(z) }{\displaystyle 1 + \frac{\eps^2}{2} \left|\partial_z U_1(z_0^*)\right|^2 - \frac{\eps^2}{2} \partial^2_z U_1(z_0^*)  - \eps^2 U_2( z_0^* )}\\
&\hspace*{-50pt} = \exp\left( U_1(z_0^*) -2 U_1(\bar z ) + U_1(z) \right) \left( 1 + \eps^2 \left( \dfrac{ 1}2 \left|\partial_z U_1(\bar z)\right|^2 - \frac{1}{2} \left|\partial_z U_1(z_0^*)\right|^2 + \frac{1}{2} \partial^2_z U_1(z_0^*) -  \dfrac{3}4 \partial^2_z U_1(\bar z) + U_2(z_0^*) -2 U_2(\bar z ) + U_2(z)  \right) \right)\, .
\end{align*}
By identifying contributions of order $\eps^2$ on both sides, we deduce the following equation for the next order correction $U_2$,
\begin{equation*}
U_2(z_0^*) -2 U_2(\bar z ) + U_2(z) =  \frac{1}{2} \left|\partial_z U_1(z_0^*)\right|^2  - \dfrac{ 1}2 \left|\partial_z U_1(\bar z)\right|^2 +  \dfrac{3}4 \partial^2_z U_1(\bar z) -  \frac{1}{2} \partial^2_z U_1(z_0^*) + \dfrac{\lambda_1 + c\partial_z U_1(z)}{1 + \G(z)}\, .
\end{equation*}
By evaluating, and differentiating at $z = z_0^*$, we deduce the following pair of identities,
\begin{equation}\label{eq_app:lambda1_sex}
\begin{cases}
0 = \dfrac14 \partial^2_z U_1(z_0^*) + \lambda_1 + c \partial_z U_1(z_0^*)\medskip\\
0 = -\dfrac12 \partial^2_z U_1(z_0^*) \partial_z U_1(z_0^*)  + \dfrac38 \partial^3_z U_1(z_0^*) + c \partial^2_z U_1(z_0^*)
\end{cases}
\end{equation}
The second identity enables to compute $p^* = \partial_z U_1(z_0^*)$:
\begin{equation} \label{eq_app:pstar}
p^* = \dfrac{3 \partial^3_z U_1(z_0^*)}{ 4 \partial^2_z U_1(z_0^*)} + 2c
= \dfrac{\partial^3_z \m(z_0^*)}{2 \partial^2_z \m(z_0^*) } + 2c \, ,
\end{equation}
where $ \partial^2_z U_1(z_0^*)$ and  $\partial^3_z U_1(z_0^*)$ are deduced straightforwardly from equation \eqref{eq_app:calcul U1 1} after multiple differentiation, or directly from \eqref{eq_app:dzU1}. This yields the missing part in \eqref{eq_app:U1 complete}. 

\subsubsection{Analytical expressions of the macroscopic corrections terms $\lambda_1$ and $z_1^*$}
\label{app:lambda1}

\paragraph*{Description of Malthus  rate $\lambda_1$.}  The   first identity in \eqref{eq_app:pstar}  provides $\lambda_1 = -\partial^2_z U_1(z_0^*)/4 - c \partial_z U_1(z_0^*)$.
The expression~\eqref{eq:U1} differentiated twice and evaluated at $z=z_0^*,$ yields $\partial^2_z U_1(z_0^*) =  2\partial^2_z \m(z_0^*).$ We conclude from the expression of $p^*$ that 
\begin{equation}
   \lambda_1  = - 2c^2 - c\dfrac{ \partial^3_z \m(z_0^*)}{ 2 \partial^2_z \m(z_0^*)} - \frac12 \partial^2_z \m(z_0^*).
\end{equation}

\paragraph*{Description of the evolutionary lag correction $z^*_1$.} 
The first order correction of the evolutionary lag $z^*_1$ is defined such that $z^*_0 + \eps z^*_1$ is the critical point of $U_0 + \eps^2 U_1$, that is $\partial_z (U_0 + \eps U_1)(z^*_0 + \eps z^*_1) = 0\, .$
Expanding this relation and keeping only the terms of order $\epsilon^2,$ we obtain using the expression of $p^*$,
\begin{equation}\label{dzU1 1} 
 z^*_1 = -\partial_z U_1(z^*_0) = -\dfrac{\partial^3_z \m(z_0^*)}{2 \partial^2_z \m(z_0^*) } - 2c \, .  
\end{equation}

\paragraph*{Description of the local shape.} 
The second derivative of $U_0 + \eps^2 U_1$ at the evolutionary lag $z^*$  is equal to $
\partial^2_z (U_0 + \eps^2 U_1)(z^*_0 + \eps^2 z^*_1) = \partial^2_z U_0(z^*_0) + \eps^2 \left( \partial^3_z U_0(z^*_0) z^*_1 +   \partial^2_z U_1(z^*_0)\right)$, up to the order $\epsilon^2$. 
Since $\partial^3_z U_0$ is equal to 0, we can deduce from the expression of $U_1$ that the local shape around $z^*$ is given by
\begin{equation*}
\partial^2_z (U_0 + \eps^2 U_1)(z^*_0 + \eps^2 z^*_1) = 1 + 2\eps^2\partial_z^2m(z^*_0) \, . 
\end{equation*}

\section{Numerical computation of the equilibrium $(\lambdab,\Fb)$}\label{app:num_sol_F_lambda}
In order to obtain numerical approximations of the pair $(\lambdab,\Fb)$, we get back to the time marching dynamics of the density $\fb(\tb,\zb)$ which satisfies the following equation:
\begin{equation}
   \partial_{\tb} \fb(\tb,\zb) - \cb \partial_{\zb} \fb(\tb,\zb) + \mub(\zb)  \fb(\tb,\zb) = \betab \mathcal{B}( \fb(\tb,\cdot) )(\zb)\,
\end{equation}
The density  $\fb(t,z)$ is expected to behave like $\exp(\lambdab \tb)\Fb(\zb)$ for large time. It is preferable to introduce the frequency of traits in population: $\pb(\tb,\zb) = \fb(\tb,\zb)/\int  \fb(\tb,\zb')\, d\zb'$. The equation for $\pb$ is:
\begin{equation}\label{eq:freq b}
   \partial_{\tb} \pb(\tb,\zb) + \left ( \betab - \bar{\mub}(\tb)\right ) - \cb \partial_{\zb} \pb(\tb,\zb) +   \mub(\zb)  \pb(\tb,\zb)   = \betab \mathcal{B}( \pb(\tb,\cdot) )(\zb)\,,
\end{equation}
where the additional  $\bar{\mub}(t)$ ensures that $\int \pb$ remains constant:
\begin{equation}\label{eq:mubar}
\bar{\mub}(\tb) = \int \mub(\zb') \pb(\tb,\zb')\, d\zb'\, .
\end{equation}
We expect that the pair $(\betab- \bar{\mub}(\tb), \pb)$ does converge to $(\lambdab,\Fb)$ as $\tb\to +\infty$. 


Classical numerical methods were used to approximate \eqref{eq:freq b}-\eqref{eq:mubar} for large time, until some error threshold is reached for $\|\partial_{\tb} \pb(\tb,\cdot)\|_{\infty}$. The transport term $- \cb \partial_{\zb} \pb(\tb,\zb)$ was handled using an upwind scheme. The convolutions involved in operator $\mathcal B$ were handled using the function \texttt{conv} in MATLAB software. The grid mesh was adapted to  the scales in Appendix \ref{app:scaled_pb} in order to capture the appropriate phenomena at the correct scale.

\end{document}